\documentclass[10pt,oneside,a4paper,leqno,final]{article}
\usepackage{latexsym}
\usepackage[latin1]{inputenc}
\usepackage{color}
\usepackage[active]{srcltx}
\usepackage{multirow,rotating}
\usepackage{a4wide}

\usepackage{amscd}
\usepackage{epsfig}
\usepackage{subfigure}
\usepackage{graphicx}
\usepackage{amsmath}
\usepackage{comment}
\usepackage{amssymb}
\usepackage{mathrsfs}
\input xy
\xyoption{all}
\newtheorem{thm}{\sc Theorem}[section]      % numbered within each section
\newtheorem{cor}[thm]{\sc Corollary}        % numbered along with Theorem
\newtheorem{lem}[thm]{\sc Lemma}            % numbered along with Theorem
\newtheorem{prop}[thm]{\sc  Proposition}     % numbered along with Theorem
\newtheorem{defn}[thm]{\sc Definition}      % numbered along with Theorem
\newtheorem{rem}[thm]{\sc Remark}       % numbered along with Theorem
\newtheorem{notation}[thm]{\sc Notation}    % numbered along with Theorem
\usepackage{subfigure}
\newcommand{\RR}{\mathbb{R}}
\newcommand{\CC}{\mathbb{C}}

\newcommand{\from}{\colon}

\newcommand{\R}{\mathbb R}
\newcommand{\Z}{\mathbb Z}

\newcommand{\sgn}{\mathrm{sign}}

\newcommand{\diag}{\mathrm{diag}}
\newcommand{\n}[1]{{\bf #1}}
\newcommand{\planare}{\mathrm{Pl}}

\newcommand{\antiprisma}{\mathrm{Tet}}

\newcommand{\proof}{{\sl Proof.}\hspace{5pt}}   % beginning of proof
\newcommand{\finedim}{\hfill $\Box\\$}            % end of proof

\newcommand{\set}[2]{\left\{{#1}\mid{#2}\right\}}       % the set
\newcommand{\fromto}{\leadsto}
\usepackage{bm}
\newcommand{\simbolovettore}[1]{{\boldsymbol{#1}}}
\newcommand{\va}{\simbolovettore{a}}
\newcommand{\vA}{\simbolovettore{A}}
\newcommand{\vb}{\simbolovettore{b}}
\newcommand{\vB}{\simbolovettore{B}}
\newcommand{\vc}{\simbolovettore{c}}
\newcommand{\vC}{\simbolovettore{C}}
\newcommand{\vd}{\simbolovettore{d}}
\newcommand{\ve}{\simbolovettore{e}}
\newcommand{\vE}{\simbolovettore{E}}

\newcommand{\vF}{\simbolovettore{F}}
\newcommand{\vg}{\simbolovettore{g}}

\newcommand{\vI}{\simbolovettore{I}}
\newcommand{\vj}{\simbolovettore{j}}
\newcommand{\vk}{\simbolovettore{k}}
\newcommand{\vl}{\simbolovettore{l}}
\newcommand{\vm}{\simbolovettore{m}}
\newcommand{\vM}{\simbolovettore{M}}
\newcommand{\vn}{\simbolovettore{n}}
\newcommand{\vN}{\simbolovettore{N}}

\newcommand{\vp}{\simbolovettore{p}}
\newcommand{\vP}{\simbolovettore{P}}
\newcommand{\vq}{\simbolovettore{q}}
\newcommand{\vQ}{\simbolovettore{Q}}
\newcommand{\vr}{\simbolovettore{r}}

\newcommand{\vs}{\simbolovettore{s}}
\newcommand{\vS}{\simbolovettore{S}}
\newcommand{\vt}{\simbolovettore{t}}
\newcommand{\vu}{\simbolovettore{u}}
\newcommand{\vU}{\simbolovettore{U}}
\newcommand{\vv}{\simbolovettore{v}}
\newcommand{\vV}{\simbolovettore{V}}
\newcommand{\vw}{\simbolovettore{w}}
\newcommand{\vW}{\simbolovettore{W}}

\newcommand{\vz}{\simbolovettore{z}}

\newcommand{\vzeta}{\simbolovettore{\zeta}}
\newcommand{\vxi}{\simbolovettore{\xi}}

\newcommand{\zero}{\boldsymbol{0}}

\newcommand{\norm}[1]{\protect\left\protect\Vert\protect#1\protect\right\protect\Vert}

\newcommand{\parabolicmanifold}{{\mathscr P}}

\begin{document}
\pagenumbering{arabic}

\title{%
Dynamics of the the dihedral four-body problem
}

\author{Davide L.~Ferrario and Alessandro Portaluri}

\date{%
\today}
\maketitle

\begin{abstract}
Consider four point particles with equal masses in the euclidean 
space,
subject to the following symmetry constraint: at each instant they
are symmetric with respect to the dihedral group $D_2$, 
that is the group
generated by two rotations of angle $\pi$ around two
orthogonal axes.
Under a 
homogeneous potential of degree $-\alpha$ for $0<\alpha<2$, 
this is a subproblem of the four-body problem, 
% 
% 
%  with three degrees of freedom,
% which is a kind of generalisation of Devaney isosceles problem, 
in
which all orbits have zero angular momentum and the configuration
space is three-dimensional. 
In this paper we
study the flow in McGehee coordinates on the collision manifold,
and 
discuss the qualitative behavior of orbits which reach or come close to a total collision.

\vspace{0.5truecm}

\noindent {\em MSC Subject Class\/}: Primary 70F10; Secondary 37C80.
\end{abstract}

%%%=========================================================================
\section{Introduction}
\label{sec:intro}

The goal of this paper is to investigate the qualitative behavior of
solutions of the dihedral symmetric $4$-body problem in space under
the action of a homogeneous potential of degree $-\alpha$. For the
Newtonian potential this problem is a kind of generalization of
Devaney planar isosceles three body problem \cite{Dev80,Dev81},
following Moeckel's approach to the study of the three body problem
in space \cite{Moe81,Moe83}. The dihedral four-body problem is a
subproblem of the full four-body problem which reduces to a
three-dimensional configuration space. Briefly, one takes
$n= 4$ equal masses whose initial position and velocity are
symmetric with respect to the $2$-dihedral group (which is isomorphic to the Klein group
$\Z_2\times \Z_2$) $D_2\subset SO(3)$. So
the masses form a (possibly planar, degenerate and non-regular) tetrahedron
in space. Because of the symmetry of the problem, the masses will
remain in such a configuration for all time. Hence we have a system
with only three degrees of freedom.

We will use McGehee coordinates 
\cite{McG74} in order to study the dynamics 
of the dihedral four body
problem for a general homogeneous potential of degree $-\alpha$
and with a slight change: we consider McGehee coordinates not only
for studying the behavior of solutions passing close to a total
collision, but also for parabolic orbits connecting central
configurations, projecting the full phase space to a codimension $1$
subspace. We replace the singularity due to total collapse with an
invariant immersed manifold in the full phase space usually called
\emph{total collision manifold} which is the immersion of the
parabolic manifold of the projected phase space. 
% In \cite{FP08} we
% computed all central configurations for this problem and we showed
% that just three types can arise: a planar regular $2l$-gon, a
% regular $l$-gonal prism and a $l$-gonal anti-prism; thus in the four
% body case the planar $2l$-gon and the $l$-gonal prism both reduces
% to a rectangle while the anti-prism reduces to a tetrahedron.
By studying this flow we are able to establish some global
results on the behavior of solutions. We discuss the qualitative
behavior of orbits which reach or come close to the total collision
and the behavior of orbits which start from total collapse,
which implies 
chaotic behavior.  
Such of problem 
includes a some of other subsystems  with one or two 
degrees of freedom studied in the past years. We also observe that
the behavior of the bodies is the same for all values of the
parameters $\alpha \in (0,2)$ and not only for the Newtonian case
$\alpha=1$.

The literature on this problem is quite broad and for this we quote
just few papers in which similar studies were carried over. In
particular some properties of the behavior of the flow on the total
collision manifold has been established for the planar case in
\cite{SL82}. Moreover, in the spatial case in a series of papers
Delgado and Vidal  in \cite{DV99, Vid99} studied the so-called
tetrahedral four body problem without and with rotation. We observe
that the potential in this case (up to a multiplicative constant
depending on the normalization of the masses) coincides with a
subcase of the dihedral four body problem.

The main results are presented in the first part of the article, 
in sections 2, 3 and 4, where we introduce the parabolic flow 
on the collision manifold for the 2-dihedral four-body problem, 
we prove the existence or non-existence of parabolic trajectories
connecting rest points or singular sets, and then analyze 
the consequences on orbits close to total collision via topological
transversality. Many tedious computations are postponed to the appendix.

%%%=========================================================================
\section{Blow up, regularization and total collision manifold}
\label{sec:mcgehhecoord}

We set $\n{n}:=\{1, \dots, n\}$. Let $V := \R^d$ denote the
Euclidean space of dimension $d$ and $n \geq 2$ an integer, and
denote with $\zero$ the origin in $\R^d$. Given $n$ point mass
particles in $V$ $P_k$, $k \in \n{n}$ with positions $\vq_k$,
momenta $\vp_k$ and masses $m_k\in \R^+$, let $\vq, \vp \in V^n$ be
the vectors $(\vq_1, \dots , \vq_n)$, $(\vp_1,\dots \vp_n)$. Let
$\alpha>0$ be a given positive real number. We consider the
potential function (the opposite of the potential energy) defined by
\[
\sum_{i< j} \dfrac{m_im_j}{|\vq_i - \vq_j|^{\alpha}}.
\]
If $M$ denotes the {\em mass matrix\/} namely the diagonal $nd
\times nd$ matrix given by
\[
M:=\diag(\underbrace{m_1, \dots, m_1}_{\textrm{$d$ times}}, \dots,
\overbrace{m_n \dots m_n}_{\textrm{$d$ times}})
\]
then Newton equations
\[
M \ddot \vq =  \dfrac{\partial U}{\partial\vq}
\]
can be written in Hamiltonian form as
\begin{equation}\label{eq:hamilton}
\left\{
\begin{aligned}
M \dot \vq &=  \vp  \\
\dot \vp &= \dfrac{\partial U}{\partial \vq},
\end{aligned}\right.
\end{equation}
where the Hamiltonian (i.e. the energy) is $H =H(\vq,\vp) =
\langle\dfrac{1}{2} M^{-1} \vp,\vp \rangle - U(\vq)$. Here we
denoted by
\[
\dfrac{\partial}{\partial \vq}:=\left(\dfrac{\partial}{\partial
\vq_1}, \dots, \dfrac{\partial}{\partial \vq_n} \right), \qquad
\dfrac{\partial}{\partial \vq_j}:=\left(\dfrac{\partial}{\partial
\vq_j^1}, \dots, \dfrac{\partial}{\partial \vq_j^n} \right).
\]
We will assume that the center of mass and the linear momentum
remains at the origin and we define
\begin{align*}
\vQ=\{\vq=(\vq_1, \dots, \vq_n)\in V^n|  \sum_{i=1}^n m_i
\vq_i=\zero\},\ \ \vP=\{\vp=(\vp_1, \dots, \vp_n)\in V^n|
\sum_{i=1}^n \vp_i=\zero\}.
\end{align*}
For each pair of indexes $i,j \in \n{n}$ let $\Delta_{i,j}$ denote
the collision set of the $i$-th and $j$-th particles $\Delta_{i,j} =
\{\vq \in \vQ|\vq_i=\vq_j\}$. For $i \neq j$, let $\Delta=
\bigcup_{i,j} \Delta_{i,j}$ be the {\em collision set\/} in $\vQ$.
The set of collision-free configurations is denoted by $\hat \vQ =
\vQ \backslash \Delta$. The differential equations
\eqref{eq:hamilton} then determine a vectorfield with singularities
on $V^n\times V^n$, or  a real analytic vectorfield without
singularities on $(V^n\setminus \Delta)\times V^n$. The vectorfield
given by \eqref{eq:hamilton} is everywhere tangent to $\vQ \times
\vP$ and so this $2d(n-1)$ dimensional linear subspace is invariant
under the flow. We henceforth restrict our attention to the flow on
the phase space $\vQ\times \vP$. Consequently, $H$ is an integral of
the system. This means that the level sets
$\Sigma_e:=H^{-1}(e)\cap(\vQ\times \vP)$ are also invariant under
the flow \eqref{eq:hamilton}. We observe that $\Sigma_e$ is a real
analytic submanifold of $(\vQ-\Delta)\times \vP$ having dimension
$2d(n-1)-1$ and the flow is not complete, however. In fact certain
solutions run off in finite time; this happens exactly in
correspondence to certain initial conditions leading to a collision
between two or more particles and the corresponding solution curve
meet $\Delta$ in  finite time. Since we are interested primarily in
solution curves leading to total collapse and since the center of
mass is fixed at the origin, this must occur at the origin of $\vQ$.
%===================================================================================
\subsection{McGehee coordinates}
The aim of this section is to recall McGehee coordinates as given in
\cite{FP08} in order to fix notations. Our basic reference will be
\cite{FP08} and references therein.

The equations \eqref{eq:hamilton} can be written in polar
coordinates by setting the mass norm in $V^n$ defined for every $\vq
\in V^n$ as
\[
\norm{\vq}^2 =  \langle M \vq, \vq \rangle
\]
and suitably rescaling the momentum as follows
\[
\rho =  \norm{\vq},\  \vs = \dfrac{\vq}{\rho},\  \vz  = \rho^{\beta}
\vp \qquad\qquad \text{with $\alpha = 2 \beta$.}
\] In these coordinates
equations~\eqref{eq:hamilton} can be read as
\begin{equation}\label{eq:mcghee1}
\left\{
\begin{aligned}
\rho' &  =  \langle \vz , \vs\rangle \rho \\
\vs' &=   M^{-1}\vz - \langle \vz,\vs\rangle  \vs  \\
\vz' &=
\beta \langle\vz, \vs\rangle \vz +
\dfrac{\partial U}{\partial \vq}(\vs),
\end{aligned}
\right.
\end{equation}
where the time has been rescaled by $dt = \rho^{1+\beta} d\tau$
(that is, $\dfrac{d}{d\tau} = \rho^{1+\beta}\dfrac{d}{dt}$); now the
energy can be written as
\begin{equation}\label{eq:energy}
H = \dfrac{1}{2} \rho^{-\alpha} \langle M^{-1}\vz,\vz\rangle  -
\rho^{-\alpha} U(\vs) = \rho^{-\alpha} \left( \dfrac{1}{2}  \langle
M^{-1}\vz,\vz\rangle  - U(\vs) \right).
\end{equation}

Let $k:=d(n-1)$ and let us consider the projection $(\vq,\vp)
\mapsto (\vs,\vz)$ from the full phase space $\vQ\times  \vP$ to the
reduced space $S^{k-1}\times \vP$ (which is the trivial
$\RR^k$-bundle  on the shape sphere $S^{k-1}$)
\[
\vQ \times \vP \to S^{k-1} \times \vP.
\]
We define the \emph{parabolic manifold} as the projection of all
zero-energy orbits (or, equivalently, of the zero-energy submanifold
of $\vQ\times \vP$) in $S\times \vP$, where $\mathbb S=S^{k-1}$,
that is
\[
 \parabolicmanifold =
\{ (\vs,\vz) \in S\times \vP : \dfrac{1}{2}  \langle
M^{-1}\vz,\vz\rangle  = U(s) \} \subset S^{k-1} \times \vP.
\]
 The next change of coordinates, due to McGehee
\cite{McG74}, is needed for defining the Sundman--Lyapunov
coordinate $v$ and for the regularization of the parabolic manifold
$\parabolicmanifold$. Let $v,\vw$ be defined by
\[
\begin{cases}
v& =\langle \vz, \vs \rangle
\\
\vw & =M^{-1} \vz - \langle \vz, \vs \rangle \vs.
\end{cases}
\]
Then $\vz = vM\vs +M\vw$ and  $\langle \vw, M\vs \rangle = 0$ (so
that $\vw$ is also tangent to $S$), and \eqref{eq:mcghee1} can be
replaced by
\begin{equation}\label{eq:mcghee3}
\left\{
\begin{aligned}
\rho'&= \rho v\\
v' &= \norm{\vw}^2 + \beta v^2 - \alpha U(\vs)\\
\vs' &=  \vw \\
\vw'
&= -\norm{\vw}^2 \vs
+(\beta-1) v\vw
+  M^{-1}\nabla_{\vs} U(\vs),
\end{aligned}
\right.
\end{equation}
where $\nabla_{\vs}$ denotes covariant derivative associated to the
Levi-Civita connection induced by the Riemannian metric
$\langle\cdot, \cdot\rangle_M:=\langle M\cdot, \cdot\rangle$,
i.e.~the component of the gradient tangent to the inertia ellipsoid
$\norm{\vq}=1$, i.e.
\[
\nabla_\vs U =
% \dfrac{\partial U}{\partial \vq}(\vs)   -
%  \langle \dfrac{\partial U}{\partial\vq}(\vs),  \vs \rangle \vs
% =
\dfrac{\partial U}{\partial \vq}(\vs)    + \alpha U(\vs) M \vs.
\]
In fact
\begin{equation*}
\begin{aligned}
\langle M^{-1}\nabla_\vs U(\vs), \vv\rangle &= \langle
\dfrac{\partial U}{\partial \vq}(\vs), \vv\rangle+ \alpha
U(\vs)\langle M\vs, \vv\rangle=\\
&= \langle \dfrac{\partial U}{\partial \vq}(\vs), \vv\rangle=
dU(\vs)[\vv]=\langle M^{-1}\nabla_\vs U(\vs),v\rangle_M
\end{aligned}
\end{equation*}
The energy relation becomes
\[
2\rho^\alpha H=  \|\vw\|^2+\|v\|^2 -2U(\vs)
\]
while the parabolic manifold $\parabolicmanifold$ is then defined by
the equation
\[
v^2 + \norm{\vw}^2   =2 U(\vs).
\]
We observe that differential equation \eqref{eq:mcghee3} gives a
real analytic vectorfield with singularities on the manifold with
boundary $\Omega:=(0, \infty)\times \R \times TS$; equivalently it
defines a real analytic vectorfield without singularities on $[0,
\infty)\times \R \times T(S\setminus \Delta)$.

In these coordinates the assumption about the center of mass gives:
\begin{align}\label{eq:centerofmass}
\vS=\{\vs=(\vs_1, \dots, \vs_n)\in V^n|  \sum_{i=1}^n m_i
\vs_i=\zero\},\ \ \vW=\{\vw=(\vw_1, \dots, \vw_n)\in V^n|
\sum_{i=1}^n m_i \vw_i=\zero\}.
\end{align}
Let $\vM(h)$ denote the subset of $(\rho, \vs,\vw)$ defined as
follows:
\begin{equation}\label{eq:defdiMh}
\vM(h):=\set{(\rho, \vs,\vw) \in [0, \infty)\times \mathbb S\times
\vW}{\eqref{eq:energy}- \eqref{eq:centerofmass}\ \textrm{hold}}.
\end{equation}
Let $\vM_+(h)$ and $\vM_0$ denote the subsets satisfying $\rho>0$
and $\rho=0$ respectively. All three are invariant submanifolds for
the vectorfield \eqref{eq:mcghee1}. Note $\vM_0$ is independent of
$h$. We shall refer to $\vM_0$ as {\em McGehee total collision
manifold.\/} We observe also that the parabolic manifold
$\parabolicmanifold$ introduced above actually is the projection of
McGehee total collision manifold (see
\cite{McG74,Dev80,Moe81,Moe83}); the manifold of $(\vs,\vw)$ here is
not considered as embedded in the space of $(\rho,\vs,\vw)$ with
$\rho=0$.
\begin{rem}
As already remarked by several authors the effect of these
transformations and of the time scaling has the effect of gluing a
boundary given by $\rho =0$ onto the phase space. Since $\rho'=0$
when $\rho=0$, this boundary is invariant under the flow generated
by \eqref{eq:mcghee3}.
\end{rem}
By the second equation in \eqref{eq:mcghee3} can be deduced the
well-known fact that for $0<\alpha <2$,  $v$ is a Lyapunov function
on the flow in the parabolic  manifold, and therefore the flow is
dissipative (gradient-like). Moreover, the {\em equilibrium
points\/} in \eqref{eq:mcghee3} are the projections of the
equilibrium points of \eqref{eq:mcghee1} (and the projection is
one-to-one in the parabolic manifold), which can be found as
solutions of
\begin{equation}\label{eq:centralconf}
\left\{
\begin{aligned}
v^2 & =  2 U(\vs)\\
\nabla_{\vs} U(\vs) & =   \zero \\
\vw  &=   \zero.
\end{aligned}
\right.
\end{equation}
Hence all equilibrium points belong to the parabolic manifold
$\parabolicmanifold$.
\subsection{The linearized flow}
As already observed in section before the vectorfield
\eqref{eq:mcghee3} has the equilibrium points given in
\eqref{eq:centralconf}; moreover as already observed all these
points lie on the total collision manifold. We now turn our
attention to the calculation of the characteristic exponents of the
various equilibrium points in $\vM_0$. Accordingly, let $\bar \vs$
be a central configuration and let $\bar v=\sqrt{2U(\bar\vs)}$. Then
using
\begin{equation}\label{eq:mcghee4}
\left\{
\begin{aligned}
\rho'&= \rho v\\
\vs' &=  \vw \\
\vw' &= -\norm{\vw}^2 \vs +(\beta-1) v\vw +  M^{-1}\nabla_{\vs}
U(\vs),
\end{aligned}
\right.
\end{equation}
one computes the characteristic exponents for the flow in $\vM(h)$.
To study the eigenvalues of the linearized vectorfield at a
restpoint $\bar \gamma =(0,\bar\vs, 0)$, we introduce coordinates
$(\zeta_\rho, \zeta_\vs, \zeta_\vw)$ in the tangent space to
$\R^{2nd+1}$. Let $\vg$ be the Riemannian manifold induced by
$\R^{2nd+1}$ on the  $[2d(n-1)-1]$-dimensional manifold $\vM_0$, $D$
and $\dfrac{D}{dt}$ be respectively the Levi-Civita connection and
the covariant derivative associated to the Riemannian metric $\vg$.
The linearization along an orbit $\gamma_0$ is represented in local
chart by the linear autonomous system
\begin{equation}\label{eq:sistemahamiltoniano}
\frac{D}{dt}\vzeta(t)\, = \, \vC \, \vzeta(t), \qquad \forall \, t
\in\R
\end{equation}
where $\vzeta \in \bar{\gamma}^*\big(T\vM(h)\big)$ is the vector
field along the curve $\bar \gamma$ and where the variational matrix
$A$ is the block matrix represented  by
\begin{equation}
\vC:=\begin{pmatrix}
\bar v & 0 & 0\\
&&\\
0 & 0 & I  \\
&&\\
0 & M^{-1}D\nabla_\vs U(\bar \vs) & (\beta-1)\bar v
\end{pmatrix}.
\end{equation}
Devaney has observed that one can guess eigenvector for problem of
this type. In fact, given the matrix below
\begin{equation}
\vB=\begin{pmatrix}
0 & I\\
&\\ M^{-1}D\nabla_\vs U(\bar \vs) & (\beta-1)\bar v
\end{pmatrix}
\end{equation}
and assuming that $\vxi$ is a $[d(n-1)-1]$-eigenvector of the matrix
$\vA:=M^{-1}D\nabla_\vs U(\bar \vs)$ (the Hessian
$D^2U(\vs)$)associated to the eigenvalue $\lambda$, then the
$[2d(n-1)-1]$-dimensional vector $\tilde \vxi=(\vxi, \mu\vxi)^T$
satisfies the following:
\begin{equation*}
\begin{pmatrix}
0 & I\\
&\\ \vA & (\beta-1)\bar v
\end{pmatrix}\begin{pmatrix}
\vxi \\
&\\ \mu \vxi
\end{pmatrix}=\begin{pmatrix}
\mu \,\vxi\\
&\\ \lambda\,\vxi +(\beta-1)\bar v\mu\,\vxi
\end{pmatrix}
\end{equation*}
By a straightforward calculation it follows that $\mu$ is an
eigenvalue of $\vB$ if
\[
\mu^2+ (1-\beta)\bar v \mu=\lambda,
\]
namely
\[
\mu_{1,2}=\dfrac{(\beta-1)\bar v\pm\sqrt{(\beta-1)^2\bar
v^2+4\lambda}}{2}.
\]
We observe that, in the Newtonian case corresponding to $\beta=1/2$
this formula agrees with that of \cite{DV99}. By this we have the
following result.
\begin{prop}
All equilibrium points are hyperbolic.
\end{prop}
\proof For the proof of this result,  see \cite[Proposition
3.8]{FP08}\finedim
\begin{prop}\label{coro:dimensions}
The dimension of the stable (unstable) manifold of $(\overline
v,\overline\vs,\zero)$ with $\overline v = \sqrt{2 U(\overline\vs)}
> 0$ is 3 (2) if $\overline\vs$ is a rectangle; it is 2
(3) if $\overline\vs$ is a tetrahedron. The dimension of the stable
(unstable) manifold of the point $(- \overline
v,\overline\vs,\zero)$ with $\overline v = \sqrt{2 U(\overline\vs)}
> 0$ is equal to the dimension of the unstable (stable) manifold of
$(\overline v,\overline\vs,\zero)$. The intersection of the stable
(unstable) manifold of $(\overline v,\overline\vs,\zero)$ with the
parabolic manifold $\parabolicmanifold$ has codimension 0 (1) in
$\parabolicmanifold$ if $\overline v >0$. It has codimension 1 (0)
in $\parabolicmanifold$ if $\overline v<0$.
\end{prop}
\proof For the proof of this result, see \cite[Proposition
3.9]{FP08}\finedim

\subsection{McGehee coordinates for the (anisotropic) Kepler problem} Consider
the classical anisotropic 
Kepler problem, in which $S\subset \RR^2$ is a
$1$-dimensional manifold. Let $(0,2\pi) \to S$ be a (maybe partially
defined) local parametrisation, which we will denote by $\theta
\mapsto \vs(\theta)$, where $\vs(\theta) = ( m_1^{-1/2}\cos \theta,
m_2^{-1/2} \sin \theta )$, if $m_1$ and $m_2$ are the diagonal
entries of $M$. If we introduce the vector $\hat\vs = (- m_1^{-1/2}
\sin \theta, m_2^{-1/2} \cos \theta)$ and the scalar $w$ such that
$\vw = w \hat\vs$,  then equations \eqref{eq:mcghee3} turn out to
become
\begin{equation}\label{eq:mcgheeonedim}
\left\{
\begin{aligned}
v' &= w^2 + \beta \left(v^2 - 2 U(\theta)\right)\\
\theta' &=  w \\
w'
&=
(\beta-1) vw
+  U'(\theta),
\end{aligned}
\right.
\end{equation}
where we use the primes for denoting the differentiation with
respect to $\tau$  except in the potential; moreover the parabolic
manifold $\parabolicmanifold$ defined by the equation
\[
v^2 + w^2 = 2U(\theta).
\]
Now consider the flow in part of the parabolic manifold
$\parabolicmanifold$
contained in the half-space $w>0$: by eliminating the term $w^2$
in the equations of $v'$ and $\theta'$, the projection of the
flow on the $(\theta,v)$-plane is contained
in the region $\{ (\theta,v) : v^2  < 2U(\theta) \}$ and
is given by the system
\begin{equation}\label{eq:batta}
\left\{
\begin{aligned}
\theta' &= \sqrt{2U(\theta) - v^2 } \\
v' & = (1-\beta) (2U(\theta) -v^2),
\end{aligned}
\right.
\end{equation}
which can be written also as %(See Figures~\Ref{Fig:Fund4} And~\Ref{Fig:Fund}
%%at page~\pageref{fig:fund4})
\begin{equation}\label{eq:prima}
\dfrac{dv}{d\theta} = (1-\beta) \sqrt{2U(\theta) -v^2}.
\end{equation}
For the projection of the part in $w<0$, the first equation of
\eqref{eq:batta} has to be changed in \( \theta' = -
\sqrt{2U(\theta) - v^2} \). 
If $U(\theta)\equiv U$ is constant, that is if the problem 
is rotationally symmetric,  then
the parabolic manifold $\parabolicmanifold$ is the  cylinder of
equation $v^2+w^2 = 2U$ and the flow in $\parabolicmanifold$ (which
is invariant up to translation in $\theta$) is given by curves
leaving the line $v=-\sqrt{2U}$, $w=0$ of equilibrium points at
$\theta_0$ and reaching the equilibrium line  $v=\sqrt{2U}$, $w=0$
at $\theta_1 = \theta_0 + \dfrac{\pi}{1-\beta}$, since
\[
(1-\beta)(\theta_1 - \theta_0) = \int_{-\sqrt{2U}}^{\sqrt{2U}}
(2U-v^2)^{-1/2} \,dv = \pi.
\]
Thus a ``bouncing'' trajectory on a collision can be seen as  a
solution of the regularized problem only for $1-\beta =
\dfrac{1}{2}$ implies that $\alpha = 1$ (this is the reason the
Levi--Civita regularization might work only in the case $\alpha=1$).
Nonetheless, the Sundman--McGehee regularization can be applied for
every $\alpha$ (yielding a less natural regularization).

%We close this section by showing that the projected differential
%system has some further symmetries. In fact, let $S$ be the map
%defined by
%\[
%S: \R^2\longrightarrow \R^2: \begin{pmatrix}\theta\\
%v\end{pmatrix}\longmapsto \begin{pmatrix}\pi/2-\theta\\
%-v\end{pmatrix}.
%\]

%%%=========================================================================
\subsection{The dihedral $4$-body problem} \label{sec:equivsetup}
Let $\RR^3 \cong \CC \times \RR$ be endowed with coordinates
$(z,y)$, $z\in \CC$, $y\in \RR$. For $l\geq 1$, let $\zeta_l$ denote
the primitive root of unity $\zeta_l = e^{2\pi i/l}$; the
\emph{dihedral} group $D_{l}\subset SO(3)$ is the group of order
$2l$ generated by the rotations
\[
\zeta_l\from (z,y) \mapsto (\zeta_l z, y )
\text{\ and \ }
\kappa\from (z,y) \mapsto (\overline{z},-y),
\]
where $\overline{z}$ is the complex conjugate of $z$. The
non-trivial elements of $D_l = \langle \zeta_l,\kappa \rangle$ are
the  $l-1$ rotations around the $l$-gonal axis $\zeta_l^j$,
$j=1,\ldots,l-1$ and  the $l$ rotations of angle $\pi$ around the
$l$ digonal  axes orthogonal to the $l$-gonal axis %%(see figure
The Newtonian potential for the $n$-body problem, homogeneous with
degree $-\alpha$ induces by restriction on the fixed subspace
$\left(\RR^{6l}\right)^{D_{l}} \cong \RR^3$ a homogeneous potential
defined for each $\vq\in \RR^3$ by
\begin{equation}\label{eq:potenziale}
U(\vq) = \sum_{g\in D_l \smallsetminus \{1\}}
\left| \vq - g \vq \right|^{-\alpha},
\end{equation}
provided we assume (without loss of generality) all masses
$m_i^{2}=1/l$. Now, the potential $U$ in \eqref{eq:potenziale} can
be re-written in terms of coordinates $\vq = (z,y) \in \CC\times
\RR$ as follows.
%\[\begin{aligned}
%U(\vq) &=
%\sum_{j=1}^{l-1} | \vq - \zeta_l^j \vq |^{-\alpha} +
%\sum_{j=1}^{l} | \vq - \zeta_l^j \kappa \vq |^{-\alpha} \\
%&=
%\sum_{j=1}^{l-1} | z - \zeta_l^j z |^{-\alpha} +
%\sum_{j=1}^{l}
%\left(  | z - \zeta_l^j \overline{z}|^2 + 4y^2 \right) ^{-\alpha/2}.
%\end{aligned}
%\]
%By definition, for each $g\in D_{l}$, $U(g\vq)=U(\vq)$.
%Further symmetries of $U$ are:
%\begin{enumerate}
%\item the reflection on the
%plane $y=0$ (given by $h\from (z,y) \mapsto (z,-y)$),
%\item
%the $l$ reflections on the  planes containing the $l$-gonal axis
%and one of the digonal axes,
%\item and the $l$ reflections on the planes
%containing the $l$-gonal axis and the points $( \zeta_l^j e^{\pi i/
%l}  , 0 )$, $j = 1, \ldots, l$.
%\end{enumerate}
%It is not difficult to prove that these are (up to conjugacy and
%multiplication with elements in $D_l$) all the elements of  the
%normaliser of $D_l$ in $O(3)$.

On the unit sphere $S\subset \RR^3$ (of equation $|z|^2 + y^2 = 1$),
parameterized by $(\varphi,\theta)\in (-\pi/2,\pi/2) \times
[0,2\pi)$ with $y=\sin  \varphi$ and $z = \cos\varphi e^{i\theta}$,
the (reduced to the $2$-sphere) potential reads
\begin{equation}\label{potential:sphere}
{U}(\theta,\varphi) \!=\! (2\cos\varphi)^{-\alpha} \left[ 1 +
\dfrac{1}{(\cos^2 \theta + \tan^2 \varphi)^{\alpha/2}} +
\dfrac{1}{(\sin^2 \theta + \tan^2 \varphi)^{\alpha/2}} \right].
\end{equation}
(See \cite{FP08} for the potential in the general case). In
spherical coordinates, the symmetry reflections of $U$ are (up to
conjugacy)
\begin{enumerate}
\item the reflection on the horizontal
plane: $h_\varphi\from (\theta,\varphi)  \mapsto  (\theta,-\varphi)$,
\item
the reflection on the vertical plane $h_\theta\from (\theta,\varphi)
\mapsto (-\theta,\varphi)$
\item and the reflection on the plane
containing the digonal axis and the point $( e^{\pi i/ 2}  , 0 )$,
defined as $h'_\theta\from (\theta,\varphi) \mapsto (\pi/2 - \theta
,\varphi)$.
\end{enumerate}
Thus we can study $U$ only in the left-upper area of the
$D_2$-fundamental domain on $S^2\subset \RR^3$, i.e. in the geodesic
triangle of the shape sphere parameterized by $(\varphi, \theta)\in
[0, \pi/2]\times [0, \pi/2]$. Since in the four body case both the
reflections $h_\varphi, h_\theta$ give arise the same potential we
 only distinguish two cases. \footnote{%
We observe that in the general dihedral $n$-body problem case, the
reflections $h_\varphi$ and $h_\theta$ induce two different
potential functions. In fact in this case we have three different
class of central configurations: $2l$-gon, prism and antiprism type.
For further details, we refer to \cite{FP08} .\/}
\begin{enumerate}
\item (Planar case) $h_\varphi$-reflection (resp.
$h_\theta$-reflection):
\begin{equation}\label{eq:potplanare}
U_\planare(\theta)=
\dfrac{1}{2^\alpha}\left\{1+\dfrac{1}{\sin^\alpha \theta}+
\dfrac{1}{\cos^\alpha\theta}\right\}\quad (\textrm{resp.}\ \
U_\planare(\varphi)=
\dfrac{1}{2^\alpha}\left\{1+\dfrac{1}{\sin^\alpha \varphi}+
\dfrac{1}{\cos^\alpha\varphi}\right\}).
\end{equation}
%\begin{equation}\label{eq:potprisma}
%
%\end{equation}
\item (Tetrahedral case) $(\theta=\pi/4)$
\begin{equation}\label{eq:potantiprisma}
U_\antiprisma(\varphi) = \dfrac{1}{(2\cos\varphi)^\alpha}
+\dfrac{2^{1-\alpha/2}}{(1+\sin^2 \varphi)^{\alpha/2}}.
\end{equation}
\end{enumerate}
%where $c_{l,\alpha}:=\sum_{j=1}^{l-1} \left(  \sin \dfrac{j\pi}{l}
%\right)^{-\alpha} $
\begin{rem}
We observe that the potential in the planar Newtonian  case with
$l=2$ is $1/4$ the potential of the problem studied by the authors
in \cite{SL82}; moreover the potential in the anti-prism  Newtonian
case with $l=2$ is $1/4$ the potential of the problem studied by the
authors in \cite{DV99}. This comes from the fact that we normalized
all masses with the value $1/l$.
\end{rem}
\subsection{Central configurations of the dihedral four body problem with homogeneous potential}
As particular case of the results proven in \cite{FP08}, we recall
that in this case we have $20$ central configurations; in particular
$12$ of all of these are rectangular and the remaining eight are
tetrahedral. %However due to simple analytic form of the potential we
%can easily explicitly compute the central configuration.
\begin{lem}[Planar central configurations] \label{nr:nagon} For any $\alpha \in
(0,2)$ there are exactly
\begin{itemize}
\item $4$ central configurations which are $h_\varphi$-symmetric, and
they are on the vertices $(e^{(2k+1)\pi i/4 },0)$ of the regular
rectangle, for $k \in \n{3}$.
\item $8$ central configurations which are $h_\theta$-symmetric (up to
conjugacy), and they are precisely on the vertices  of a prism:
$(1/\sqrt{2} e^{k\pi i/2 },\pm 1/\sqrt{2})$.
\end{itemize}
\end{lem}
It is left to compute critical points for $\theta = (2k+1)\pi/4$,
that is, to find $h'_\theta$-symmetric central configurations.
\begin{lem}[Tetrahedral central configurations]
\label{nr:antiprism} For any $\alpha \in (0,2)$ there are exactly
$8$ central configurations which are $h'_\theta$ symmetric (up to
conjugacy). They are on the vertices of a prism:
\[
(\cos\hat\varphi e^{(2k+1)\pi i/4 },\pm \sin \hat\varphi).
\]
\end{lem}
Since there are no other central configurations, we can summarize
the results in the following  proposition.
\begin{prop}
\label{propo:main} All central configurations in the dihedral
$4$-body problem are symmetric for one of the three types of
reflections $h_\varphi$, $h_\theta$ (cfr. Lemma \ref{nr:nagon})
 or $h'_\theta$ (cfr. Lemma
\ref{nr:antiprism}).
\end{prop}
\subsection{The linearized flow}
As already observed in section before the vectorfield
\eqref{eq:mcghee3} has the equilibrium points given in
\eqref{eq:centralconf} and all these points lie on the total
collision manifold. Moreover we

We now turn our attention to the calculation of the characteristic
exponents of the various equilibrium points in $\vM_0$. Accordingly,
let $\bar \vs$ be a central configuration and let $\bar
v=\sqrt{2U(\bar\vs)}$. Then using \eqref{eq:mcghee3} one computes
the characteristic exponents for the flow in $\vM(h)$. In fact ,

\begin{table}\centering
\begin{tabular}{|l|c|cc|cc|}
\hline & $\overline v $&
\begin{sideways}$\dim  W^s$\end{sideways} &
\begin{sideways}$\dim  W^u$\end{sideways} &
\begin{sideways}$\dim W^s \cap \parabolicmanifold$\end{sideways}&
\begin{sideways}$\dim W^u \cap \parabolicmanifold$\end{sideways}
\\ \hline
\multirow{2}{*}{rectangular} & $>0$  &
3 & 2 & 3 & 1 \\
 & $<0$
 & 2 & 3 & 1 & 3 \\
\hline \multirow{2}{*}{tethraedral}  & $>0$&
2 &3& 2& 2   \\
 &$<0$ & 3 &2& 2 & 2   \\\hline
\end{tabular}
\caption{Dimensions of stable and unstable manifolds.}
\label{tb:tabella1}
\end{table}
Now we use the above results to describe the set of orbits which
begin or end in quadruple collision. Orbits which begin at collision
are called {\em ejection\/} orbits; orbits which end at triple
collision are called {\em collision\/} orbits; and orbits which do
both are known as {\em ejection-collision\/} orbits. It is
well-known result that any such orbit must be asymptotic to one of
the equilibria associated to a central configuration. That is, they
lie on the stable and unstable manifolds of these equilibria. We
denote by $\vE(c.c)$ and $\vC(c.c)$ the set of ejection and
collision orbits at the central configuration $c.c.$. As a direct
consequences of the dimension of the stable and unstable manifolds,
we have the following:
\begin{prop}
In any energy level $H=h$, both $\vC$ and $\vE$ consist of union of
twenty submanifolds, where twelve (corresponding to the planar
configurations) are bi-dimensional and the others are
three-dimensional (the ones corresponding to tetrahedral central
configurations). All ejection orbits emanate from the equilibria
$v>0$ , whereas all collision orbits are asymptotic to the
equilibria with $v<0$.
\end{prop}

%%===================================================================================
\subsection{Regularization of double collisions}

It is not know if a $\mathscr C^1$-regularization of simultaneous
binary collision is possible in general, but in  this problem, due
to the symmetry, simultaneous double collisions can be regularized
using the type of transformation used by Devaney in \cite{Dev80}.
The aim of this section is to provide the regularization of double
singularities both in the planar and tetrahedral cases.
\subsubsection*{Regularization in the planar case}
We define in the planar and prism case the regularized potential as
follows
\[
W_\planare(\theta):=\sin^\alpha (2\theta)\, U_\planare(\theta)
\]
and we consider the new variables
\[
u:= \dfrac{\sin^\alpha (2\theta)}{\sqrt{W_\planare(\theta)}}\, w,
\qquad \dfrac{d}{d\sigma}:=
\dfrac{\sin^\alpha(2\theta)}{\sqrt{W_\planare(\theta)}}\dfrac{d}{d\tau}.
\]
Then the equation of motions become
\begin{equation}\label{eq:mcgheeonedimplanar}
\left\{
\begin{aligned}
\rho'&=\dfrac{\sin^\alpha(2\theta)\rho
v}{\sqrt{W_\planare(\theta)}} \\
v' &= \dfrac{\sin^\alpha(2\theta)}{\sqrt{W_\planare(\theta)}}
\left[\dfrac{W_\planare(\theta)u^2}{\sin^{2\alpha}(2\theta)} + \beta
\left(v^2 - 2 \dfrac{W_\planare(\theta)}{\sin^\alpha(2\theta)}
\right)\right]\\
\theta' &=  u \\
u' &=\dfrac{\sin^\alpha(2\theta)(\beta-1)v
u}{\sqrt{W_\planare(\theta)}}+
\dfrac{W'_\planare(\theta}{2W_\planare(\theta)}(2\sin^\alpha(2\theta)-u^2)-\dfrac{2\alpha}{\tan(2\theta)}(\sin^\alpha(2\theta)-u^2)
,
\end{aligned}
\right.
\end{equation}
where we still use the primes for denoting differentiation with
respect to $\sigma$ except in $W_\planare'(\theta)$. Moreover in
these new coordinates the energy relation becomes
\begin{equation}\label{eq:energyrelregularizedplanar}
u^2+(v^2-2\rho^\alpha
H)\dfrac{\sin^{2\alpha}(2\theta)}{W_\planare(\theta)}=2
\sin^\alpha(2 \theta),
\end{equation}
and by taking into account of the energy relation, the equation of
motions could be written as follows
\begin{equation}\label{eq:mcgheeonedimplanar2}
\left\{
\begin{aligned}
\rho'&=\dfrac{\sin^\alpha(2\theta)\rho
v}{\sqrt{W_\planare(\theta)}} \\
v' &= \sqrt{W_\planare(\theta)}\left[2(1-\beta)+
\dfrac{(\beta-1)v^2+2\rho^\alpha H}{W_\planare(\theta)}\sin^\alpha(2\theta)\right]\\
\theta' &=  u \\
u' &=\dfrac{\sin^\alpha(2\theta)(\beta-1)v
u}{\sqrt{W_\planare(\theta)}}+
\dfrac{W'_\planare(\theta}{2W_\planare(\theta)}(2\sin^\alpha(2\theta)-u^2)+\dfrac{2\alpha}{\tan(2\theta)}
\left[1+\dfrac{(2\rho^\alpha
H-v^2)\sin^\alpha(2\theta)}{W_\planare(\theta)}\right] ,
\end{aligned}
\right.
\end{equation}
(We recall that the Newtonian case corresponds to $\alpha=1$).
Explicitly
\[
W_\planare(\theta)=\sin^\alpha\theta\cos^\alpha\theta+ \sin^\alpha
\theta + \cos^\alpha\theta, \qquad W'_\planare(\theta)=
\alpha\left(\dfrac{\sin^\alpha(2\theta)}{\tan(2\theta)}+
\dfrac{\sin^\alpha\theta}{\tan\theta}-\cos^\alpha\theta\tan\theta\right).
\]
\begin{rem}
There is a strict analogy between these equation of motions and the
differential system given in \cite[Eqn.(1.12)]{Dev80} and in
\cite[Eqn. (2.3)]{DV99}. We also observe that the above
regularization holds equally well for the parabolic manifold. In
fact, the variable $\rho$ is not essential neither to regularize the
total collision manifold nor to investigate the flow on the total
collision manifold.
\end{rem}
By taking into account the energy relation and by putting $\rho=0$,
it follows
\begin{equation}\label{eq:parabolicregularizedplanar}
u^2+v^2\dfrac{\sin^{2\alpha}(2\theta)}{W_\planare(\theta)}=2
\sin^\alpha(2 \theta),
\end{equation}

%=========================================================================================
\subsubsection*{Regularization in the tetrahedral case}
We define in the tetrahedral case the regularized potential as
follows
\[
W_\antiprisma(\varphi):=(2\cos\varphi)^\alpha\,
U_\antiprisma(\varphi),
\]
and we consider the new variables
\[
u:= \dfrac{(2\cos\varphi)^\alpha}{\sqrt{W_\antiprisma(\varphi)}}\,
w, \qquad \dfrac{d}{d\sigma}:=
\dfrac{(2\cos\varphi)^\alpha}{\sqrt{W_\antiprisma(\varphi)}}\dfrac{d}{d\tau}.
\]
Then the equation of motions become
\begin{equation}\label{eq:mcgheeonedimtetrahedral}
\left\{
\begin{aligned}
\rho'&=\dfrac{(2\cos\varphi)^\alpha\rho
v}{\sqrt{W_\antiprisma(\varphi)}} \\
v' &= \dfrac{(2\cos\varphi)^\alpha}{\sqrt{W_\antiprisma(\varphi)}}
\left[\dfrac{W_\antiprisma(\varphi)u^2}{(2\cos\varphi)^{2\alpha}} +
\beta \left(v^2 - 2
\dfrac{W_\antiprisma(\varphi)}{(2\cos\varphi)^\alpha}
\right)\right]\\
\varphi' &=  u \\
u' &=\dfrac{(2\cos\varphi)^\alpha(\beta-1)v
u}{\sqrt{W_\antiprisma(\varphi)}}+
\dfrac{W'_\antiprisma(\varphi}{2W_\antiprisma(\varphi)}[2(2\cos\varphi)^\alpha-u^2]
+\alpha\tan\varphi[(2\cos\varphi)^\alpha-u^2] ,
\end{aligned}
\right.
\end{equation}
where we still use the primes for denoting differentiation with
respect to $\sigma$ except in $W'_\antiprisma(\varphi)$. Moreover in
these new coordinates the energy relation becomes
\begin{equation}\label{eq:energyrelregularizedtetrahedral}
u^2+(v^2-2\rho^\alpha
H)\dfrac{(2\cos\varphi)^{2\alpha}}{W_\antiprisma(\varphi)}=2
(2\cos\varphi)^\alpha,
\end{equation}
and by taking into account of the energy relation, the equation of
motions could be written as follows
\begin{equation}\label{eq:mcgheeonedimtetrahedral2}
\left\{
\begin{aligned}
\rho'&=\dfrac{(2\cos\varphi)^\alpha\rho
v}{\sqrt{W_\antiprisma(\varphi)}} \\
v' &= \sqrt{W_\antiprisma(\varphi)}\left[2(1-\beta)+
\dfrac{(\beta-1)v^2+2\rho^\alpha H}{W_\antiprisma(\varphi)}(2\cos\varphi)^\alpha\right]\\
\varphi' &=  u \\
u' &=\dfrac{(2\cos\varphi)^\alpha(\beta-1)v
u}{\sqrt{W_\antiprisma(\varphi)}}+
\dfrac{W'_\antiprisma(\varphi}{2W_\antiprisma(\varphi)}[2(2\cos\varphi)^\alpha-u^2]+\\
&-\alpha\tan\varphi(2\cos\varphi)^\alpha\left[1+\dfrac{(2\rho^\alpha
H-v^2))(2\cos\varphi)^\alpha} {W_\antiprisma(\varphi)}\right],
\end{aligned}
\right.
\end{equation}
Explicitly
\[
W_\antiprisma(\varphi)=1+ \dfrac{2^{1+\beta}\cos^\alpha
\varphi}{(1+\sin^2\varphi)^\beta}, \qquad
W'_\antiprisma=-\dfrac{2^{2+\beta}\,\alpha\,\sin\varphi\,\cos^{\alpha-1}\varphi}{(1+\sin^2\varphi)^{\beta+1}}.
\]
\begin{rem}
There is a strict analogy between these equation of motions and the
differential system given in \cite[Eqn.(1.12)]{Dev80} and in
\cite[Eqn. (2.3)]{DV99}. We also observe that the above
regularization holds equally well for the parabolic manifold. In
fact, the variable $\rho$ is not essential neither to regularize the
total collision manifold nor to investigate the flow on the total
collision manifold.
\end{rem}
By taking into account the energy relation
\eqref{eq:energyrelregularizedtetrahedral}, the {\em regularized
total collision manifold\/} reduces to
\begin{equation}\label{eq:energyrelregularizedtetrahedral2}
u^2+v^2\dfrac{(2\cos\varphi)^{2\alpha}}{W_\antiprisma(\varphi)}=2
(2\cos\varphi)^\alpha,\end{equation}

We shall discuss the properties of this flow in the next section.

%================================================================================================================
\section{Colliding and non-colliding parabolic connections}
\label{sec:parabolicconnections}

%\section{Flow on the total collision manifold}\label{sec:flowontotalcollisionmnfld}

The object of this section is to describe in details the flow
generated by the differential system \eqref{eq:mcgheeonedimplanar}
and \eqref{eq:mcgheeonedimtetrahedral2} on the total collision
manifold in the planar and tetrahedral case, respectively. A direct
consequence of this study is the proof of the existence of
connecting orbits on the invariant subset of the
$\parabolicmanifold$ fixed by the the reflections $h_\theta,
h_\varphi, h'_\theta$. The idea to perform our analysis is based
upon a study of the intersection between the stable and unstable
manifold of the equilibria on the parabolic manifold
$\parabolicmanifold$.

\begin{figure}\centering
\includegraphics[width=0.4\textwidth]{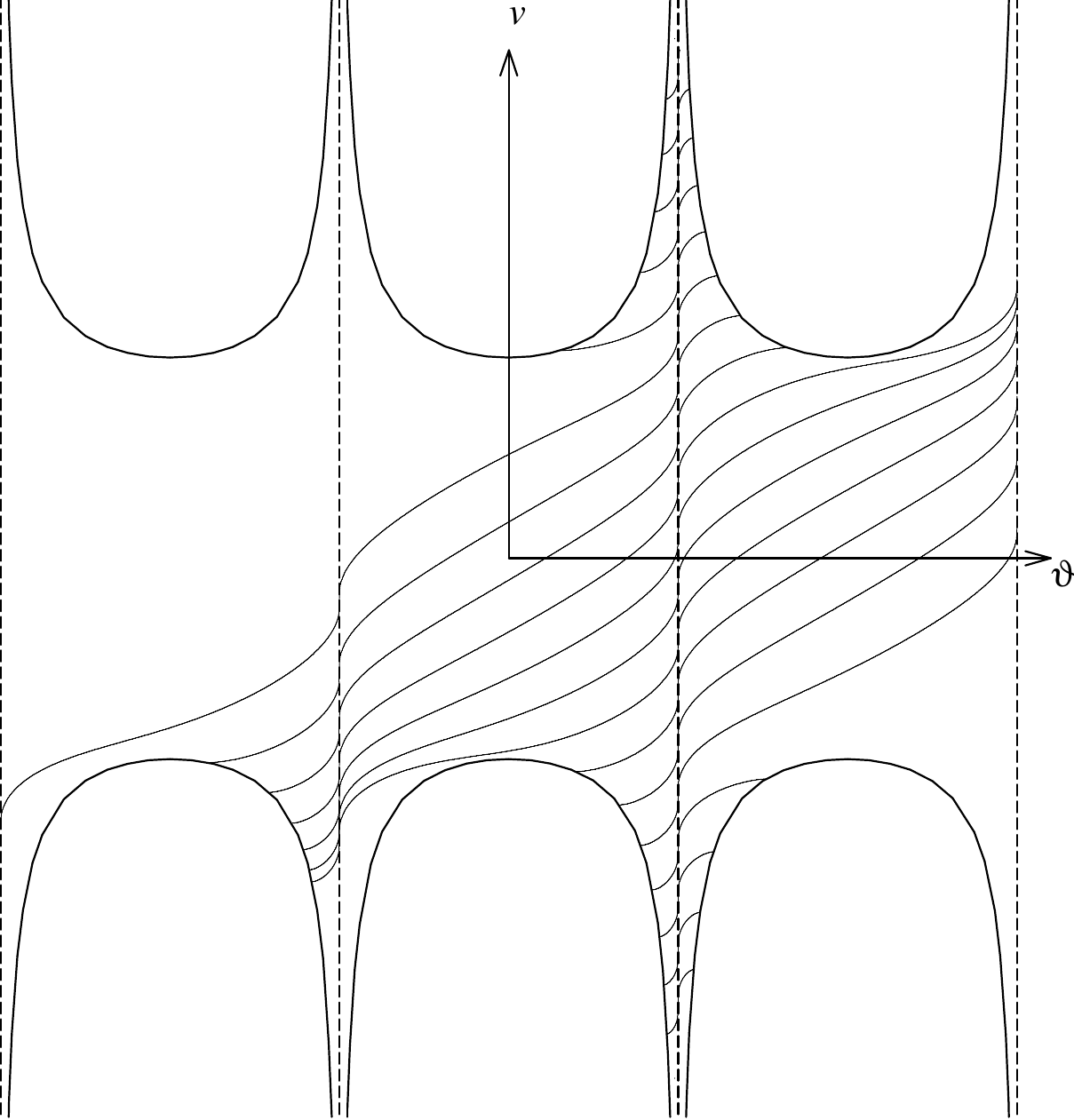}\label{fig:prismsection}
\caption{The flow on the covering of a regularised component of the
parabolic manifold for the planar section.}
\end{figure}

\begin{figure}\centering
\includegraphics[width=0.4\textwidth]{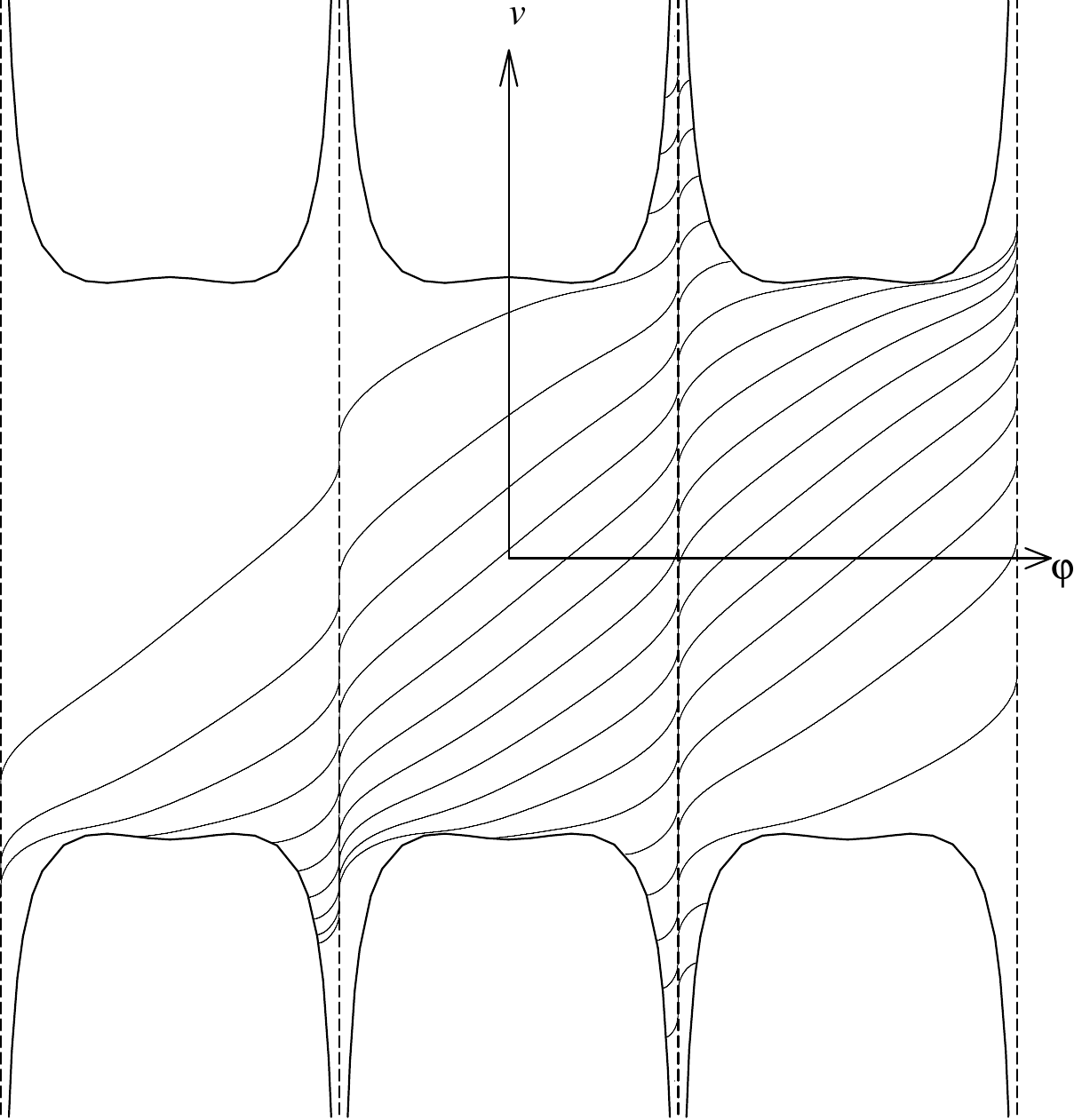}
\caption{The flow on the covering of a regularised component of the
parabolic manifold for the tetrahedral section.}
\label{fig:antiprismsection}
\end{figure}

Due to the gradient-like character of the (regularised) flow on
$\parabolicmanifold$, each orbit  in $\subset
\parabolicmanifold$ either tends to a rest-point or $v \to + \infty$
as $t \to + \infty$. Moreover by taking into account the inequality
$v^2 \leq 2U(\vs),$ it is clear that if $v \to +\infty$ then $U(\vs)
\to + \infty$ and so two or more particles have to collide. Thus if
$\vs$ is parameterized in cartesian coordinates by $(s_1,s_2, s_3)$,
this means that one of the distance $s_i^2+s_j^2\to 0$ for $1\leq
i<j \leq 3$, that is the particles tend to have a binary
configuration. Therefore, by following \cite[Section 6]{Vid99}, we
can define the {\em stable binary escape sets} as follows:
\begin{equation}\label{eq:binescapesetsvplus}
\mathscr B_j^{s,\pm}:= \left\{(v,\vs,\vw) \in \parabolicmanifold
\colon \ v \to +\infty \ \  \textrm{and}  \ \ s_j \to \pm 1\ \
\textrm{for}  \ \ t \to +\infty \right\},
\end{equation}
for $j \in \n{3}$. In an analogous way we define
\begin{equation}\label{eq:binscapesetsvminus}
\mathscr B_j^{u,\pm}:= \left\{(v,\vs,\vw) \in \parabolicmanifold
\colon \ v \to -\infty \ \  \textrm{and}  \ \ s_j \to \pm 1\ \
\textrm{for}  \ \ t \to +\infty \right\}.
\end{equation}
The binary sets together with the $40$ critical points are the only
possible asymptotic {\em alpha\/} and {\em omega limits\/} for an
orbit on $\parabolicmanifold$.
\begin{notation}
If $X$ and $Y$ are two equilibrium points or even two binary escape
sets on $\parabolicmanifold$ we shall write $X \fromto Y$ in order
to indicate the existence of an orbit in $\parabolicmanifold$ with
asymptotic behavior $X$ for $t \to -\infty$ and $Y$ as $t \to +
\infty$ i.e. if there exists an orbit on the total collision
manifold which lies in $W^u(X) \cap W^s(Y)$. Since at each central
configuration corresponds two value of the function $v$, we
introduce the following convention. We denote with superscript $+$
(resp. $-$), central configurations corresponding to the positive
(resp. negative) value of the coordinate $v$.
\end{notation}
\begin{rem}
In order to simplify the study of the connecting orbits between
critical points of the function $v$  some remarks are in order.
\begin{enumerate}
\item At first we observe that the presence of the symmetry of the
vector field \eqref{eq:mcghee3}, namely
\[
(v,\vs,\vw) \mapsto (-v, \vs, -\vw)
\]
implies that for each connection $X^+ \fromto Y^+$ there is a
symmetric one $Y^- \fromto X^-$ which will be term \emph{dual
connection}. Moreover this transformation sends the stable and
unstable manifold of the point $X^+$ respectively on the unstable
and stable manifold of the point $X^-$.

In conclusion the transformation $(\vs, \vw)\mapsto (\vs, -\vw)$,
together with a time reversal $t\mapsto -t$, takes orbits into
orbits. This transformation carries $\mathscr B_j^{s, \pm}$ to
$\mathscr B_j^{u, \pm}$ and viceversa.
\item
The second observation is that the increasing character of the
function $v$ imposes some restriction about the existence of
connecting orbits between central configurations. In fact, since $v$
is non-decreasing and strictly increasing away from the central
configurations, connecting orbits between $X$ and $Y$ cannot occur
if $v(Y) \leq v(X)$. Denoting by $\planare, \antiprisma$
respectively the planar and tetrahedral type central configurations,
as already observed,  we have
\[
v(\planare) \geq v(\antiprisma)>0.
\]
In fact, in the Newtonian case ($\alpha=1$) the inequality above,
readily follows since
\[
v^+(\planare):=\sqrt{1+2\sqrt2}>\dfrac{\sqrt{6\sqrt6}}{2}:=v^+(\antiprisma).
\]
However in the homogeneous case, it is enough to observe that the
function $\vl: [0,2)\to \R$ defined by
\[
\vl(\alpha):=v^+(\planare)-v^+(\antiprisma)=\sqrt{2^{1-\alpha}+2^{2-\alpha/2}}-\sqrt{2^{1-\alpha}3^\alpha6^{-\alpha/2}+
4^{1-\alpha}3^{\alpha/2}2^{\alpha/2}},
\]
is a non-negative increasing function such that $\vl(0)=0$; thus it
is strictly positive for $\alpha\neq0$.
\end{enumerate}
\end{rem}
In order to study the possible connections in the dihedral four body
problem for any homogenous potential of degree $-\alpha$, for
$\alpha \in (0,2)$ we shall analyze two cases that it contains as
subproblems.
\begin{itemize}
\item {\em The Planar problem.\/} This particular case appears when
we consider the invariant set $\{\varphi'=\varphi=0\}$ (or
$\{\theta= 0, \theta'=0\}$) which corresponds to the subset
$\parabolicmanifold^{\planare}$ of the parabolic manifold fixed by
the symmetry $h_\varphi$. In this case the parabolic manifold is
two-dimensional and in the Newtonian case $\alpha=1$ it was studied
by \cite{SL82}; in Figure \ref{fig:prismsection} it is given a
qualitative description of the flow on the covering of the
regularized component.
\item {The Tetrahedral problem.\/} This particular case appears when
we consider the invariant set $\{\theta=\pi/4, \theta'=0\}$. Also in
this case the parabolic manifold is  also two-dimensional. This
problem in the Newtonian case was studied by Delgado \& Vidal in
\cite{DV99}, and the flow on the covering of the regularized
component is qualitative depicted in Figure
\ref{fig:antiprismsection}.
\end{itemize}
The main goals of this section is to completely study the flow in
the above two cases. In Theorem \ref{thm:diedraleplanare} are
described all connecting orbits corresponding to the planar section
and in  Theorem \ref{thm:diedraletetracolliding} and Theorem
\ref{thm:diedraletetranoncolliding} are described all
connecting orbits corresponding to the tetrahedral section.\\
Before proceeding further with the description of the flow on the
parabolic manifold some comments are in order. First of all in
figure \ref{fig:prismsection} it is represented the flow obtained by
integrating the corresponding $1$-dimensional Cauchy problem for
different values of the initial conditions. This flow is represented
on the covering over each piece of $1$-sphere between the binary
collision and the binary collision. 
Figure \ref{fig:prismsection}
actually represents the flow on the regularized component of the
parabolic manifold over the $1$-sphere corresponding to the planar
section. Figure
\ref{fig:antiprismsection} represents the flow on the covering
of a regularized component of the parabolic manifold for the
tetrahedral section. In this case the situation is more complicated
due to the presence of two different types of central
configurations. By changing coordinates, locally in the neighborhood
of binary collisions, the regularized component of the parabolic
manifold in the planar section resemble that studied by
McGehee in \cite{McG74} for the collinear three body problem, while
the regularized component of the parabolic manifold in the
tetrahedral resemble that studied by Devaney in \cite{Dev80} in the
isosceles three body problem.

%%%%%%===================================================================================================

\subsection{Flow on the parabolic equation: the planar case}\label{subsec:planarflow}

The projection of the equation of motions on the parabolic manifold
are
\begin{equation}\label{eq:mcgheeonedimplanar3}
\left\{
\begin{aligned}
v' &= (1-\beta)\sqrt{W_\planare(\theta)}\left[2-
\dfrac{v^2\sin^\alpha(2\theta)}{W_\planare(\theta)}\right]\\
\theta' &=  u \\
u' &=\dfrac{\sin^\alpha(2\theta)(\beta-1)v
u}{\sqrt{W_\planare(\theta)}}+
\dfrac{W'_\planare(\theta}{2W_\planare(\theta)}(2\sin^\alpha(2\theta)-u^2)+\dfrac{2\alpha}{\tan(2\theta)}
\left[1-\dfrac{v^2\sin^\alpha(2\theta)}{W_\planare(\theta)}\right] ,
\end{aligned}
\right.
\end{equation}
where the regularized parabolic manifold $\parabolicmanifold$ is
given by
\begin{equation}
u^2+v^2\dfrac{\sin^{2\alpha}(2\theta)}{W_\planare(\theta)}=2
\sin^\alpha(2 \theta),
\end{equation}
Now consider the flow in part of the parabolic manifold
$\parabolicmanifold$ contained in the half-space $u>0$: the
projection of the flow on the $(\theta,v)$-plane is contained in the
region $\{ (\theta,v) : v^2  < 2U_\planare(\theta) \}$ and is given
by the system
\begin{equation}\label{eq:battaplanare}
\left\{
\begin{aligned}
\theta' &= \sqrt{2U_\planare(\theta) - v^2 } \\
v' & = (1-\beta) (2U_\planare(\theta) -v^2),
\end{aligned}
\right.
\end{equation}
which can be written also as
\begin{equation}\label{eq:primaplanare}
\dfrac{dv}{d\theta} = (1-\beta) \sqrt{2U_\planare(\theta) -v^2}.
\end{equation}
For the projection of the part in $u<0$, the first equation of
\eqref{eq:batta} has to be changed in \( \theta' = -
\sqrt{2U(\theta) - v^2} \).
\begin{notation}
In order to give the complete description of the flow on the
parabolic manifold $\parabolicmanifold$ including orbits that {\em
run out along the arms of binary collisions,\/} we use the following
convention. Let $W^s(\pi/2)$ ($W^u(\pi/2)$) denote the set of points
on $\parabolicmanifold$ whose forward (backward) orbit run out along
the upper (lower) arm at $\theta=\pi/2$. Similar definitions for
$W^s(-\pi/2)$ ($W^u(-\pi/2)$ and $W^s(0)$ ($W^u(0)$).
\end{notation}
Before stating the main result of this section, we remark that the
projected  differential system \eqref{eq:battaapr} has the following
symmetries:
\begin{equation}\label{eq:simmetriepl}
\begin{aligned}
\mathscr R_1&: (\theta,v,u,\sigma)\longmapsto(\theta, -v,-u,-\sigma)\\
\mathscr R_2&:(\theta,
v,u,\sigma)\longmapsto(\pi/2-\theta,v,-u,\sigma)
\end{aligned}
\end{equation}
and the composition of both.

\begin{prop}\label{prop:flussoplanare}
For any $\alpha \in (0, 1]$, the global flow in the parabolic
manifold for the dihedral planar four body problem is described by
the following relationship among the stable and unstable manifolds
on $\parabolicmanifold$:
\begin{enumerate}
\item $W^u(\vc_{-})\subset W^s(0)\cup W^s(\pi/2)$.
\item $W^u(\epsilon\pi/2)\cap W^s(\sigma \pi/2)$ is an open set
for $\epsilon, \sigma \in \{0,1\}$,
\end{enumerate}
where we denoted by $\vc_{-} $ the planar central configuration
corresponding to the negative value of $v$.
\end{prop}
\proof The proof of this proposition is deferred to the Appendix
\ref{app:planar}. \finedim
\begin{rem}
Proposition \ref{prop:flussoplanare} implies that there are no
saddle connections on $\parabolicmanifold$, thus preventing
regularizability of the total collision.

In our case item 1 in Proposition \ref{prop:flussoplanare} says that
orbits passing close to the planar homothetic orbit escape from a
neighborhood of total collapse with successive binary collisions;
item 2 says that any combination of orbits running from one arm of
binary collision to another can occur, and there is an open set of
initial conditions in $\parabolicmanifold$ whose orbit have this
property. In particular no saddle connections on
$\parabolicmanifold$, thus preventing regularizability of the total
collision.
\end{rem}
\begin{prop}\label{prop:flussoplanareb}
For any $\alpha \in (1, 2)$, the global flow in the parabolic
manifold for the dihedral planar four body problem is described by
the following relationship among the stable and unstable manifolds
on $\parabolicmanifold$. There exist $\alpha_{0*}, \alpha_* \in
(1,2)$ with $1<\alpha_{0*}< \alpha_*<2$, such that
\begin{enumerate}
\item $W^u(\vc_{-})\subset W^s(0)\cup W^s(\pi/2)$, for any $\alpha \in (1,2)\setminus \{\alpha_{0*},
\alpha_*\}$.
\item For $\alpha=\alpha_{0*}$ the right-hand branch of $W^u(\vc_-)$
coincides with the right-hand one of  $W^s(\vc_+)$.
\item For $\alpha=\alpha_{*}$ the right-hand branch of $W^u(\vc_-)$
coincides with the left-hand one of  $W^s(\vc_+)$.
\item $W^u(\epsilon\pi/2)\cap W^s(\sigma \pi/2)$ is an open set for
$\epsilon, \sigma \in \{0,1\}$,
\end{enumerate}
where as before we denoted by $\vc_{-} (resp. \vc_{+})$ the planar
central configuration corresponding to the negative (resp. positive)
value of $v$.
\end{prop}
\proof The proof of this proposition is deferred to the Appendix
\ref{app:planar}.\finedim

%============================================================================================================

\subsection{Flow on the parabolic equation: the  tetrahedral case}\label{subsec:tetrahedralflow}

The projection of the equation of motions on the parabolic manifold
and in the tetrahedral case are given by
\begin{equation}\label{eq:mcgheeonedimtetrahedral3}
\left\{
\begin{aligned}
v' &= (1-\beta)\sqrt{W_\antiprisma(\varphi)}\left[2-
\dfrac{v^2(2\cos\varphi)^\alpha}{W_\antiprisma(\varphi)}\right]\\
\varphi' &=  u \\
u' &=\dfrac{(2\cos\varphi)^\alpha(\beta-1)v
u}{\sqrt{W_\antiprisma(\varphi)}}+
\dfrac{W'_\antiprisma(\varphi}{2W_\antiprisma(\varphi)}[2(2\cos\varphi)^\alpha-u^2]
-\alpha\tan\varphi(2\cos\varphi)^\alpha\left[1-\dfrac{v^2(2\cos\varphi)^\alpha}
{W_\antiprisma(\varphi)}\right],
\end{aligned}
\right.
\end{equation}
and the regularized parabolic manifold takes the form
\begin{equation}
u^2+v^2\dfrac{(2\cos\varphi)^{2\alpha}}{W_\antiprisma(\varphi)}=2
(2\cos\varphi)^\alpha.
\end{equation}

Now consider the flow in part of the parabolic manifold
$\parabolicmanifold$ contained in the half-space $u>0$: the
projection of the flow on the $(\theta,v)$-plane is contained in the
region $\{ (\theta,v) : v^2  < 2U_\planare(\theta) \}$ and is given
by the system
\begin{equation}\label{eq:battaapr}
\left\{
\begin{aligned}
\varphi' &= \sqrt{2U_\antiprisma(\varphi) - v^2 } \\
v' & = (1-\beta) (2U_\antiprisma(\varphi) -v^2),
\end{aligned}
\right.
\end{equation}
which can be written also as %(See Figures~\Ref{Fig:Fund4} And~\Ref{Fig:Fund}
%%at page~\pageref{fig:fund4})
\begin{equation}\label{eq:primaapr}
\dfrac{dv}{d\varphi} = (1-\beta) \sqrt{2U_\antiprisma(\varphi)
-v^2}.
\end{equation}
For the projection of the part in $u<0$, the first equation of
\eqref{eq:batta} has to be changed in \( \varphi' = -
\sqrt{2U_\antiprisma(\varphi) - v^2} \).

The differential system \eqref{eq:battaapr} has the following
symmetries:
\begin{equation}\label{eq:simmetrieapr}
\begin{aligned}
\mathscr S_1&: (\varphi,
v,u,\sigma)\longmapsto(\varphi,-v,-u,-\sigma)\\
\mathscr S_2&:(\varphi, v,u,\sigma)\longmapsto(-\varphi,v,-u,\sigma)
\end{aligned}
\end{equation}
and the composition of both. The first symmetry comes from the
reversibility while the second comes from the reflection with
respect to the horizontal plane.
\begin{prop}\label{prop:flussotetra}
The global flow in the parabolic manifold for the dihedral planar
four body problem is described by the following relationship among
the stable and unstable manifolds on $\parabolicmanifold$:
\begin{enumerate}
\item $W^u(\vc_{-})\subset W^s(-\pi/2)\cup W^s(\pi/2)$.
\item $W^u(\epsilon\pi/2)\cap W^s(\sigma \pi/2)$ is an open set
for $\epsilon, \sigma \in \{0,1\}$.
\item $W^s(\vr_{-})\subset W^u(\vc_{-})\cup W^u(\pi/2)$
\end{enumerate}
where we denoted by $\vc_{-} $ the planar central configuration
corresponding to the negative value of $v$ and $\vr_{-}$ is the
tetrahedral central configuration.
\end{prop}
\proof The proof of this proposition is deferred to the Appendix
\ref{app:tetra}.\finedim

\begin{rem}
In our case item 1 in Proposition \ref{prop:flussotetra} says that
orbits passing close to the planar homothetic orbit escape from a
neighborhood of total collapse with successive binary collisions;
item 2 says that any combination of orbits running from one arm of
binary collision to another can occur, and there is an open set of
initial conditions in $\parabolicmanifold$ whose orbit have this
property. The last item says that the stable branch of $\vr_{-}$ die
in $\vc_{-}$ and in the the lower arm of a double collision at
$\theta=\pi/2$.
\end{rem}
\begin{prop}\label{prop:flussotetrab}
For any $\alpha \in (1, 2)$, the global flow in the parabolic
manifold for the dihedral planar four body problem is described by
the following relationship among the stable and unstable manifolds
on $\parabolicmanifold$. There exist $\alpha_{0*}, \alpha_* \in
(1,2)$ with $1<\alpha_{0*}< \alpha_*<2$, such that
\begin{enumerate}
\item $W^u(\vc_{-})\subset W^s(-\pi/2)\cup W^s(\pi/2)$, for any $\alpha \in (1,2)\setminus \{\alpha_{0*},
\alpha_*\}$.
\item For $\alpha=\alpha_{0*}$ the right-hand branch of $W^u(\vr_-)$
coincides with the right-hand one of  $W^s(\vl_+)$.
\item For $\alpha=\alpha_{*}$ the right-hand branch of $W^u(\vr_-)$
coincides with the left-hand one of  $W^s(\vr_+)$.
\item $W^u(\epsilon\pi/2)\cap W^s(\sigma \pi/2)$ is an open set for
$\epsilon, \sigma \in \{-1,1\}$,
\end{enumerate}
where as before we denoted by $\vc_{-} (resp. \vc_{+})$ the planar
central configuration corresponding to the negative (resp. positive)
value of $v$, $\vr_{-} (resp. \vr_{+})$ the tetrahedral central
configuration corresponding to the negative (resp. positive) value
of $v$ for $\varphi>0$ and finally $\vl_{-} (resp. \vl_{+})$ the
tetrahedral central configuration corresponding to the negative
(resp. positive) value of $v$ for $\varphi<0$.
\end{prop}
\proof The proof of this proposition is deferred to the Appendix
\ref{app:tetra}.\finedim
%flowoncollision_2
\begin{notation}
We denote the central configurations in the first octant of the
shape sphere corresponding to the fundamental domain by:
\begin{enumerate}
\item (tetrahedral c.c.): $\ve_{k1}:=(k\pi/4, \theta_c), \ve_{k2}:=(k\pi/4, -\theta_c)$
\item (planar c.c.): $\vp_{1k}:=(k\pi/4,0)$,  $\vp_{2k}:=(\pi/2, \pi/4)$, $\vp_{3k}:=(0,
\pi/4)$
\end{enumerate}
where $\theta_c:=\arctan(1/\sqrt2)$ and $k \in \n{3}$.
\end{notation}
\begin{thm}\label{thm:diedraleplanare} There exist $\alpha_{0,*} , \alpha_* \in
(1,2)$ with $\alpha_{0*}<\alpha_*$ such that on the parabolic
manifold $\parabolicmanifold$ the following connections and their
dual ones occur:
\begin{enumerate}
\item for $\alpha \in (0,\alpha)\setminus\{\alpha_{0*}\} $
\begin{itemize}
\item[(P1)] $\vp_{11}^+ \fromto \mathscr B^{s,+}_1;$
\item[(P2)] $\vp_{11}^+ \fromto \mathscr B^{s,+}_2;$
\item[(P3)] $ \mathscr B^{u,+}_1\fromto \vp_{11}^+;$
\item[(P4)] $ \mathscr B^{u,+}_2\fromto \vp_{11}^+.$
\item[(P5)] There is no saddle connections between two planar central configurations
$\planare^{\pm}.$
\end{itemize}
Moreover in the case (P4) and (P5) the last binary collision occur
for a positive value of $v$.
\item for $\alpha=\alpha_{0*}$ there are saddle connections $\vp_{11}^- \fromto
\vp_{11}^+$ and $\vp_{11}^- \fromto \vp_{12}^+$ passing through one
binary collision.
\item For $\alpha=\alpha_{*}$ there are   saddle connections $\vp_{11}^- \fromto
\vp_{11}^+$, $\vp_{11}^- \fromto \vp_{12}^+$ and finally $\vp_{11}^-
\fromto \vp_{13}^+$ passing through two binary collisions.
\item For $\alpha\in (\alpha_{*},2)$
\begin{itemize}
\item[(P3)] $ \mathscr B^{u,+}_1\fromto \vp_{11}^+;$
\item[(P4)] $ \mathscr B^{u,+}_2\fromto \vp_{11}^+.$
\end{itemize}
However in this last two cases the last binary collision occur for a
negative value of $v$.
\end{enumerate}
%Let $P_1$ be the planar central configuration in the fundamental
%domain and for $j=0,1$ let $B^{s}_j$ be the stable manifold of the
%binary collision sets. Then the following connections and the dual
%ones occur in $\parabolicmanifold$.
\end{thm}
\begin{rem}
We observe that analogous connections and saddle connections occur
between the remaining  central configurations and (or) central
configurations/binary escape sets.
\end{rem}
\begin{proof}
The proof of this result immediately follows by the discussion at
the beginning of this section as well as Proposition
\ref{prop:flussoplanare} and Proposition \ref{prop:flussoplanareb}.
\end{proof}
\begin{thm}\label{thm:diedraletetracolliding}
There exist $\gamma_{0,*} , \gamma_* \in (1,2)$ with
$\gamma_{0*}<\gamma_*$ such that on the parabolic manifold
$\parabolicmanifold$ the following connections and their dual ones
occur:
\begin{enumerate}
\item for $\gamma \in (0,\gamma_*)\setminus\{\gamma_{0*}\} $
\begin{itemize}
\item[(T1)] $\ve_{11}^+ \fromto  \vp_{11}^+, \ve_{12}^+ \fromto \vp_{11}^+
$ (non-colliding parabolic orbit);
\item[(T2)] $\ve_{11}^+ \fromto \mathscr B_3^{s,+},\ve_{12}^+ \fromto \mathscr B_3^{s,-};$
\item[(T3)] $ \mathscr B^{u,-}_3\fromto \ve_{11}^+,  \mathscr B^{u,-}_3\fromto \ve_{12}^+;$
\item[(T4)] $ \mathscr B^{u,-}_3\fromto \ve_{11}^+.$
\item[(T5)] There is no saddle connections between two central
configurations.
\end{itemize}
\item for $\gamma=\gamma_{0*}$ there are saddle connections $\ve_{11}^- \fromto
\ve_{12}^+$, $\ve_{11}^- \fromto \ve_{21}^+$ and finally $\ve_{11}^-
\fromto \ve_{31}^+$ passing through two binary collision.
\item For $\gamma=\gamma_{*}$ there  are  saddle connections $\ve_{11}^- \fromto
\ve_{11}^+$ and $\ve_{11}^- \fromto \ve_{21}^+$ passing through one
binary collision.
\item For $\alpha\in (\gamma_{*},2)$
\begin{itemize}
\item[(T6)] $ \mathscr B^{u,-}_3\fromto \ve_{11}^+,  \mathscr B^{u,-}_3\fromto \ve_{12}^+;$
\item[(T7)] $ \mathscr B^{u,-}_3\fromto \ve_{11}^+.$
\end{itemize}
However in this last two cases the last binary collision occur for a
negative value of $v$.
\end{enumerate}
\end{thm}
\proof The proof of this result immediately follows by the
discussion at the beginning of this section as well as Proposition
\ref{prop:flussotetra} and Proposition \ref{prop:flussotetrab}.
\finedim
\begin{rem}
Analogous connections holds between the remaining six tetrahedral
type central configurations and the remaining three planar central
configurations located on the equator $\varphi=0$ of the shape
sphere.
\end{rem}

%%%=========================================================================
\subsection{Non-colliding parabolic connections }
\label{sec:noncollidingparabolic}

The scope of this section is to prove the existence of non-colliding
parabolic connections. This will be achieved by using topological
arguments similar to those used in \cite{Moe83}. Before, we recall
some basic definitions about {\em topological transversality\/}.

Given two submanifolds $ \vP, \vQ$ of the manifold $ \vN$, we say
that their are {\em transverse at $\vp$\/} if $\vp \in  \vP \cap
\vQ$ and if $T_\vp  \vP + T_\vp  \vQ = T  \vN$. Of course by the
rank theorem it follows that the sum of the dimension of $ \vP$ and
of $ \vQ$ can be greater than or equal to the dimension of the
manifold $ \vN$. Moreover we say that $\vP$ meets $ \vQ$ {\em
transversely\/} if either $ \vP \cap \vQ=\emptyset$ or else $ \vP$
is transverse to $ \vQ$ at each $\vp \in  \vP \cap  \vQ$. The
important fact is that the transverse intersection cannot be removed
by a $\mathscr C^1$- small change; thus in other words, two
submanifolds which meet transversely do so stably. This statement
can be easily proven since the transversality is an open condition.

Let $ \vP$ be a $\vp$-dimensional submanifold of the
$\vn$-dimensional manifold $ \vN$ and let us denote by $D^k$ the
$k$-dimensional disk. Given a continuous map of pairs $\sigma \colon
(D^{n-p}, D^{n-p}\backslash \{0\}) \to ( \vN , \vN \backslash \vP)$
such that $\sigma(0) =\vp \in \vP$, we will say that it is a
$(n-p)$-{\em complementary cell\/} at $\vp$ to $ \vP$ if in a small
enough neighborhood $ U$ of the point $\vp$ we have $H_{n-p}\big(
\vU, \vU \backslash ( \vU \cap  \vP)\big) \cong \Z$ and $\sigma$ is
a generator of the $(n-p)$-dimensional relative homology group. Now
we are in position to state the following useful definition.
\begin{defn}
Let $ \vP,\vQ$ be a $\vp$ and $\vq$-dimensional submanifold
respectively of the $n$-dimensional manifold $\vN$.  We say that
they are {\em topologically transverse\/} at the point $\vp \in  \vP
\cap  \vQ$ and we write $ \vP \pitchfork \vQ$ if there exists a
$(n-p)$-complementary cell at $\vp$ to $ \vP$ in $( \vQ,  \vQ
\backslash \vQ \cap \vP)$ and a $(n-q)$-complementary cell at $\vp$
to $ \vQ$ in $(\vP, \vP \backslash  \vQ \cap \vP)$.
\end{defn}
\begin{thm}\label{thm:diedraletetranoncolliding}
The following non-colliding  connections and their dual ones occur
on the parabolic manifold $\parabolicmanifold$:
\begin{enumerate}
\item[(Q1)] $\ve_{11}^+\fromto \vp_{11}^+, \ve_{11}^+\fromto\vp_{21}^+,\ve_{11}^+\fromto \vp_{31}^+$;
\item[(Q2)] $\ve_{11}^+\fromto \mathscr B_{1}^{s,+}, \ve_{11}^+\fromto \mathscr  B_{2}^{s,+},\ve_{11}^+\fromto
\mathscr B_{3}^{s,+}$.
\end{enumerate}
Furthermore for each of these connections the intersection between
the stable and unstable manifold is topologically transverse.
\end{thm}
\begin{rem}
We observe that analogous connections occur between the remaining
central configurations. Item (Q1) was already proved in Theorem
\ref{thm:diedraletetracolliding} Item (T1) in complete different
ways.
\end{rem}
\proof  We shall denote by the symbol $\parabolicmanifold_\lambda$
the level set of  $v$, namely
\[
\parabolicmanifold_\lambda\, :=\, \{(v,\vs,\vw) \in \parabolicmanifold; \ \
v \geq \lambda\}.
\]
Clearly since $v$ is gradient-like,  it follows that for all
$\lambda \in \R$ each of these sets are invariant for the flow
induced by the function $v$.

Let $\vp \colon \parabolicmanifold \longrightarrow \mathscr D$
(defined by $(v,\vs,\vw) \mapsto \vs$) be the canonical projection
onto the fundamental domain of the action of $D_2$ on the shape
sphere and we let $S_\lambda:= \vp(\parabolicmanifold_\lambda)$.
Taking into account the inequality $v^2 \leq 2 U(\vs)$ it follows
that
\[
\mathscr D_\lambda = \, \{\vs \in  \mathscr D\colon U(\vs) \geq
\lambda^2/2\}.
\]
Denoting by $c:=v(\ve_{11}^+)$ the level set of the function $v$ at
the tetrahedral type central configuration then for all positive
$\varepsilon$ the set given by:
\[
\overline {\mathscr D}_{c+\varepsilon}  \, := \, \big\{s \in
\mathscr D\colon \ \ U(\vs) \geq (c+ \varepsilon)^2/2\big\}
\]
is just the fundamental domain  with a small disks about $\ve_{11}$
removed. The proof of this result is completed by using the
following lemma.
\begin{lem}\label{thm:lemmettinofurbo}
Let $\alpha \colon S^1 \to \parabolicmanifold_{c+\varepsilon}$ be an
analytic closed curve such that $\vp \circ \alpha \colon S^1 \to
\mathscr D_{c+\varepsilon}$ is not contractible in
$\overline{\mathscr D}_{c+\varepsilon}$. Then $\alpha$ intersects
$\mathscr B^{s,+}_j$ and contains one-cells complementary to
$W^s(\vp_{i1}^+)\cap\parabolicmanifold$ for each $i \in \n{3}$.
\end{lem}
\proof The proof of this result follows, up to trivial adaptation,
the proof of by \cite[Lemma 4.5 ]{Moe83}.\finedim
\begin{rem}
Another way to proceed in order to prove the existence of the
heteroclinic connections predicted in item (Q1) (up to
transversality) is to argue as in the proof of \cite[Theorem
3]{Vid99}.
\end{rem}

\newcommand{\cerchio}[1]{\save[]*\frm{#1}\restore}
\begin{figure}
% \entrymodifiers={++[o][F-]}
\begin{center}\SelectTips{xy}{}
\[\xymatrix@C-0pt@R-8pt{%
*+[F-:<3pt>]{\text{$\mathscr B_j^{s,-}$}} & & &  *+[F-:<3pt>]{\text{$\mathscr B_j^{s,+}$}} \\
& *+[F-:<3pt>]{\text{$\vp_{11}^+$}} & & &  \\
& & *+[F-:<3pt>]{\text{$\vp_{21}^+, \vp_{31}^+$}} \ar@/_2ex/[ruu] \ar@/_4ex/[lluu] & &  \\
& &  *+[F-:<3pt>]{\text{$\ve_{11}^+, \ve_{12}^+$}} \ar[u] \ar@/_3ex/[ruuu] \ar@/^2ex/[luu] \ar@/^3ex/[lluuu]&  \\
& &  *+[F-:<3pt>]{\text{$\ve_{11}^-, \ve_{12}^-$}} \ar@/_4ex/[ruuuu]  &  \\
& & *+[F-:<3pt>]{\text{$\vp_{21}^-, \vp_{31}^-$}} \ar[u] \ar@/^4ex/[lluuuuu] \ar@/_6ex/[ruuuuu] & &  \\
& *+[F-:<3pt>]{\text{$\vp_{11}^-$}} \ar@/^2ex/[ruu] \ar@/^4ex/[luuuuuu] & & & \\
*+[F-:<3pt>]{\text{$\mathscr B_j^{u,-}$}}
\ar@/^12ex/[uuuuuuurrr] \ar@/^12ex/[uuuuuur] \ar@/^8ex/[uuuuurr]
\ar@/^6ex/[rruu] \ar@/^2ex/[ur] \ar@/^6ex/[uuurr] & & &
*+[F-:<3pt>]{\text{$\mathscr B_j^{u,-}$}} \ar@/_6ex/[luu] \ar@/_6ex/[luuu]
\ar@/_24ex/[llluuuuuuu] \ar@/_8ex/[luuuuu] \ar@/_8ex/[luuuu] }
\]
\end{center}
\caption{All colliding and non-colliding connecting orbits on the
parabolic manifold with their dual ones between the central
configurations contained in the fundamental domain and the alpha and
omega limit sets.} \label{fig:mf}
\end{figure}
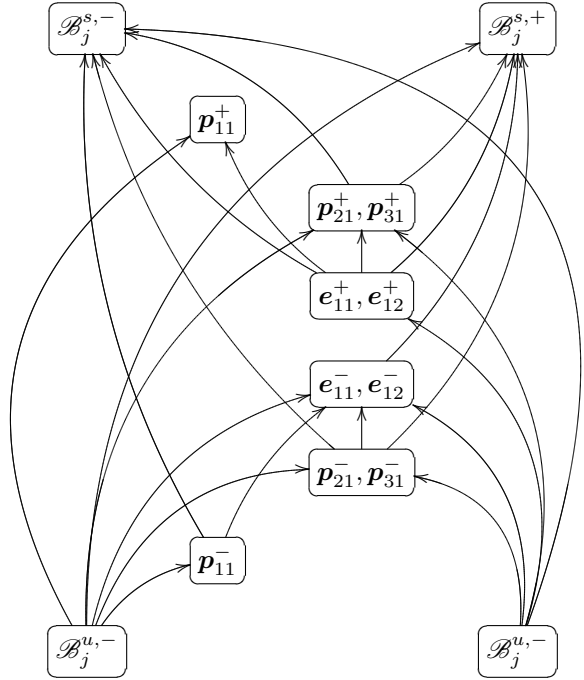

%%%=========================================================================
\section{Behavior of orbits close to quadruple collision}
\label{sec:topologicaltransversality}

In this section we show how to use the connecting orbits on the
$\parabolicmanifold$ in order to establish the motion of the
dihedral four-body problem. We briefly recall the classical  {\em
Palis inclination's lemma\/}.
\begin{lem}
\label{thm:lambdalemma}
Let $p$ be a hyperbolic rest-point of a flow on a manifold $ N$. If a
positively invariant manifold $ P$ has a non-empty topologically
transverse intersection with $W^s(p)$, then $W^u(p) \subset
\overline{ P}$.
\end{lem}
The following result stated by Moeckel in the case of the three body
problem, (cfr. \cite[Proposition]{Moe83} and references therein) has
an analogous in the situation we are dealing with.
\begin{prop}\label{thm:prop2.1}
Let $r,h$ be given constants, with $r>0$. Then there is a constant
$c(r,h)$ such that any orbit of energy $h$ beginning in the region
where $\rho\leq r$ and $v\geq c(r,h)$ is of the following type: two
opposite sides of the rectangle formed by the four bodies are always
the shortest and remains bounded while $\rho$ tends monotonically to
infinity.
\end{prop}
\proof The proof of this result follows by the planar and
tetrahedral energy relations \eqref{eq:energyrelregularizedplanar}
and \eqref{eq:energyrelregularizedtetrahedral2} respectively and by
the fact that at each instant the bodies form an orbit of the
dihedral group $D_2$. In fact, for large values of $v$ the
conservation of the energy forces a small value of the
multiplicative factor of
$\sin^{2\alpha}(2\theta)/W_\planare(\theta)$ in the planar case and
of $(2\cos\varphi)^{2\alpha}/W_\antiprisma(\varphi)$. In each of
these case this condition  correspond to the fact that the
configuration in which two particles are close together. Now, by
taking into account the symmetry of the problem, also the remaining
two are to be closed, and this conclude the proof. \finedim

%For the proof see for instance Lemma 4.2 in \cite{Moe83}.
\begin{prop}\label{thm:primaconclusione} For  any $j \in \n{3}$, the following inclusions between
invariant manifolds hold.
\begin{enumerate}
\item[(i)] $W^s(\vp^-_{j1}) \subset \overline{W^s(\ve^-_{11})}$;
\item[(ii)] $W^u(\vp^+_{j1}) \subset  \overline{W^u(\ve^+_{11})}$.
\end{enumerate}

%\item[(ii)]$W^s(\prisma_i^-) \subset \overline{W^s(\antiprisma^-)}$ and
%$W^u(\prisma_i^+) \subset  \overline{W^u(\antiprisma^+)}$ for
%$i=1,2$.
%\end{enumerate}
\end{prop}
\proof We observe that, by symmetry, it is suffices to show that
$W^u(\vp^+_{j1}) \subset  \overline{W^u(\ve^+_{11})}$, and this will
follow from $W^u(\vp^+_{j1})\cap \parabolicmanifold \subset
\overline{W^u(\ve^+_{11})\cap \parabolicmanifold}$, since the subset
of the full phase space in McGehee coordinates corresponding to
$\{r>0\}$ is open. By arguing as in the proof of \cite[Proposition
5.1]{Moe83}, the result follows by using Lemma \ref{thm:lambdalemma}
and Theorem \ref{thm:diedraletetranoncolliding}.\finedim

\begin{cor}
Arbitrarily close to every planar $\vp^-_{j1}$ type collision orbit
there are tetrahedral $\ve^-_{11}$ type collision orbits. Moreover,
arbitrarily close to every tetrahedral $\ve^+_{11}$ type ejection
orbit there are planar $\vp^+_{j1}$ type ejection orbits.
\end{cor}

%%%=========================================================================
To investigate the behavior of orbits which do not actually collide
we consider the  connections between an  asymptotic set for a
negative value of $v$ and one for a positive value of $v$.
\begin{prop}
Suppose the following connection occur $\mathscr B_i^{u,+}\fromto
\mathscr B_j^{s,+}$ for $i,j \in\n{3}$. Then there is an open set of
orbits with the following behavior: $r(t)\to \infty$ as $t \to \pm
\infty$, but as $t \to -\infty$ two opposite edges of the rectangle
remains bounded (and the remaining ones become unbounded) while for
$t\to + \infty$ the behavior of the opposite edges interchange.
\end{prop}
\proof Any orbit in $\vM_+(h)$ passing close to an orbit in $\vM_0$
connecting $\mathscr B_i^{u,+}$ to $\mathscr B_j^{u,+}$ enters the
domain of applicability of Proposition \ref{thm:prop2.1} in forward
and backward time. The conclusion easily follows. \finedim
% For instance, let us
%denote the vertices of the regular $2l$-gon by $P_k $ for $k=0, 1,
%\dots, 2l-1$.
%\begin{prop}
%Suppose that  there exists the following connection $B^{u}_j\fromto
%B^{s}_i$ on the parabolic manifold $\parabolicmanifold$. Then there
%is an open set of orbits in the interior of $\Sigma$ with the
%following behavior $r(t) \to \infty \ \ \textrm{as}\ \ t \to \pm
%\infty,$ but as $t \to -\infty$ the configuration of the system is
%an irregular $n$-gon in which the edge with vertices $P_{j-1}, P_j$
%let us say $\overline{P_{j-1}P_j}$ remains bounded and as $t \to +
%\infty$ the edge with vertices $P_{i-1}, P_i$ remains bounded.
%
%Suppose that the following connection $B^{u}_j\fromto P^{s}_\pm$
%occur on $\parabolicmanifold$. Then there is an open set of orbits
%in the interior of $\Sigma$ with the following behaviour $r(t) \to
%\infty \ \ \textrm{as}\ \ t \to \pm \infty, $ but as $t \to -\infty$
%the configuration of the system is an irregular $n$-gon in which the
%edge with vertices $P_{j-1}, P_j$ remains bounded and as $t \to +
%\infty$ the configuration of the system is a prism in which the top
%and bottom faces are regular $p$-gons which remains bounded.
%\end{prop}
%The proof of this result follows by integrating the first equation
%of the McGehee equations.
%The last dynamical feature of the problem is a direct consequence of
%the connecting orbits on the parabolic manifold between a central
%configuration for a negative value of $v$ and the asymptotic escape
%sets.
\begin{prop}
Suppose $\vt^-\fromto \mathscr B_j^{s,\pm} $ where $\vt^-$ is any
rectangular or central configuration. Then arbitrarily close to
every $\vt^-$ collision orbit in $\{r>0\}$, there are orbits with
the following behavior: after reaching a minimum, $r(t)$ tends
monotonically to infinity while in a rectangular type configuration
in which two edges become bounded and the others longer and longer.
\end{prop}
\proof Let $\mathscr U$ be any neighborhood of a point $p$ on a
$\vt^-$ collision orbit and let us consider the forward orbits of
the points in $\mathscr U$, namely $\mathscr U\cdot [0, \infty)$,
the union of the forward orbits  of the points on $\mathscr U$.
Since $\mathscr U\cdot [0, \infty)$ is topologically transverse to
$W^s(\vt^-)$ and positively invariant, by using Lemma
\ref{thm:lambdalemma} we conclude that $W^u(\vt^-)\subset
\overline{\mathscr U\cdot [0,\infty)}$ and therefore there are
points in $\mathscr U$ which pass arbitrarily close to a
$\vt^-\fromto \mathscr B_j^{s,\pm}$ connecting orbit. But to any
orbit passing close enough to such an orbit it is possible to apply
Proposition \ref{thm:prop2.1} and the conclusion easily follows.
\finedim

We turn now to a consideration of the more delicate phenomena
associated to a topologically transverse connections between a
central configuration corresponding to a negative value of $v$ and
one corresponding to a  positive value.
\begin{prop}\label{thm:prop5.5}
Suppose that the topologically transverse connection $\vt^-\fromto
\vu^+$ on the parabolic manifold $\parabolicmanifold$ occurs. Let
$\vp\in W^s(\vt^-)$ and $\vq \in W^u(\vu^+)$ and let $\mathscr
U_\vp$ and $\mathscr U_\vq$ be any neighborhoods of these points in
$\vM_+(h)$. Then there is an orbit segments in $\vM_+(h)$ which
begins in $\mathscr U_\vp$, passes close to  a triple collision, and
ends in $\mathscr U_\vq$.
\end{prop}
\proof Since the dimension of the stable and unstable manifolds are
as in the case of the spatial three body problem, the proof of this
result is a trivial adaptation of \cite[Proposition 5.5]{Moe83}.
\finedim
\begin{rem}
Some remarks are in order. First of all we observe that all these
results reveal in a certain sense a chaotic behavior of the total
collision orbits in the sense that all these orbits are extremely
sensitive to small changes in initial conditions. For instance in
the first result of this section we proved that after passing close
to a planar type configuration, the system emerges with arbitrarily
large kinetic energy. In particular in every neighborhood of every
planar type collision orbit there are tetrahedral collision orbits
as well as orbits which avoid collision and emerge from a
neighborhood of the singularity in any of the irregular rectangle
configuration.

Moreover, when a total collapse is reached (asymptotical to a planar
central configuration), the solution oscillates very rapidly passing
through an arbitrarily number of instantaneous planar configurations
(this is due to the spiralling character of the sink) until a double
collision occurs, followed by a passage through a planar central
configuration and further escape from a neighborhood of total
collapse with multiple binary collisions. An approach to total
collision asymptotically to a tetrahedral type central configuration
is very unstable in the sense that after a close encounter is a
double collision followed by a planar configuration and then an
escape from total collision with multiple binary collisions.
\end{rem}

\begin{figure}\centering
\hfill\subfigure[Four bodies with dihedral symmetry and anti-symmetric constraint.]{%
\includegraphics[width=0.3\textwidth]{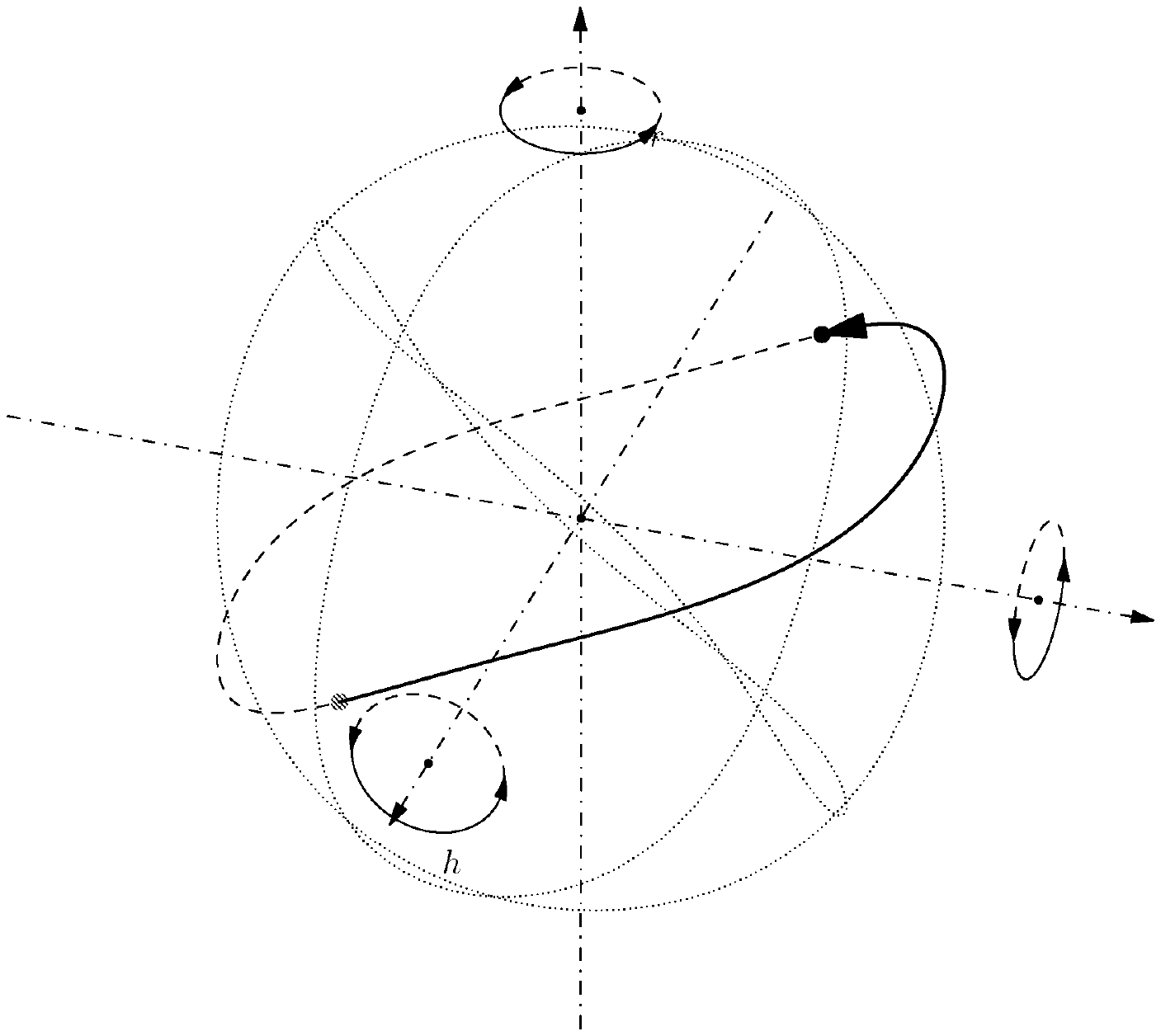}\label{fig:figs/klein4}
}\hfill
\subfigure[Four bodies with dihedral symmetry and semi-coreography constraint.]{%
\includegraphics[width=0.3\textwidth]{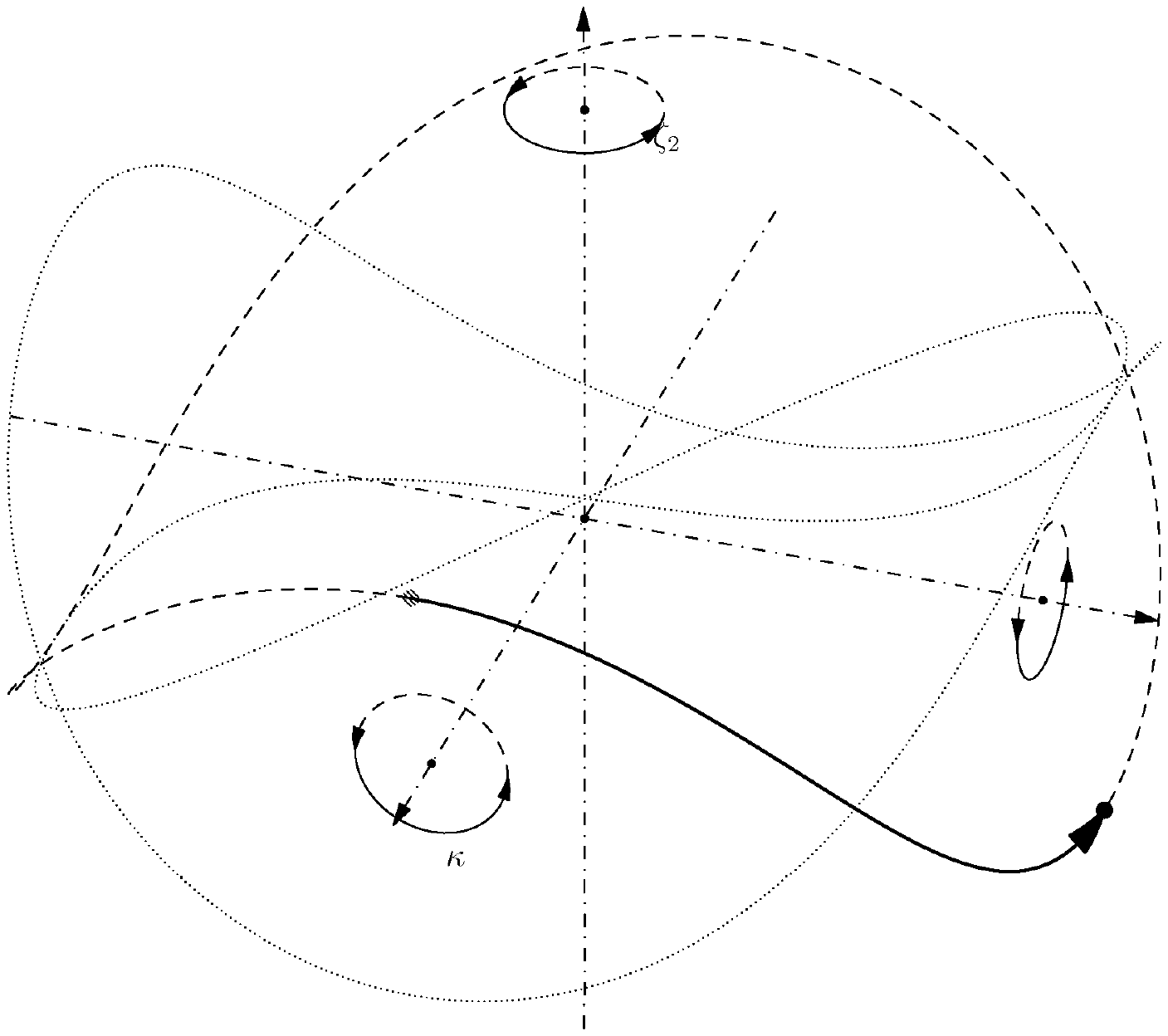}\label{figs/dihedral4}
} \hfill
\subfigure[Twelve bodies with dihedral symmetry and 3-coreography constraint.]{%
\includegraphics[width=0.3\textwidth]{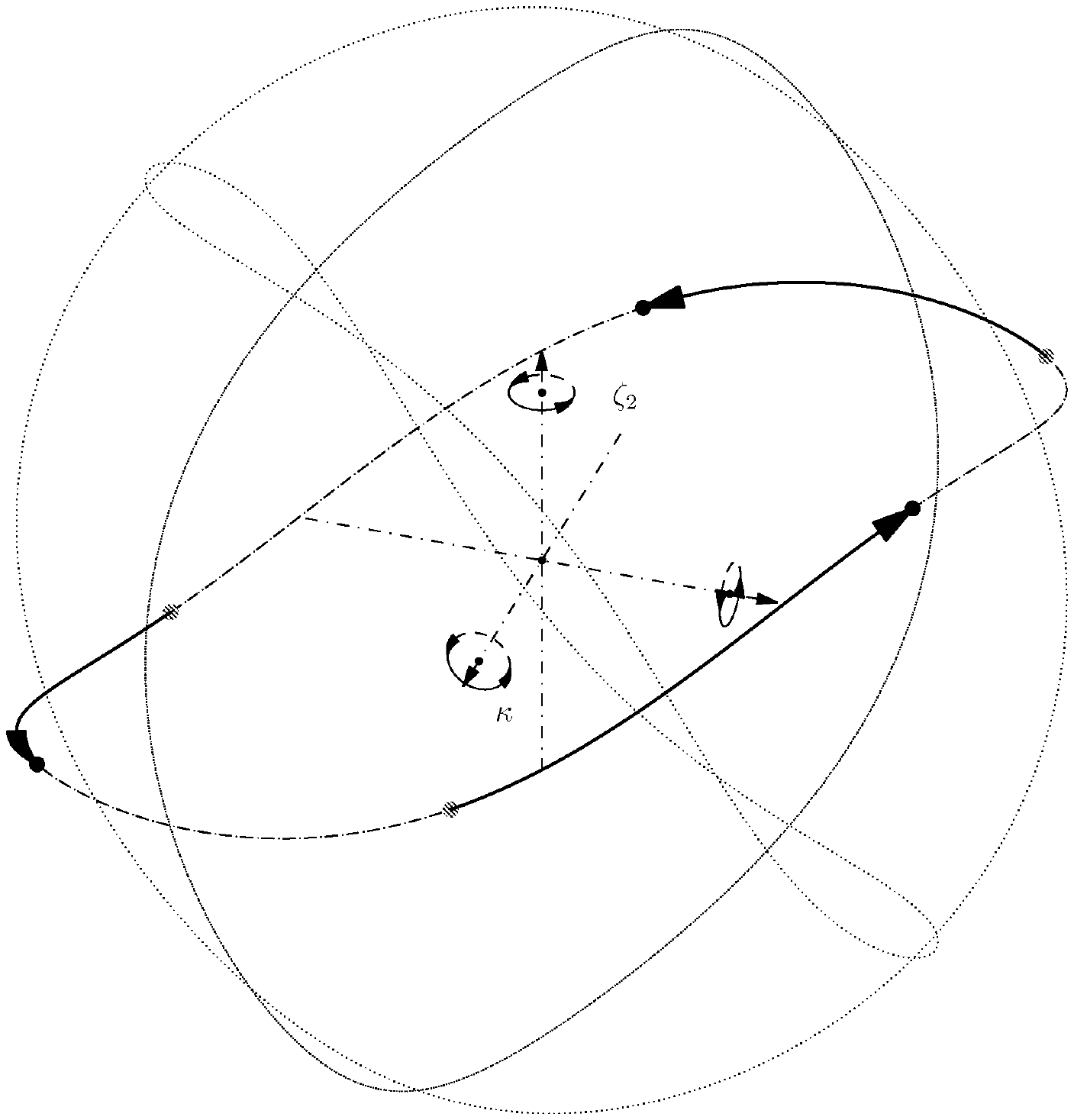}\label{figs/dihedral4_12b_choreography}
} \hfill
\end{figure}

\begin{rem}
Variational methods under symmetry constraints have been popular in
the recent years, starting from Chenciner and Montgomery ``figure
eight'' orbit \cite{Che02,CM00}. Symmetric periodic minimisers have
been found by many authors (to name a few, \cite{ACZ94, Che02,
Che03a, Che03b,SX05}). Periodic orbits dihedrally symmetric can be
shown to exist as well: from \cite{Fer06,Fer07,FT04} it follows that
if $x(t)$ is a minimiser of the Lagrangean action functional,
constrained to the Sobolev space of $G$-equivariant periodic
trajectories (where $G$ is a cyclic extension of $D_l$), then $x(t)$
is collision-less. In the Figures \ref{fig:figs/klein4},
\ref{figs/dihedral4} and \ref{figs/dihedral4_12b_choreography}, some
of the resulting orbits are shown, which are symmetric with respect
to some cyclic extensions of the group $D_l$.
\end{rem}

\appendix

\section{Some pointwise estimates for the planar case}\label{app:planar}

\subsection{Proof of Proposition \ref{prop:flussoplanare} }

The proof of this result will be divided into two main steps. The
first correspond to the Newtonian case $\alpha=1$ and the second to
the case $\alpha \in (0,1)$.

Moreover we refer to the {\em right (left) unstable branch of
$W^u(\vc_{-})$ as that having $u>0$ ($u<0$) in a small neighborhood
of $\vc_{-}$ in $\parabolicmanifold$.\/}

Using uniqueness of solutions, the symmetries and the gradient-like
character of the flow with respect to $v$, it is sufficient to prove
the following fact.
\begin{enumerate}
{\bf \item[(1)] The right branch of $W^u(\vc_{-})$ intersects
$\theta=\pi/2$ ($\theta=-\pi/2$, resp.) with $v<0$ and then
intersects the section $v=0$ with an angle $\theta \in
(0,\pi/4)$.\/}
\end{enumerate}
\subsubsection*{Proof in the Newtonian case.} The proof of this claim
will be divided into some parts:
\begin{enumerate}
\item[(1a)] Upper bound for $v(\pi/2)$ in the case $u>0$ and $\theta \in (\pi/4,
\pi/2)$.
\begin{equation*}
\begin{aligned}
\dfrac{dv}{d\theta}&=
\sqrt{\dfrac14\left[1+\dfrac{1}{\cos\theta}+\dfrac{1}{\sin\theta}\right]-\dfrac{v^2}{4}}\leq\\
&\leq \dfrac12\sqrt{\dfrac{1}{\cos \theta}+\va_2}
\end{aligned}
\end{equation*}
where
\[
\va_2:=\max_{\theta \in [\pi/4, \pi/2]}f(\theta)=1+\sqrt 2, \qquad
f(\theta)= 1 + \dfrac{1}{\sin \theta}.
\]
Thus
\begin{equation*}
\begin{aligned}
v(\pi/2)&\leq v(\pi/4)+ \int_{\pi/4}^{\pi/2}\dfrac12
\sqrt{\dfrac{1}{\cos\theta}+\va_2}\ d\theta  =_{\Phi:=\pi/2-\theta}\\
&=v(\pi/4) + \dfrac12 \int_0^{\pi/4}
\sqrt{\dfrac{1}{\sin\Phi}+\va_2}\ d\Phi
\leq\\
&\leq v(\pi/4) +\dfrac12\int_0^{\pi/4}
\sqrt{\dfrac{\va_1}{\Phi}+\va_2} \ d\Phi
\end{aligned}
\end{equation*}
where $\va_1:=\max_{[0,\pi/4]}\dfrac{\Phi}{\sin
\Phi}=\dfrac{\pi}{4}\sqrt 2$ since the function
$g(\Phi):=\dfrac{\Phi}{\sin\Phi}$ is increasing on $[0, \pi/4]$.
Since $v(\pi/4)=-\sqrt{1+2\sqrt2}$ and by using the quadrature
\[
\int \sqrt{\dfrac{a}{x}+b}\, dx= \sqrt{x(a+bx)}+ \dfrac{a}{2\sqrt
b}\log[2\sqrt b\sqrt{(x(a+bx)}+ 2bx+a],
\]
it follows that
\begin{equation*}
\begin{aligned}
v(\pi/2)\leq    \vv_1:= \dfrac12\int_0^{\pi/4} \sqrt{\dfrac{\pi
\sqrt 2}{4\,x}+1+\sqrt 2}-\sqrt{1+2\sqrt2}%=\\
%1/32\,{\dfrac {\pi
%\, \left( -\sqrt {2}\ln  \left( 2 \right) -2\,\sqrt {2}\ln  \left(
%\pi \right) +4\,\sqrt {1+2\,\sqrt {2}}\sqrt {\sqrt {2}
%+1}\right)}{\sqrt
%{\sqrt {2 }+1}}}+\\
%+\dfrac{\left(2\,\sqrt {2}\ln  \left( 3\,\pi \,\sqrt {2}+2\,\pi
%+2\,\pi \,\sqrt {1+2\,\sqrt {2}}\sqrt {\sqrt {2}+1} \right)  \right)
%}{\sqrt {\sqrt {2
%}+1}}}-\sqrt {1+2\,\sqrt {2}}\\
\approx {\bf -0.8014\/}<0,
\end{aligned}
\end{equation*}
therefore $v$ intersects $\theta=\pi/2$ at a negative value.
%%%%%%%%%%%%%%%%%%%%%5%%%%%%%%%%%%%%%%%%%%%%%%%%%%%%%%%%%%%%%%%%%%%%%%%%%%%%%%%%%%%%%%%%%%%%%%%%%%%%%%%%%%%%%%%
\item[(1b)] Lower bound for $v(\pi/2)$ in the case $u>0$ and $\theta \in (\pi/4,
\pi/2)$.\\
Next we calculate a lower estimate for $v(\pi/2)$. For this we first
observe that
\[
\dfrac{dv}{d\theta}=\dfrac12
\sqrt{1+\dfrac{1}{\sin\theta}+\dfrac{1}{\cos\theta}-v^2}\geq
\dfrac12\sqrt{\dfrac{1}{\cos\theta}+\va_3-v^2},
\]
where $\va_3:= \min_{\theta \in [\pi/4, \pi/2]}f(\theta)=2$. Now we
consider the interval $J:=[3/8\pi, \pi/2]\subset[\pi/4,
\pi/2]$\footnote{%
We observe that the value $\bar \theta:=3/8\pi$ has been chosen
suitable in order the estimates to work and in such a way the
expression under square root is nonnegative.\/} and by a direct
calculation:
\begin{equation*}
\begin{aligned}
\dfrac{dv}{d\theta}\geq \dfrac12\sqrt{\dfrac{1}{\cos\theta}+2-
v^2(\pi/4)}\ \Rightarrow v(\pi/2) &\geq v(3/8\pi) +
\dfrac12\int_{3/8\pi}^{\pi/2}\sqrt{\dfrac{1}{\cos\theta}+2-
v^2(\pi/4)}\ d\theta\geq\\
& \geq v(\pi/4) +
\dfrac12\int_{0}^{\pi/8}\sqrt{\dfrac{1}{\Phi}+\va_4}\ d\Phi\\
&\approx {\bf-1.4164}:=\vv_2<0,
\end{aligned}
\end{equation*}
for $\Phi$ as above and $\va_4:=2-v^2(\pi/4)=1-2\sqrt 2$; moreover
in the last inequality we used the fact that $\sin \Phi \leq \Phi$.
From (1a) \&(1b) it follows that $\vv_2 \leq v(\theta)\leq \vv_1$,
for $\theta \in [\pi/4, \pi/2]$ and $u>0$.
%%%%%%%%%%%%%%%%%%%%%%%%%%%%%%%%%%%%%%%%%%%%%%%%%%%%%%%%%%%%%%%%%%%%%%%%%%%%%%%%%%%
\item[(1c)] The continuation of the right branch through double
collision at $\theta=\pi/2$ for $u<0$. We affirm that $v$ remains
negative in the interval $J$.

In fact, for $u<0$
\[
-\dfrac{dv}{d\theta} \leq \dfrac12\sqrt{\dfrac{1}{\cos\theta}+\va_5}
\]
for $\va_5:=f(3/8\pi)=1+2\,{\dfrac {1}{\sqrt {2+\sqrt {2}}}}\approx
{2.0824}$ and$f$ as above. Let $\va_6:=g(\pi/8)=3/4\,{\dfrac {\pi
}{\sqrt {2+\sqrt {2}}}}\approx {1.0260}$. Integrating we get:
\begin{equation}
\begin{aligned}
v(3/8\pi)&\leq v(\pi/2)
+\dfrac12\int_{3/8\pi}^{\pi/2}\sqrt{\dfrac{1}{\cos\theta}+\va_5}\ d\theta \leq\\
&\leq
v(\pi/2)+\dfrac12\int_0^{\pi/8}\sqrt{\dfrac{\va_6}{\Phi}+\va_5}\
d\Phi \leq \vv_1
+\dfrac12\int_0^{\pi/8}\sqrt{\dfrac{\va_6}{\Phi}+\va_5}\
d\Phi \\
&\approx \vv_1+ {0.7111}\approx {\bf-0.0903}<0
\end{aligned}
\end{equation}
Therefore
\[
\vv_2 \leq v(\pi/2)\leq v(\theta)\leq v(3/8\pi)<0, \qquad
\textrm{for}\ \ \theta \in [3/8\pi, \pi/2], \ \ \textrm{and}\ \ u<0.
\]
Therefore
\[
\vv_2 \leq v(\pi/2)\leq v(\theta)\leq v(3/8\pi)<0, \qquad
\textrm{for}\ \ \theta \in [3/8\pi, \pi/2], \ \ \textrm{and}\ \ u<0.
\]
%%%%%%%%%%%%%%%%%%%%%%%%%%%%%%%%%%%%%%%%%%%%%%%%%%%%%%%%%%%%%%%%%%%%%%%%%%%%%%%%%%%
\item[(1d)] The continuation of the right branch through double
collision at $\theta=\pi/2$ for $u<0$. We affirm that $v$ remains
negative at the point $\theta=\pi/4$.

In fact, for $u<0$
\[
-\dfrac{dv}{d\theta} \leq \dfrac12\sqrt{\dfrac{1}{\cos\theta}+\va_5}
\]
for $\va:=f(\pi/4)=1+2\sqrt2$ and$f$ as above. Let
$\vb:=g(\pi/4)=\pi/4\sqrt2$. Integrating we get:
\begin{equation}
\begin{aligned}
v(\pi/4)&\leq v(\pi/2)
+\dfrac12\int_{\pi/4}^{\pi/2}\sqrt{\dfrac{1}{\cos\theta}+\va}\ d\theta \leq\\
&\leq v(\pi/2)+\dfrac12\int_0^{\pi/4}\sqrt{\dfrac{\vb}{\Phi}+\va}\
d\Phi \leq \vv_1
+\dfrac12\int_0^{\pi/4}\sqrt{\dfrac{\vb}{\Phi}+\va}\
d\Phi \\
&\approx \vv_1+ {0.7445}={\bf-0.5630}<0
\end{aligned}
\end{equation}

%%%%%%%%%%%%%%%%%%%%%%%%%%%%%%%%%%%%%%%%%%%%%%%%%%%%%%%%%%%%%%%%%%%%%%%%%%%%%%%%%%%
\item[(1e)] From the above inequality we get the estimates $\vv_2^2\geq
v(\theta)^2$, for all $\theta \in [3\pi/8, \pi/2]$ and $u<0$. Using
this estimate  we can write:
\[
-\dfrac{dv}{d\theta}=\dfrac12
\sqrt{1+\dfrac{1}{\sin\theta}+\dfrac{1}{\cos\theta}-v^2}\geq\dfrac12\sqrt{\dfrac{1}{\cos\theta}+\va_7
-\vv_2^2}
\]
where $\va_7:=\min_{\theta \in [3/8\pi,
\pi/2]}f(\theta)=f(\pi/2)=2$. By a direct integration we have:
\begin{equation*}
\begin{aligned}
v(3/8\pi) &\geq v(\pi/2) +
\dfrac12\int_{3/8\pi}^{\pi/2}\sqrt{\dfrac{1}{\cos\theta}+\va_7
-\vv_2^2}\ d\theta \geq\\
&\geq v(\pi/2) + \dfrac12 \int_0^{\pi/8}
\sqrt{\dfrac{1}{\Phi}+2-\vv_2^2}\
d\Phi\\
&\geq \vv_2+ \dfrac12 \int_0^{\pi/8}
\sqrt{\dfrac{1}{\Phi}+2-\vv_2^2}\ d\Phi \approx{\bf-0.7900}.
\end{aligned}
\end{equation*}
\item[(1f)]To finish the estimate for the continuation of the right
branch through double collision, we now consider $\theta \in
[0,3\pi/8]$ and $u<0$. We observe that in $\theta=\pi/4$,
$U(\theta)$ has a relative minimum and therefore for $\theta \in
[0,3\pi/8]$ we have $U(\theta)\geq U(\pi/4)$. Now let $\tilde
v:=\sqrt{2U(\pi/4)}$; thus
\[
-\dfrac{dv}{d\theta}\geq \dfrac12\sqrt{\tilde v^2-v^2}
\]
and by a direct integration, we have
\[
-\arcsin\left(\dfrac{v}{\tilde v}\right)\Big\vert_{v(0)}^{v(3\pi/8)}
\geq \dfrac12 \int_{0}^{3\pi/8} d \theta=3\pi/16.
\]
Therefore
\[
v(0)\geq \tilde v\sin\left(3\pi/16+\arcsin(v(3\pi/8)/\tilde
v)\right)\geq \tilde v\sin\left(3\pi/16+\arcsin(-0.7900/\tilde
v)\right) \approx{\bf0.3376}.
\]
We observe that this estimate could be refined by considering the
interval $[0,\pi/4]$ instead of $[0, 3\pi/8]$; however to do so we
need an estimate from below for $v(\pi/4)$.
\end{enumerate}
By symmetry, the same arguments holds also for the left branch of
$W^u(\vc_{-})$. Now the thesis easily follows by the fact that the
flow is of gradient-type with restpoints $\vc_{-}$ and $\vc_{+}$.
\finedim

\subsubsection*{Proof in the  case $\alpha \in (0,1)$}

Also in this case the proof of this claim will be divided into some
parts:
\begin{enumerate}
\item[(1a)] Upper bound for $v(\pi/2)$ in the case $u>0$ and $\theta \in (\pi/4,
\pi/2)$.
\begin{equation*}
\begin{aligned}
\dfrac{dv}{d\theta}&=\dfrac12
\sqrt{2^{1-\alpha}\left[1+\dfrac{1}{\cos^\alpha\theta}+\dfrac{1}{\sin^\alpha\theta}\right]-v^2}\leq
\dfrac{1}{2^{1/2+\beta}}\sqrt{\left[1+\dfrac{1}{\cos^\alpha\theta}+\dfrac{1}{\sin^\alpha\theta}\right]}\leq\\
&\leq \dfrac{1}{2^{1/2+\beta}}\sqrt{\dfrac{1}{\cos^\alpha
\theta}+\vb_2}
\end{aligned}
\end{equation*}
where
\[
\vb_2:=\max_{\theta \in [\pi/4, \pi/2]}h(\theta)=1+2^\beta, \qquad
h(\theta):= 1 + \dfrac{1}{\sin^\alpha \theta}.
\]
Thus
\begin{equation*}
\begin{aligned}
v(\pi/2)&\leq v(\pi/4)+ \int_{\pi/4}^{\pi/2}\dfrac{1}{2^{1/2+\beta}}
\sqrt{\dfrac{1}{\cos^\alpha\theta}+\vb_2}\ d\theta  =_{\Phi:=\pi/2-\theta}\\
&=v(\pi/4) + \dfrac{1}{2^{1/2+\beta}}\int_0^{\pi/4}
\sqrt{\dfrac{1}{\sin^\alpha\Phi}+\vb_2}\ d\Phi
\leq\\
&\leq v(\pi/4) +\dfrac{1}{2^{1/2+\beta}}\int_0^{\pi/4}
\sqrt{\dfrac{\vb_1}{\Phi^\alpha}+\vb_2} \ d\Phi\leq\\
&\leq -\dfrac{1}{2^{\beta-1/2}}\sqrt{1+2^{\beta+1}}
+\dfrac{1}{2^{1/2+\beta}}\int_0^{\pi/4}
\sqrt{\dfrac{\vb_1}{\Phi^{2\beta}}+\vb_2} \ d\Phi:=\vj_1(\beta),
\end{aligned}
\end{equation*}
where $\vb_1:=\max_{[0,\pi/4]}\left(\dfrac{\Phi}{\sin
\Phi}\right)^\alpha=\left(\dfrac{\pi}{4}\sqrt 2\right)^\alpha$ since
$g(\Phi):=\left(\dfrac{\Phi}{\sin\Phi}\right)^\alpha$ is an
increasing function. Since  $(0,1/2)\ni\beta\mapsto \vj_1(\beta)$ is
an increasing function on the interval $(0, 1/2)$, then we conclude
that
\[
v(\pi/2)\leq \vj_1(1/2)\approx {\bf-0.8013}:=\vw_1.
\]

Therefore $v$ intersects $\theta=\pi/2$ at a negative value.
%%%%%%%%%%%%%%%%%%%%%5%%%%%%%%%%%%%%%%%%%%%%%%%%%%%%%%%%%%%%%%%%%%%%%%%%%%%%%%%%%%%%%%%%%%%%%%%%%%%%%%%%%%%%%%%
\item[(1b)] Lower bound for $v(\pi/2)$ in the case $u>0$ and $\theta \in (\pi/4,
\pi/2)$.

Next we calculate a lower estimate for $v(\pi/2)$. For this we first
observe that
\[
\dfrac12
\sqrt{2^{1-\alpha}\left[1+\dfrac{1}{\cos^\alpha\theta}+\dfrac{1}{\sin^\alpha\theta}\right]-v^2}\geq
\dfrac{1}{2^{1/2+\beta}}\sqrt{\dfrac{1}{\cos^\alpha\theta}+\vb_3-v^2},
\]
where $\vb_3:= \min_{\theta \in [\pi/4, \pi/2]}h(\theta)=2$. Now we
consider the interval $J:=[\pi/4, 3\pi/8]$ and by a direct
calculation, we have:
\begin{equation*}
\begin{aligned}
\dfrac{dv}{d\theta}&\geq
\dfrac{1}{2^{1/2+\beta}}\sqrt{\dfrac{1}{\cos^\alpha\theta}+2-
v^2(\pi/4)}\ \Rightarrow\\ v(\pi/2) &\geq v(3\pi/8) +
\int_{3\pi/8}^{\pi/2}\dfrac{1}{2^{1/2+\beta}}\sqrt{\dfrac{1}{\cos^\alpha\theta}+2-
v^2(\pi/4)}\ d\theta\geq\\
&\geq v(3\pi/8) +
\int_{0}^{\pi/8}\dfrac{1}{2^{1/2+\beta}}\sqrt{\dfrac{1}{\Phi^{2\beta}}+1-2^{\beta+1}}\
d\Phi:=\vj_2(\beta).
\end{aligned}
\end{equation*}
Since:
\[
v(3\pi/8):=v_\beta(3\pi/8)=-\sqrt{2^{1-2\beta}\left[1+
\left(\dfrac{2\sqrt2}{1+\sqrt2}\right)^\beta+\left(\dfrac{2\sqrt2}{\sqrt2-1}\right)^\beta
\right]};
\]
and since both the functions
\[
(0,1/2)\ni\beta\mapsto v_\beta(3\pi/8)\ \textrm{and}\ (0, 1/2)\ni
\beta\mapsto
\int_{0}^{\pi/8}\dfrac{1}{2^{1/2+\beta}}\sqrt{\dfrac{1}{\Phi^{2\beta}}+1-2^{\beta+1}}
\]
are strictly increasing it immediately follows that also $\vj_2$ is
an increasing function and therefore
%\[
%v_{1/2}(3\pi/8)\approx \textcolor{red}{-2.449}\leq v(3\pi/8)\leq
%v_0(3\pi/8) \approx \textcolor{red}{-2.167}.
%\]
%To conclude the estimate, we observe that also the function:
%\[
%(0,1/2)\ni\beta\mapsto
%\vi(\beta):=\int_{0}^{\pi/8}\dfrac{1}{2^{1/2+\beta}}\sqrt{\dfrac{1}{\theta^{2\beta}}+1-
%2^{(\beta+1)}}\ d\theta\in \R^+
%\]
%is a positive and increasing function, with $\vi(0)=0$. Therefore,
we can  conclude that
\[
v(\pi/2)\geq \vj_1(1/2):=\vw_2\approx{\bf-1.6267}.
\]
%%%%%%%%%%%%%%%%%%%%%%%%%%%%%%%%%%%%%%%%%%%%%%%%%%%%%%%%%%%%%%%%%%%%%%%%%%%%%%%%%%%
\item[(1c)] The continuation of the right branch through double
collision at $\theta=\pi/2$ for $u<0$. We affirm that $v$ remains
negative in the interval $J$.

In fact, for $u<0$
\[
-\dfrac{dv}{d\theta} \leq
\dfrac12\sqrt{2^{1-\alpha}\left[\dfrac{1}{\cos^\alpha\theta}+\vb_5\right]}
\]
for $\vb_5:=h(3\pi/8)$ and $h$ as above. Let $\vb_6:=g(\pi/8)$.
Integrating we get:
\begin{equation}
\begin{aligned}
v(3\pi/8)&\leq v(\pi/2)
+\dfrac12\int_{3\pi/8}^{\pi/2}\dfrac{1}{2^{1/2+\beta}}\sqrt{\left[\dfrac{1}{\cos^\alpha\theta}+\vb_5\right]}\ d\theta \leq\\
&\leq
v(\pi/2)+\dfrac12\int_0^{\pi/8}\dfrac{1}{2^{1/2+\beta}}\sqrt{\left[\dfrac{\vb_6}{\Phi^{2\beta}}+\vb_5\right]}\
d\Phi \leq\\
\leq& \vw_1
+\dfrac12\int_0^{\pi/8}\dfrac{1}{2^{1/2+\beta}}\sqrt{\left[\dfrac{\vb_6}{\Phi^{2\beta}}+\vb_5\right]}\
d\Phi\\
%v_1
%+\dfrac{1}{\sqrt2}\int_0^{\pi/16}\sqrt{\dfrac{\va_6}{\theta}+\va_5}\
%d\theta \\
&\approx {\bf-0.0902}.
\end{aligned}
\end{equation}
Therefore
\[
\vw_2\leq v(\pi/2)\leq v(\theta)\leq v(3\pi/8)<0, \qquad
\textrm{for}\ \ \theta \in [3\pi/8, \pi/2], \ \ \textrm{and}\ \ u<0.
\]

%%%%%%%%%%%%%%%%%%%%%%%%%%%%%%%%%%%%%%%%%%%%%%%%%%%%%%%%%%%%%%%%%%%%%%%%%%%%%%%%%%%

\item[(1d)] From the above inequality we get the estimates $\vw_2^2\geq
v(\theta)^2$, for all $\theta \in [3\pi/8, \pi/2]$ and $u<0$. Using
this estimate  we can write:
\[
-\dfrac{dv}{d\theta}=\dfrac12
\sqrt{2^{1-\alpha}\left[1+\dfrac{1}{\sin^\alpha\theta}+\dfrac{1}{\cos^\alpha\theta}\right]-v^2}
\geq\dfrac12\sqrt{2^{1-\alpha}\left[\dfrac{1}{\cos\theta}+\va_7\right]
-\vw_2^2}
\]
where $\va_7:=\min_{\theta \in [3\pi/8,
\pi/2]}h(\theta)=h(\pi/2)=2$. By a direct integration we have:
\begin{equation*}
\begin{aligned}
v(3\pi/8) &\geq v(\pi/2) +
\dfrac12\int_{3\pi/8}^{\pi/2}\sqrt{2^{1-\alpha}\left[\dfrac{1}{\cos^\alpha\theta}+\va_7\right]
-\vw_2^2}\ d\theta \geq\\
&\geq v(\pi/2) + \dfrac12 \int_0^{\pi/8}
\sqrt{2^{1-\alpha}\left[\dfrac{1}{\Phi^\alpha}+2\right]-\vw_2^2}\
d\Phi\geq\\ &\geq \vw_2 + \dfrac12 \int_0^{\pi/8}
\sqrt{2^{1-\alpha}\left[\dfrac{1}{\Phi^\alpha}+2\right]-\vw_2^2}\
d\Phi\geq \\
\geq&\vw_2 + \dfrac12 \int_0^{\pi/8}
\sqrt{\left[\dfrac{1}{\Phi}+2\right]-\vw_2^2}\ d\Phi \approx {\bf
-1.0276}.
\end{aligned}
\end{equation*}
%%%%%%%%%%%%%%%%%%%%%%%%%%%%%%%%%%%%%%%%%%%%%%%%%%%%%%%%%%%%%%%%%%%%%%%%%%%%%%%%%%%%%%%%%%%%%%%%%%%%%%%%%%

\item[(1e)] The continuation of the right branch through double
collision at $\theta=\pi/2$ for $u<0$. We affirm that $v$ remains
negative at the point $\theta=\pi/4$.

In fact, for $u<0$
\[
-\dfrac{dv}{d\theta} \leq
\dfrac12\sqrt{2^{1-\alpha}\left[\dfrac{1}{\cos^\alpha\theta}+\vb_5\right]}
\]
for $\vb_5:=h(\pi/4)$ and $h$ as above. Let $\vb_6:=g(\pi/4)$.
Integrating we get:
\begin{equation}
\begin{aligned}
v(\pi/4)&\leq v(\pi/2)
+\dfrac12\int_{\pi/4}^{\pi/2}\dfrac{1}{2^{1/2+\beta}}\sqrt{\left[\dfrac{1}{\cos^\alpha\theta}+\vb_5\right]}\ d\theta \leq\\
&\leq
v(\pi/2)+\dfrac12\int_0^{\pi/4}\dfrac{1}{2^{1/2+\beta}}\sqrt{\left[\dfrac{\vb_6}{\Phi^{2\beta}}+\vb_5\right]}\
d\Phi \leq\\
\leq& \vw_1
+\dfrac12\int_0^{\pi/4}\dfrac{1}{2^{1/2+\beta}}\sqrt{\left[\dfrac{\vb_6}{\Phi^{2\beta}}+\vb_5\right]}\
d\Phi\\
%v_1
%+\dfrac{1}{\sqrt2}\int_0^{\pi/16}\sqrt{\dfrac{\va_6}{\theta}+\va_5}\
%d\theta \\
&\approx {\bf-0.2237}.
\end{aligned}
\end{equation}

\item[(1f)]To finish the estimate for the continuation of the right
branch through double collision, we now consider $\theta \in [0,
\pi/8]$ and $u<0$. We observe that in $\theta=\pi/8$, $U(\theta)$
has a relative minimum and therefore for $\theta \in [0, \pi/8]$ we
have $U(\theta)\geq U(\pi/8)$. Now let $\tilde v:=\sqrt{2U(\pi/8)}$;
thus
\[
-\dfrac{dv}{d\theta}\geq \dfrac12\sqrt{\tilde v^2-v^2}
\]
and by a direct integration, we have
\[
-\arcsin\left(\dfrac{v}{\tilde v}\right)\Big\vert^{v(\pi/8)}_{v(0)}
\geq \dfrac12 \int_{0}^{\pi/8} d \theta=\pi/16.
\]
Therefore
\[
v(0)\geq \tilde v\sin\left(\pi/16+
\arcsin\left(\dfrac{v(\pi/8)}{\tilde v}\right)\right)\geq \tilde
v\sin\left(\pi/16+ \arcsin\left(\dfrac{-0.3741}{\tilde
v}\right)\right) \approx {\bf 0.3810}.
\]
\end{enumerate}
By symmetry, the same arguments holds also for the left branch of
$W^u(\vc_{-})$. Now the thesis easily follows by the fact that the
flow is of gradient-type with restpoints $\vc_{-}$ and
$\vc_{+}$.\finedim

%%%================================================================================================================

\subsection{Proof of Proposition \ref{prop:flussoplanareb}}

The proof of this result will be divided into several steps.
\begin{enumerate}
\item {\bf There exists $\alpha_0 \in (1,2)$ such that for each $\alpha
\in (1, \alpha_0)$ the right branch of $W^u(\vc_{-})$ intersects
$\theta=\pi/2$ ($\theta=-\pi/2$, resp.) with $v<0$ and then
intersects the section $v=0$ with an angle $\theta \in (0,\pi/4)$.}
\item {\bf There exists $\alpha_1 \in (1,2)$ such that for each $\alpha
\in (\alpha_1, 2)$ the right branch of $W^u(\vc_{-})$ intersects
$\theta=\pi/2$ ($\theta=-\pi/2$, resp.) with $v<0$ and then
intersects the section $v=0$ with an angle $\theta \in (\pi/4,
\pi/2)$.}
\item {\bf There exists $\alpha_* \in (\alpha_0,\alpha_1)$ such that
the right branch of $W^u(\vc_{-})$ intersects the right branch of
$W^u(\vc_{+})$.}
\end{enumerate}
In order to prove item 1 we shall proceed as in the proof of the
previous proposition.
\begin{enumerate}
\item[(1a)] {\bf There exists $\alpha_0 \in (1,2)$ such that for
each $\alpha \in (1, \alpha_0)$, $v(\pi/2)<0$.\/}
\end{enumerate}
\begin{itemize}
\item  Upper bound for $v(\pi/2)$ in the case $u>0$ and $\theta \in (\pi/4,
\pi/2)$.
\begin{equation*}
\begin{aligned}
\dfrac{dv}{d\theta}&=\dfrac12
\sqrt{2^{1-\alpha}\left[1+\dfrac{1}{\cos^\alpha\theta}+\dfrac{1}{\sin^\alpha\theta}\right]-v^2}\leq
\dfrac{1}{2^{1/2+\beta}}\sqrt{\left[1+\dfrac{1}{\cos^\alpha\theta}+\dfrac{1}{\sin^\alpha\theta}\right]}\leq\\
&\leq \dfrac{1}{2^{1/2+\beta}}\sqrt{\dfrac{1}{\cos^\alpha
\theta}+\vb_2}
\end{aligned}
\end{equation*}
where
\[
\vb_2:=\max_{\theta \in [\pi/4, \pi/2]}h(\theta)=1+2^\beta, \qquad
h(\theta):= 1 + \dfrac{1}{\sin^\alpha \theta}.
\]
Thus
\begin{equation*}
\begin{aligned}
v(\pi/2)&\leq v(\pi/4)+ \int_{\pi/4}^{\pi/2}\dfrac{1}{2^{1/2+\beta}}
\sqrt{\dfrac{1}{\cos^\alpha\theta}+\vb_2}\ d\theta  =_{\Phi:=\pi/2-\theta}\\
&=v(\pi/4) + \dfrac{1}{2^{1/2+\beta}}\int_0^{\pi/4}
\sqrt{\dfrac{1}{\sin^\alpha\Phi}+\vb_2}\ d\Phi
\leq\\
&\leq v(\pi/4) +\dfrac{1}{2^{1/2+\beta}}\int_0^{\pi/4}
\sqrt{\dfrac{\vb_1}{\Phi^\alpha}+\vb_2} \ d\Phi=\\
&= -\dfrac{1}{2^{\beta-1/2}}\sqrt{1+2^{\beta+1}}
+\dfrac{1}{2^{1/2+\beta}}\int_0^{\pi/4}
\sqrt{\dfrac{\vb_1}{\Phi^{2\beta}}+\vb_2} \ d\Phi:=\vk_1(\beta),
\end{aligned}
\end{equation*}
where $\vb_1:=\max_{[0,\pi/4]}\left(\dfrac{\Phi}{\sin
\Phi}\right)^\alpha=\left(\dfrac{\pi}{4}\sqrt 2\right)^\alpha$ since
$g(\Phi):=\left(\dfrac{\Phi}{\sin\Phi}\right)^\alpha$ is an
increasing function. We observe that  $(1/2,1)\ni\beta\mapsto
\vk_1(\beta)$ is a continuous  increasing function on the interval
$(1/2, 1)$. However it changes sign in the interval $A_1:=(1.4,
1.6)$;
%%$A:=(0.73068, 0.73029)$;
in particular if $\alpha_0$ is any number in the open interval $(1,
1.46136]$, then we can conclude that
\[
v(\pi/2)\leq \vk_1(1/2)\approx {\bf-0.1659}:=\vw_1.
\]

\item Lower bound for $v(\pi/2)$ in the case $u>0$ and $\theta \in (\pi/4,
\pi/2)$. For this we first observe that
\[
\dfrac12
\sqrt{2^{1-\alpha}\left[1+\dfrac{1}{\cos^\alpha\theta}+\dfrac{1}{\sin^\alpha\theta}\right]-v^2}\geq
\dfrac{1}{2^{1/2+\beta}}\sqrt{\dfrac{1}{\cos^\alpha\theta}+\vb_3-2^{\alpha-1}v^2},
\]
where $\vb_3:= \min_{\theta \in [\pi/4, \pi/2]}h(\theta)=2$. Now we
consider the interval $J:=[\pi/4, 3\pi/8]$ and by a direct
calculation, we have:
\begin{equation*}
\begin{aligned}
\dfrac{dv}{d\theta}&\geq
\dfrac{1}{2^{1/2+\beta}}\sqrt{\dfrac{1}{\cos^\alpha\theta}+2-
v^2(\pi/4)}\ \Rightarrow\\ v(\pi/2) &\geq v(3\pi/8) +
\int_{3\pi/8}^{\pi/2}\dfrac{1}{2^{1/2+\beta}}\sqrt{\dfrac{1}{\cos^\alpha\theta}+2-
v^2(\pi/4)}\ d\theta\geq\\
&\geq v(3\pi/8) +
\int_{0}^{\pi/8}\dfrac{1}{2^{1/2+\beta}}\sqrt{\dfrac{1}{\Phi^{2\beta}}+1-2^{\beta+1}}\
d\Phi:=\vj_2(\beta).
\end{aligned}
\end{equation*}
Since:
\[
v(3\pi/8):=v_\beta(3\pi/8)=-\sqrt{2^{1-2\beta}\left[1+
\left(\dfrac{2\sqrt2}{1+\sqrt2}\right)^\beta+\left(\dfrac{2\sqrt2}{\sqrt2-1}\right)^\beta
\right]};
\]
and since both the functions
\[
(0,1/2)\ni\beta\mapsto v_\beta(3\pi/8)\ \textrm{and}\ (0, 1/2)\ni
\beta\mapsto
\int_{0}^{\pi/8}\dfrac{1}{2^{1/2+\beta}}\sqrt{\dfrac{1}{\Phi^{2\beta}}+1-2^{\beta+1}}
\]
are strictly increasing it immediately follows that also $\vj_2$ is
an increasing function and therefore
%\[
%v_{1/2}(3\pi/8)\approx \textcolor{red}{-2.449}\leq v(3\pi/8)\leq
%v_0(3\pi/8) \approx \textcolor{red}{-2.167}.
%\]
%To conclude the estimate, we observe that also the function:
%\[
%(0,1/2)\ni\beta\mapsto
%\vi(\beta):=\int_{0}^{\pi/8}\dfrac{1}{2^{1/2+\beta}}\sqrt{\dfrac{1}{\theta^{2\beta}}+1-
%2^{(\beta+1)}}\ d\theta\in \R^+
%\]
%is a positive and increasing function, with $\vi(0)=0$. Therefore,
we can  conclude that
\[
v(\pi/2)\geq \vj_1(1/2):=\vw_2\approx{\bf-1.6267}.
\]
We observe that for $\alpha \geq 1.1$ we have $ v(\pi/2)\geq
\vj_1(0.55):=\vw_3\approx{\bf-1.5285}$
Therefore $v$ intersects $\theta=\pi/2$ at a negative value.
%%%%%%%%%%%%%%%%%%%%%5%%%%%%%%%%%%%%%%%%%%%%%%%%%%%%%%%%%%%%%%%%%%%%%%%%%%%%%%%%%%%%%%%%%%%%%%%%%%%%%%%%%%%%%%%

\end{itemize}
%%%%%%%%%%%%%%%%%%%%%%%%%%%%%%%%%%%%%%%%%%%%%%%%%%%%%%%%%%%%%%%%%%%%%%%%%%%%%%%%%%%
\begin{enumerate}
\item[(1b)] {\bf There exists $\alpha_{00} \in (1, \alpha_0)$ such that the continuation of the right branch through double
collision at $\theta=\pi/2$ for $u<0$, remains negative in the
interval $[\pi/4, \pi/2]$.\/}
\end{enumerate}
In fact, for $u<0$
\[
-\dfrac{dv}{d\theta} \leq
\dfrac12\sqrt{2^{1-\alpha}\left[\dfrac{1}{\cos^\alpha\theta}+\vb_5\right]}
\]
for $\vb_5:=h(49\pi/100)$ and $h$ as above. Let $\vb_6:=g(\pi/100)$.
Integrating we get:
\begin{equation}
\begin{aligned}
v(49\pi/100)&\leq v(\pi/2)
+\dfrac12\int_{49\pi/100}^{\pi/2}\dfrac{1}{2^{1/2+\beta}}\sqrt{\left[\dfrac{1}{\cos^\alpha\theta}+\vb_5\right]}\ d\theta \leq\\
&\leq
v(\pi/2)+\dfrac12\int_0^{\pi/100}\dfrac{1}{2^{1/2+\beta}}\sqrt{\left[\dfrac{\vb_6}{\Phi^{2\beta}}+\vb_5\right]}\
d\Phi \leq\\
\leq& \vw_1
+\dfrac12\int_0^{\pi/100}\dfrac{1}{2^{1/2+\beta}}\sqrt{\left[\dfrac{\vb_6}{\Phi^{2\beta}}+\vb_5\right]}\
d\Phi\leq\\
&\leq -0.1659+
\dfrac12\int_0^{\pi/100}\dfrac{1}{2^{1/2+\beta}}\sqrt{\left[\dfrac{\vb_6}{\Phi^{2\beta}}+\vb_5\right]}\
d\Phi\\
%v_1
%+\dfrac{1}{\sqrt2}\int_0^{\pi/16}\sqrt{\dfrac{\va_6}{\theta}+\va_5}\
%d\theta \\
&\approx {\bf-0.0057},
\end{aligned}
\end{equation}
which holds for $\alpha \in (1,1.1]$. Thus the thesis follows by
choosing $\alpha_{00}=1.1$. Therefore
\[
\vw_2\leq v(\pi/2)\leq v(\theta)\leq v(49\pi/100)<0, \qquad
\textrm{for}\ \ \theta \in [49\pi/100, \pi/2], \ \ \textrm{and}\ \
u<0.
\]
From the above inequality we get the estimates $\vw_2^2\geq
v(\theta)^2$, for all $\theta \in [49\pi/100, \pi/2]$ and $u<0$.
Using this estimate  we can write:
\[
-\dfrac{dv}{d\theta}=\dfrac12
\sqrt{2^{1-\alpha}\left[1+\dfrac{1}{\sin^\alpha\theta}+\dfrac{1}{\cos^\alpha\theta}\right]-v^2}
\geq\dfrac12\sqrt{2^{1-\alpha}\left[\dfrac{1}{\cos\theta}+\va_7\right]
-\vw_2^2}
\]
where $\va_7:=\min_{\theta \in [49\pi/100,
\pi/2]}h(\theta)=h(\pi/2)=2$. By a direct integration we have:
\begin{equation*}
\begin{aligned}
v(49\pi/100) &\geq v(\pi/2) +
\dfrac12\int_{49\pi/100}^{\pi/2}\sqrt{2^{1-\alpha}\left[\dfrac{1}{\cos^\alpha\theta}+\va_7\right]
-\vw_2^2}\ d\theta \geq\\
&\geq v(\pi/2) + \dfrac12 \int_0^{\pi/100}
\sqrt{2^{1-\alpha}\left[\dfrac{1}{\Phi^\alpha}+2\right]-\vw_2^2}\
d\Phi\geq\\ &\geq \vw_2 + \dfrac12 \int_0^{\pi/100}
\sqrt{2^{1-\alpha}\left[\dfrac{1}{\Phi^\alpha}+2\right]-\vw_2^2}\
d\Phi\geq \\
&\geq%\vw_2 + \dfrac12 \int_0^{\pi/100}
%\sqrt{\left[\dfrac{1}{\Phi}+2\right]-\vw_2^2}\ d\Phi \approx
{\bf-1.1804}.
\end{aligned}
\end{equation*}
\begin{itemize}
\item Now we prove that the continuation of the right
branch through double collision is negative at the point
$\theta=\pi/4$.
\end{itemize}
In fact, for $u<0$
\[
-\dfrac{dv}{d\theta} \leq
\dfrac12\sqrt{2^{1-\alpha}\left[\dfrac{1}{\cos^\alpha\theta}+\vb_5\right]}
\]
for $\vb_5:=h(\pi/4)$ and $h$ as above. Let $\vb_6:=g(\pi/4)$.
Integrating we get:
\begin{equation}
\begin{aligned}
v(\pi/4)&\leq v(\pi/2)
+\dfrac12\int_{\pi/4}^{\pi/2}\dfrac{1}{2^{1/2+\beta}}\sqrt{\left[\dfrac{1}{\cos^\alpha\theta}+\vb_5\right]}\ d\theta \leq\\
&\leq
v(\pi/2)+\dfrac12\int_0^{\pi/4}\dfrac{1}{2^{1/2+\beta}}\sqrt{\left[\dfrac{\vb_6}{\Phi^{2\beta}}+\vb_5\right]}\
d\Phi \leq\\
\leq& \vw_1
+\dfrac12\int_0^{\pi/}\dfrac{1}{2^{1/2+\beta}}\sqrt{\left[\dfrac{\vb_6}{\Phi^{2\beta}}+\vb_5\right]}\
d\Phi\leq\\
&\leq \vk_1(\beta)+
\dfrac12\int_0^{\pi/4}\dfrac{1}{2^{1/2+\beta}}\sqrt{\left[\dfrac{\vb_6}{\Phi^{2\beta}}+\vb_5\right]}\
d\Phi\leq\\&\leq -0.6857 +
\dfrac12\int_0^{\pi/4}\dfrac{1}{2^{1/2+\beta}}\sqrt{\left[\dfrac{\vb_6}{\Phi^{2\beta}}+\vb_5\right]}\
d\Phi\leq\\
%v_1
%+\dfrac{1}{\sqrt2}\int_0^{\pi/16}\sqrt{\dfrac{\va_6}{\theta}+\va_5}\
%d\theta \\
&\leq {\bf-0.0713}.
\end{aligned}
\end{equation}

%
%
%For this we consider $\theta \in [\pi/4,49\pi/100]$ and $u<0$. We
%observe that in $\theta=49\pi/100$, $U(\theta)$ has a relative
%maximum and therefore for $\theta \in [\pi/4,49\pi/100 ]$ we have
%$U(\theta)\leq U(49\pi/100)$. Now let $\bar
%v:=-\sqrt{2U(49\pi/100)}$; thus
%\[
%-\dfrac{dv}{d\theta}\leq \dfrac12\sqrt{\bar v^2-v^2}
%\]
%and by a direct integration, we have
%\[
%-\arcsin\left(\dfrac{v}{\bar
%v}\right)\Big\vert^{v(49\pi/100)}_{v(\pi/4)} \leq \dfrac12
%\int_{\pi/4}^{49\pi/100} d \theta=3\pi/25.
%\]
%Therefore
%\[
%v(\pi/4)\leq \bar v\sin(3\pi/25+\arcsin(v(49\pi/100)/\bar v))\leq
%\textcolor{blue}{-0.4502}.
%\]
%Moreover  in $\theta=\pi/4$, $U(\theta)$ has a minimum and therefore
%for $\theta \in [\pi/8,3\pi/8 ]$ we have $U(\theta)\geq U(\pi/4)$.
%Now let $\bar v':=\sqrt{2U(\pi/4)}$; thus
%\[
%-\dfrac{dv}{d\theta}\geq \dfrac12\sqrt{\bar v'^2-v^2}
%\]
%and by a direct integration, we have
%\[
%-\arcsin\left(\dfrac{v}{\bar
%v}\right)\Big\vert^{v(3\pi/8)}_{v(\pi/8)} \geq \dfrac12
%\int_{\pi/8}^{3\pi/8} d \theta=\pi/8.
%\]
%Therefore
%\[
%v(\pi/8)\geq \bar v'\sin(\pi/8+\arcsin(v(3\pi/8)/\bar v'))\geq {\bf
%-0.3741}.
%\]
%These two estimates in particular imply that $v(\pi/4)<0$ by taking
%into account that for $u<0$ the solution $v$ is strictly decreasing
%on the interval $[\pi/8, 3\pi/8]$.

To finish the estimate for the continuation of the right branch
through double collision, we now consider $\theta \in [0, \pi/8]$
and $u<0$. We observe that in $\theta=\pi/8$, $U(\theta)$ has a
relative minimum and therefore for $\theta \in [0, \pi/8]$ we have
$U(\theta)\geq U(\pi/8)$. Now let $\tilde v:=\sqrt{2U(\pi/8)}$; thus
\[
-\dfrac{dv}{d\theta}\geq \dfrac12\sqrt{\tilde v^2-v^2}
\]
and by a direct integration, we have
\[
-\arcsin\left(\dfrac{v}{\tilde v}\right)\Big\vert^{v(\pi/8)}_{v(0)}
\geq \dfrac12 \int_{0}^{\pi/8} d \theta=\pi/16.
\]
Therefore
\[
v(0)\geq \tilde v\sin\left(\pi/16+
\arcsin\left(\dfrac{v(\pi/8)}{\tilde v}\right)\right)\geq \tilde
v\sin\left(\pi/16+ \arcsin\left(\dfrac{-0.3741}{\tilde
v}\right)\right) \approx {\bf0.4948}.
\]
By symmetry, the same arguments holds also for the left branch of
$W^u(\vc_{-})$. Now the thesis easily follows by the fact that the
flow is of gradient-type with restpoints $\vc_{-}$ and $\vc_{+}$.
\begin{enumerate}
\item[(1c)] {\bf There exists $\alpha_{01} \in (\alpha_{00}, \alpha_0)$
such that the continuation of the right branch through double
collision at $\theta=\pi/2$ for $u<0$, changes sign in the interval
$[\pi/4, \pi/2]$ for any $\alpha \in [\alpha_{01}, \alpha_0)$.\/}
\end{enumerate}
In fact, using this estimate  we can write:
\[
-\dfrac{dv}{d\theta}=\dfrac12
\sqrt{2^{1-\alpha}\left[1+\dfrac{1}{\sin^\alpha\theta}+\dfrac{1}{\cos^\alpha\theta}\right]-v^2}
\geq\dfrac12\sqrt{2^{1-\alpha}\left[\dfrac{1}{\cos\theta}+\va_7\right]
-\vw_2^2}
\]
where $\va_7:=\min_{\theta \in [49\pi/100,
\pi/2]}h(\theta)=h(\pi/2)=2$. By a direct integration we have:
\begin{equation*}
\begin{aligned}
v(49\pi/100) &\geq v(\pi/2) +
\dfrac12\int_{49\pi/100}^{\pi/2}\sqrt{2^{1-\alpha}\left[\dfrac{1}{\cos^\alpha\theta}+\va_7\right]
-\vw_2^2}\ d\theta \geq\\
&\geq v(\pi/2) + \dfrac12 \int_0^{\pi/100}
\sqrt{2^{1-\alpha}\left[\dfrac{1}{\Phi^\alpha}+2\right]-\vw_2^2}\
d\Phi\geq\\ &\geq \vw_2 + \dfrac12 \int_0^{\pi/100}
\sqrt{2^{1-\alpha}\left[\dfrac{1}{\Phi^\alpha}+2\right]-\vw_2^2}\
d\Phi\geq \\
&\geq%\vw_2 + \dfrac12 \int_0^{\pi/100}
%\sqrt{\left[\dfrac{1}{\Phi}+2\right]-\vw_2^2}\ d\Phi \approx
{\bf-1.1804}.
\end{aligned}
\end{equation*}
To finish the estimate for the continuation of the right branch
through double collision, we now consider $\theta \in [3\pi/8,
\pi/2]$ and $u<0$. We observe that in $\theta=\pi/4$, $U(\theta)$
has a relative minimum and therefore for $\theta \in [\pi/4,
49\pi/100]$ we have $U(\theta)\geq U(\pi/4)$. Now let $\tilde
v:=\sqrt{2U(\pi/4)}$; thus
\[
-\dfrac{dv}{d\theta}\geq \dfrac12\sqrt{\tilde v^2-v^2}
\]
and by a direct integration, we have
\[
-\arcsin\left(\dfrac{v}{\tilde
v}\right)\Big\vert^{v(49\pi/100)}_{v(\pi/4)} \geq \dfrac12
\int_{\pi/4}^{49\pi/100} d \theta=3\pi/25.
\]
Therefore
\[
v(\pi/4)\geq \tilde v\sin\left(3\pi/25+
\arcsin\left(\dfrac{v(3\pi/25)}{\tilde v}\right)\right)\geq \tilde
v\sin\left(3\pi/25+ \arcsin\left(\dfrac{-1.1804}{\tilde
v}\right)\right) \approx {\bf0.0511}.
\]

By symmetry, the same arguments holds also for the left branch of
$W^u(\vc_{-})$. Now the thesis easily follows by the fact that the
flow is of gradient-type with restpoints $\vc_{-}$ and $\vc_{+}$.
\begin{enumerate}
\item[(1d)] {\bf There exists $\alpha_{0*} \in (\alpha_{00}, \alpha_{01})$
such that for the continuation of the right branch through double
collision at $\theta=\pi/2$ for $u<0$, we have $v(\pi/4)=0$.\/}
\end{enumerate}
In order to prove this fact we consider the one parameter family of
initial value problems:
\begin{equation}
(P_\alpha) \qquad \left\{
\begin{array}{ll}
\dfrac{dv}{d\theta}= \dfrac12 \sqrt{2^{1-\alpha}\left[1+
\dfrac{1}{\cos^\alpha \theta}+ \dfrac{1}{\sin^\alpha
\theta}\right]-v^2}:=\vF_\alpha(\theta, v)\\
\\
v(\pi/4)=-\sqrt{2^{1-\alpha}(1+2^{\beta+1})},
\end{array}\right.
\end{equation}
which can be written in integral form as follows
\[
v_\alpha(\theta)=v(\pi/4)+\int_{\pi/4}^\theta \vF_\alpha(\bar
\theta, v(\bar\theta))\, d\bar\theta, \qquad \forall\, \theta \in
(\pi/4, \pi/2).
\]
If $v$ is a classical ($\mathscr C^1$)-solution of the ivp above,
then it is a $\mathscr C^0$-solution of the corresponding Volterra
integral equation; moreover the function
\[
\vI:[1,\alpha_0)\longrightarrow  \R: \alpha \longmapsto
\vI(\alpha):=\int_{\pi/4}^{\pi/2} \vF_\alpha(\theta, v(\theta))\,
d\theta,
\]
is continuous, by taking into account the theorem on integrals
depending on parameters. Thus the function
\[
\vV:(1,\alpha_0)\rightarrow \R:\alpha\mapsto
\vV(\alpha):=v_\alpha(\pi/2):= v(\pi/4)+ \vI(\alpha)
\]
is continuous. Moreover by the monotonicity of the integral of a
nonnegative function with respect to the domain of integration and
by the fact that the map $\alpha\mapsto \vF_\alpha(\theta, v)$
increases on $(1,2)$, it follows that the map $\vV$ also increases.
Furthermore $\vV(\alpha_{00})<0$ and $\vV(\alpha_{01})>0$; thus by
taking into account the theorem of zeros of a continuous and
increasing function, there exists a unique $\alpha_*$ in between
such that $V(\alpha_{0*})=0$, which means nothing but that
$v_\alpha(\pi/4)=0$.

\smallskip

In order to prove item 2, we proceed as follows.
\begin{enumerate}
\item[(2)] {\bf There exists $\alpha_1 \in (1,2)$ such that for
each $\alpha \in (\alpha_1, 2)$, $v(\pi/2)>0$.\/}
\end{enumerate}
\begin{itemize}
\item Lower bound for $v(\pi/2)$ in the case $u>0$ and $\theta \in (\pi/4,
\pi/2)$.
\end{itemize}

Next we calculate a lower estimate for $v(\pi/2)$. For this we first
observe that
\[
\dfrac12
\sqrt{2^{1-\alpha}\left[1+\dfrac{1}{\cos^\alpha\theta}+\dfrac{1}{\sin^\alpha\theta}\right]-v^2}\geq
\dfrac{1}{2}\sqrt{2^{1-\alpha}\left[\dfrac{1}{\cos^\alpha\theta}+\vb_3\right]-v^2},
\]
where $\vb_3:= \min_{\theta \in [\pi/4, \pi/2]}h(\theta)=2$. Now we
consider the interval $J:=[3\pi/8, \pi/2]$ and by a direct
calculation, we have:
\begin{equation*}
\begin{aligned}
\dfrac{dv}{d\theta}&\geq
\dfrac{1}{2}\sqrt{2^{1-\alpha}\left[\dfrac{1}{\cos^\alpha\theta}+2\right]-v^2(\pi/4)}\
\Rightarrow\\ v(\pi/2) &\geq v(3\pi/8) +
\int_{3\pi/8}^{\pi/2}\dfrac{1}{2^{1/2+\beta}}\sqrt{\dfrac{1}{\cos^\alpha\theta}+2-
2^{\alpha-1}v^2(\pi/4)}\ d\theta\geq\\
&\geq v(3\pi/8) +
\int_{0}^{\pi/8}\dfrac{1}{2^{1/2+\beta}}\sqrt{\dfrac{1}{\Phi^{2\beta}}+1-2^{\beta+1}}\
d\Phi:=\vk_2(\beta).
\end{aligned}
\end{equation*}
We observe that:
\[
v(3\pi/8):=v_\beta(3\pi/8)=-\sqrt{2^{1-2\beta}\left[1+
\left(\dfrac{2\sqrt2}{1+\sqrt2}\right)^\beta+\left(\dfrac{2\sqrt2}{\sqrt2-1}\right)^\beta
\right]}.
\]
Moreover, the function below
\[
(1/2,1)\ni\beta\mapsto v_\beta(3\pi/8)\ \textrm{and}\ (0, 1/2)\ni
\beta\mapsto
\int_{0}^{\pi/8}\dfrac{1}{2^{1/2+\beta}}\sqrt{\dfrac{1}{\Phi^{2\beta}}+1-2^{\beta+1}}
\]
as well as $\vk_2$ are strictly increasing functions. %it immediately
%follows that also $\vk_2$ is an increasing function and therefore
%\[
%v_{1/2}(3\pi/8)\approx \textcolor{red}{-2.449}\leq v(3\pi/8)\leq
%v_0(3\pi/8) \approx \textcolor{red}{-2.167}.
%\]
%To conclude the estimate, we observe that also the function:
%\[
%(0,1/2)\ni\beta\mapsto
%\vi(\beta):=\int_{0}^{\pi/8}\dfrac{1}{2^{1/2+\beta}}\sqrt{\dfrac{1}{\theta^{2\beta}}+1-
%2^{(\beta+1)}}\ d\theta\in \R^+
%\]
%is a positive and increasing function, with $\vi(0)=0$. Therefore,
However it changes sign and in particular it is positive in the
interval  $A_2:=[1.7,2)$. Therefore by choosing any number $\alpha_1
\in A_2$ we can  conclude that
\[
v(\pi/2)\geq \vk_2(1.7):=\vw_2\approx {\bf 0.1055}>0.
\]

By symmetry, the same arguments holds also for the left branch of
$W^u(\vc_{-})$. Now the thesis easily follows by the fact that the
flow is of gradient-type with restpoints $\vc_{-}$ and $\vc_{+}$.

\smallskip
\begin{enumerate}
\item[(3)] {\bf There exists a unique $\alpha_* \in (\alpha_0,\alpha_1)$ such that
$v_\alpha(\pi/2)=0$.\/}
\end{enumerate}
In order to prove this fact we consider the one parameter family of
initial value problems:
\begin{equation}
(P_\alpha) \qquad \left\{
\begin{array}{ll}
\dfrac{dv}{d\theta}= \dfrac12 \sqrt{2^{1-\alpha}\left[1+
\dfrac{1}{\cos^\alpha \theta}+ \dfrac{1}{\sin^\alpha
\theta}\right]-v^2}:=\vF_\alpha(\theta, v)\\
\\
v(\pi/4)=-\sqrt{2^{1-\alpha}(1+2^{\beta+1})},
\end{array}\right.
\end{equation}
which can be written in integral form as follows
\[
v_\alpha(\theta)=v(\pi/4)+\int_{\pi/4}^\theta \vF_\alpha(\bar
\theta, v(\bar\theta))\, d\bar\theta, \qquad \forall\, \theta \in
(\pi/4, \pi/2).
\]
If $v$ is a classical ($\mathscr C^1$)-solution of the ivp above,
then it is a $\mathscr C^0$-solution of the corresponding Volterra
integral equation; moreover the function
\[
\vI:[1,2)\longrightarrow  \R: \alpha \longmapsto
\vI(\alpha):=\int_{\pi/4}^{\pi/2} \vF_\alpha(\theta, v(\theta))\,
d\theta,
\]
is continuous, by taking into account the theorem on integrals
depending on parameters. Thus the function
\[
\vV:(1,2)\rightarrow \R:\alpha\mapsto \vV(\alpha):=v_\alpha(\pi/2):=
v(\pi/4)+ \vI(\alpha)
\]
is continuous. Moreover by the monotonicity of the integral of a
nonnegative function with respect to the domain of integration and
by the fact that the map $\alpha\mapsto \vF_\alpha(\theta, v)$
increases on $(1,2)$, it follows that the map $\vV$ also increases.
Furthermore $\vV(\alpha_0)<0$ and $\vV(\alpha_1)>0$; thus by taking
into account the theorem of zeros of a continuous and increasing
function, there exists a unique $\alpha_*$ in between such that
$V(\alpha_*)=0$, which means nothing but that $v_\alpha(\pi/2)=0$.
\finedim

%%%%%===============================================================================================
\section{Some pointwise estimates for the tetrahedral case}\label{app:tetra}
\subsection{Proof of Proposition \ref{prop:flussotetra}}

The proof of this result will be divided into two main steps. The
first correspond to the Newtonian case $\alpha=1$ and the second to
the case $\alpha \in (0,1)$.

Moreover we refer to the {\em right (left) unstable branch of
$W^u(\ve^-_{11})$ as that having $u>0$ ($u<0$) in a small
neighborhood of $\ve^-_{11}$ in $\parabolicmanifold$.\/}

Using uniqueness of solutions, the symmetries and the gradient-like
character of the flow with respect to $v$, it is sufficient to prove
the following fact.
\begin{enumerate}
{\bf \item[(1)] The right branch of $W^u(\ve^{-}_{11})$ intersects
$\varphi=\pi/2$ ($\varphi=-\pi/2$, resp.) with $v<0$ and then
intersects the section $v=0$ with an angle $\varphi >0$ (resp.
$\varphi<0$).\/}
\end{enumerate}
\subsubsection*{Proof in the Newtonian case.} The proof of this claim
will be divided into some parts:
\begin{enumerate}
\item[(1a)] Upper bound for $v(\pi/2)$ in the case $u>0$ and $\varphi \in (\varphi_1,
\pi/2)$, where $\varphi_1=\arctan(1/\sqrt2)$.
\begin{equation*}
\begin{aligned}
\dfrac{dv}{d\theta}&=\dfrac12
\sqrt{\dfrac{1}{\cos\varphi}+\dfrac{2\sqrt2}{\sqrt{1+\sin^2\varphi}}-v^2}\leq\\
&\leq \dfrac12\sqrt{\dfrac{1}{\cos \varphi}+\vd_1}
\end{aligned}
\end{equation*}
where
\[
\vd_1:=\max_{\varphi \in [\varphi_1, \pi/2]}\vm(\varphi)=\sqrt 6,
\qquad \vm(\varphi):= \dfrac{2\sqrt2}{\sqrt{1+\sin^2\varphi}}.
\]
We observe that the maximum is achieved at the point $\varphi_1$
since the function $\vm$ is strictly decreasing on that interval; in
fact
\[
\vm'(\varphi)= -\dfrac{2\sqrt2\sin \varphi
\cos\varphi}{(1+\sin^2\varphi)^{3/2}}.
\]

Thus
\begin{equation*}
\begin{aligned}
v(\pi/2)&\leq v(\varphi_1)+ \dfrac12\int_{\varphi_1}^{\pi/2}
\sqrt{\dfrac{1}{\cos\varphi}+\vd_1}\ d\varphi  =_{\Phi:=\pi/2-\varphi}\\
&=v(\varphi_1) +  \dfrac12\int_0^{\pi/2-\varphi_1}
\sqrt{\dfrac{1}{\sin\Phi}+\vd_1}\ d\Phi
\leq\\
&\leq v(\varphi_1) +\dfrac12 \int_0^{\pi/2-\varphi_1}
\sqrt{\dfrac{\vd_2}{\Phi}+\vd_1} \
d\Phi=\dfrac12\left[-\sqrt{6\sqrt6}+\int_0^{\pi/2-\varphi_1}
\sqrt{\dfrac{\vd_2}{\Phi}+\vd_1} \ d\Phi\right]
\end{aligned}
\end{equation*}
where
$\vd_2:=\max_{[0,\pi/2-\varphi_1]}\vn(\Phi)=\dfrac{\sqrt6}{2}\left(\dfrac{\pi}{2}-\arctan\dfrac{1}{\sqrt
2}\right)$, since the function $\vn(\Phi):=\dfrac{\Phi}{\sin\Phi}$
is increasing on $[0, \pi/2-\varphi_1]$. Since
$v(\varphi_1)=-\dfrac{\sqrt{6\sqrt6}}{2}$ and by using the
quadrature
\[
\int \sqrt{\dfrac{a}{x}+b}\, dx= \sqrt{x(a+bx)}+ \dfrac{a}{2\sqrt
b}\log[2\sqrt b\sqrt{(x(a+bx)}+ 2bx+a],
\]
it follows that
\begin{equation*}
\begin{aligned}
v(\pi/2)\leq    \vv_1:= {\bf -0.5727\/}<0,
\end{aligned}
\end{equation*}
therefore $v$ intersects $\theta=\pi/2$ at a negative value.
%%%%%%%%%%%%%%%%%%%%%5%%%%%%%%%%%%%%%%%%%%%%%%%%%%%%%%%%%%%%%%%%%%%%%%%%%%%%%%%%%%%%%%%%%%%%%%%%%%%%%%%%%%%%%%%
\item[(1b)] Lower bound for $v(\pi/2)$ in the case $u>0$ and $\varphi \in (\varphi_1,
\pi/2)$.\\
Next we calculate a lower estimate for $v(\pi/2)$. For this we first
observe that
\[
\dfrac{dv}{d\varphi}=\dfrac12
\sqrt{\dfrac{1}{\cos\varphi}+\dfrac{2\sqrt2}{\sqrt{1+\sin^2\varphi}}-v^2}\geq
\dfrac12\sqrt{\dfrac{1}{\cos\varphi}+\vd_3-v^2},
\]
where $\vd_3:= \min_{\phi \in [\varphi_1\, \pi/2]}\vm(\varphi)=2$.
Now we consider the interval $J:=[7/16\pi, \pi/2]\subset[\varphi_1,
\pi/2]$ and by a direct calculation:
\begin{equation*}
\begin{aligned}
\dfrac{dv}{d\varphi}\geq \dfrac12\sqrt{\dfrac{1}{\cos\varphi}+2-
v^2(\varphi_1)}\ \Rightarrow v(\pi/2) &\geq v(7/16\pi) +
\dfrac12\int_{7/16\pi}^{\pi/2}\sqrt{\dfrac{1}{\cos\varphi}+2-
v^2(\varphi_1)}\ d\varphi\geq\\
& \geq v(\varphi_1) +
\dfrac12\int_{0}^{\pi/16}\sqrt{\dfrac{1}{\Phi}+\vd_4}\ d\Phi\\
&\geq {\bf-1.4994}:=\vv_2<0,
\end{aligned}
\end{equation*}
for $\Phi$ as above and $\vd_4:=2-v^2(\varphi_1)=2-3\sqrt6/2$;
moreover in the last inequality we used the fact that $\sin \Phi
\leq \Phi$. From (1a) \&(1b) it follows that $\vv_2 \leq
v(\theta)\leq \vv_1$, for $\varphi \in [\varphi_1, \pi/2]$ and
$u>0$.
%%%%%%%%%%%%%%%%%%%%%%%%%%%%%%%%%%%%%%%%%%%%%%%%%%%%%%%%%%%%%%%%%%%%%%%%%%%%%%%%%%%
\item[(1c)] The continuation of the right branch through double
collision at $\theta=\pi/2$ for $u<0$. We affirm that $v$ remains
negative in the interval $J$.

In fact, for $u<0$
\[
-\dfrac{dv}{d\theta} \leq \dfrac12\sqrt{\dfrac{1}{\cos\theta}+\vd_5}
\]
for $\va_5:=\vm(7\pi/16)\approx {2.0193}$. Let
$\vd_6:=\vn(\pi/16)\approx {1.0064}$. Integrating we get:
\begin{equation}
\begin{aligned}
v(3/8\pi)&\leq v(\pi/2)
+\dfrac12\int_{7\pi/16}^{\pi/2}\sqrt{\dfrac{1}{\cos\varphi}+\vd_5}\ d\varphi \leq\\
&\leq
v(\pi/2)+\dfrac12\int_0^{\pi/16}\sqrt{\dfrac{\vd_6}{\Phi}+\vd_5}\
d\Phi \leq \vv_1
+\dfrac12\int_0^{\pi/16}\sqrt{\dfrac{\vd_6}{\Phi}+\vd_5}\
d\Phi \\
&\leq \vv_1+ {0.4723}\leq {\bf-0.1004}<0
\end{aligned}
\end{equation}
Therefore
\[
\vv_2 \leq v(\pi/2)\leq v(\theta)\leq v(7\pi/16)<0, \qquad
\textrm{for}\ \ \theta \in [7\pi/16, \pi/2], \ \ \textrm{and}\ \
u<0.
\]
Therefore
\[
\vv_2 \leq v(\pi/2)\leq v(\theta)\leq v(7\pi/16)<0, \qquad
\textrm{for}\ \ \theta \in [7\pi/16, \pi/2], \ \ \textrm{and}\ \
u<0.
\]

%%%%%%%%%%%%%%%%%%%%%%%%%%%%%%%%%%%%%%%%%%%%%%%%%%%%%%%%%%%%%%%%%%%%%%%%%%%%%%%%%%%
\item[(1d)] From the above inequality we get the estimates $\vv_2^2\geq
v(\theta)^2$, for all $\theta \in [7\pi/16, \pi/2]$ and $u<0$. Using
this estimate  we can write:
\[
-\dfrac{dv}{d\varphi}=\dfrac12
\sqrt{2U(\varphi)-v^2}\geq\dfrac12\sqrt{\dfrac{1}{\cos\varphi}+\vd_7
-\vv_2^2}
\]
where $\va_7:=\min_{\varphi \in [7\pi/16,
\pi/2]}\vm(\varphi)=f(\pi/2)=2$. By a direct integration we have:
\begin{equation*}
\begin{aligned}
v(7\pi/16) &\geq v(\pi/2) +
\dfrac12\int_{7\pi/16}^{\pi/2}\sqrt{\dfrac{1}{\cos\varphi}+\vd_7
-\vv_2^2}\ d\theta \geq\\
&\geq v(\pi/2) + \dfrac12 \int_0^{\pi/16}
\sqrt{\dfrac{1}{\Phi}+2-\vv_2^2}\
d\Phi\\
&\geq \vv_2+ \dfrac12 \int_0^{\pi/16} \sqrt{\dfrac{1}{\Phi}+\vd_8}\
d\Phi\geq\vv_3:={\bf-1.0600},
\end{aligned}
\end{equation*}
where $\vd_8:=2-\vv_2^2$.
\item[(1e)]To finish the estimate for the continuation of the right
branch through double collision, we now consider $\varphi \in
[0,7\pi/16]$ and $u<0$. We observe that in $\varphi=\varphi_1$,
$U(\theta)$ has a relative minimum and therefore for $\varphi \in
[0,7\pi/16]$ we have $U(\varphi)\geq U(\varphi_1)$. Now let $\tilde
v:=\sqrt{2U(\varphi_1)}$; thus
\[
-\dfrac{dv}{d\varphi}\geq \dfrac12\sqrt{\bar v^2-v^2}
\]
and by a direct integration, we have
\[
-\arcsin\left(\dfrac{v}{\bar v}\right)\Big\vert_{v(0)}^{v(7\pi/16)}
\geq \dfrac12 \int_{0}^{7\pi/16} d \theta=7\pi/32.
\]
Therefore
\[
v(0)\geq \tilde v\sin\left(7\pi/32+\arcsin(v(7\pi/16)/\bar
v)\right)\geq \bar v\sin\left(7\pi/32+\arcsin(\vv_3/\bar v)\right)
\geq {\bf0.1937}.
\]
\end{enumerate}
By symmetry, the same arguments holds also for the left branch of
$W^u(\vc_{-})$. Now the thesis easily follows by the fact that the
flow is of gradient-type with restpoints $\vc_{-}$ and $\vc_{+}$.
\begin{enumerate}
{\bf \item[(2)] The left branch of $W^u(\ve^1_{11})$ intersects
$\varphi=-\pi/2$ with at a negative value of $v$.\/}
\end{enumerate}
In order to prove this claim we follow the idea used by \cite{DV99}.

Now we consider the left branch of $W^u(\ve^-_{11})$. It is clear
that for the same value of $\varphi$, $v$ is greater along the left
branch of $W^u(\ve^-_{11})$ that  the corresponding value for the
left branch of $W^u(\ve_{12}^-)$, then due to the symmetries of the
problem it follows that the left branch of $W^u(\ve_{11}^-)$ meets
the section $v=0$ for some value of $\varphi$ which is negative. It
remains to prove that it actually meets $\varphi=-\pi/2$ at a
negative value of $v$. For $u<0$, we have:
\[
-\dfrac{dv}{d\varphi}=\sqrt{\va(\varphi)-\dfrac{v^2}{4}}
\]
where
$\va(\varphi):=\dfrac{U_\antiprisma(\varphi)}{2}=\dfrac{1}{4\cos\varphi}+\dfrac{\sqrt2}{2\sqrt{1+\sin^2\varphi}}$.
If $I_n:=[\tilde \varphi_{n-1}, \tilde \varphi_n]$ and
$\va_n:=\max_{\varphi \in I_n}\va(\varphi)$, then for $\varphi \in
I_n$, we have:
\[
-\dfrac{dv}{d\varphi}\leq\sqrt{\va_n-\dfrac{v^2}{4}}.
\]
By a direct integration, we can conclude that
\begin{equation}\label{eq:disricorsiva}
v(\tilde \varphi_n)\leq 2\sqrt{\va_n}\sin\left(\dfrac12|\tilde
\varphi_n-\tilde \varphi_{n-1}|+\arcsin\left(\dfrac{v(\tilde
\varphi_{n-1})}{2\sqrt{\va_n}}\right)\right).
\end{equation}
It is readily verified that $\va_n=\dfrac{v(\tilde \varphi_n)}{2}$
for any interval $I_n\subset [-\pi/2, -\varphi]\cup [0,
\varphi_1]:=\vA$, since $\varphi\mapsto \va(\varphi)$ is positive
and non-increasing on that subset; moreover $\va_n=\va(\tilde
\varphi_{n-1})$ in $[\pi/2, \pi/2]\setminus \vA:=\vB$ since on this
set this function is positive and non-decreasing. We define
\[
\tilde \varphi_n= \dfrac{(12-n)\varphi_1}{12}, \qquad \textrm{for
each}\ \ n\, \in\, \n{36}.
\]
Solving recursively the inequality, we obtain
\[
v(\tilde \varphi_{32})=v(-5\varphi_1/3)\leq{\bf -1.1452.\/}
\]
On the other side, for $ \varphi\in [-\pi/2, 5\varphi_1/3]$,
\[
-\dfrac{dv}{d\varphi} \leq
\dfrac12\sqrt{\dfrac{1}{\cos\varphi}+\vd_9}
\]
where  $\vd_9:=\vm (-5\varphi_1/3)$. Again integrating and letting
$\vd_{10}:= \vn(\pi/2-5\varphi_1/3)$, we obtain:
\begin{equation}\label{eq:quasifine}
\begin{aligned}
-v(5\varphi_1/3)+ v(-\pi/2)&\leq
\dfrac12\int_{-\pi/2}^{-5\varphi_1/3}
\sqrt{\dfrac{1}{\cos\varphi}+\vd_9}\,d\varphi\\
&\leq \dfrac12\int_0^{\pi/2-5\varphi_1/3}
\sqrt{\dfrac{1}{\sin\Phi}+\vd_9}\,d\Phi\\
&\leq \dfrac12\int_0^{\pi/2-5\varphi_1/3} \sqrt{\dfrac{\vd_{10}}{\Phi}+\vd_9}\,d\Phi\\
&\leq {\bf 0.8803}
\end{aligned}
\end{equation}
Now suppose that $v(-\pi/2)>0$; then from the inequality
\eqref{eq:quasifine} it  follows:
\[
v(-5\varphi_1/3)\geq {\bf -0.8803}
\]
which  is  a contradiction. We conclude that $v(-\pi/2)<0$ for the
left branch of $W^u(\ve_{11}^-)$.

\finedim

\subsubsection*{Proof in the  case $\alpha \in (0,1)$}

The proof of this claim will be divided into some parts:
\begin{enumerate}
\item[(1a)] Upper bound for $v(\pi/2)$ in the case $u>0$ and $\varphi \in (\varphi_1,
\pi/2)$, where $\varphi_1=\arctan(1/\sqrt2)$.
\begin{equation*}
\begin{aligned}
\dfrac{dv}{d\theta}&=\dfrac12
\sqrt{\dfrac{1}{(2\cos\varphi)^\alpha}+
\dfrac{2^{1-\alpha/2}}{(1+\sin^2
\varphi)^{\alpha/2}}-v^2}\leq\dfrac12
\sqrt{\dfrac{1}{(2\cos\varphi)^\alpha}+ \dfrac{2^{1-\alpha/2}}{(1+\sin^2 \varphi)^{\alpha/2}}}\\
&=\dfrac{1}{2^{1/2+\beta}}\sqrt{\dfrac{1}{\cos
\varphi^{2\beta}}+\dfrac{2^{1+\beta}}{(1+\sin^2\varphi)^\beta}} \leq
\dfrac{1}{2^{1/2+\beta}}\sqrt{\dfrac{1}{\cos\varphi^{2\beta}}+\ve_1}
\end{aligned}
\end{equation*}
where
\[
\ve_1:=\max_{\varphi \in [\varphi_1, \pi/2]}\vm(\beta,
\varphi)=2^{1-\beta}3^\beta \qquad \vm(\beta, \varphi):=
\dfrac{2^{1+\beta}}{(1+\sin^2\varphi)^\beta}.
\]
In fact for each $\beta \in (0, 1/2)$, the map $\varphi\mapsto
m(\beta, \varphi)$ is strictly decreasing on the interval
$[\varphi_1, \pi/2]$; thus the maximum is achieved at the point
$\varphi=\varphi_1$. Moreover
\[
\dfrac{\partial \vm}{\partial
\beta}=\dfrac{2^{1+\beta}}{(1+\sin^2)^\beta}[\log 2-
\log(1+\sin^2\varphi)] \qquad \dfrac{\partial \vm}{\partial
\varphi}=-\dfrac{2^{2+\beta}\beta\sin\varphi\cos\varphi}{(1+\sin^2\varphi)^{\beta+1}}
\]
Thus
\begin{equation*}
\begin{aligned}
v(\pi/2)&\leq v(\varphi_1)+
\dfrac{1}{2^{\beta+1/2}}\int_{\varphi_1}^{\pi/2}
\sqrt{\dfrac{1}{\cos\varphi^{2\beta}}+\ve_1}\ d\varphi  =_{\Phi:=\pi/2-\varphi}\\
&=v(\varphi_1) +  \dfrac{1}{2^{\beta+1/2}}\int_0^{\pi/2-\varphi_1}
\sqrt{\dfrac{1}{\sin\Phi^{2\beta}}+\ve_1}\ d\Phi
\leq\\
&\leq v(\varphi_1) +\dfrac{1}{2^{\beta+1/2}}
\int_0^{\pi/2-\varphi_1}
\sqrt{\dfrac{\ve_2}{\Phi^{2\beta}}+\ve_1} \ d\Phi=\\
&= -\sqrt 6\sqrt{3^\beta 8^{-\beta}}+\dfrac{1}{2^{\beta+1/2}}
\int_0^{\pi/2-\varphi_1} \sqrt{\dfrac{\ve_2}{\Phi^{2\beta}}+\ve_1} \
d\Phi:=\vj_1(\beta)
\end{aligned}
\end{equation*}
where $\ve_2:=\max_{[0,\pi/2-\varphi_1]}\vn(\beta,\Phi)$ where
$\vn(\beta, \Phi)= \left(\dfrac{\Phi}{\sin\Phi}\right)^{2\beta}$. We
observe that
\[
\dfrac{\partial \vn}{\partial \beta}=
2\left(\dfrac{\Phi}{\sin\Phi}\right)^{2\beta}\log
\dfrac{\Phi}{\sin\Phi}, \qquad \dfrac{\partial \vn}{\partial \Phi}=
\dfrac{2\beta\sin\Phi}{\Phi}\left[\left(\dfrac{\Phi}{\sin\Phi}\right)^{2\beta}\left(\dfrac{1}{\sin\Phi}-
\dfrac{\Phi\cos\Phi}{\sin^2\Phi}\right) \right].
\]
Now since the function $[0, \pi/2-\varphi1]\ni\Phi\mapsto
\vg(\Phi):=\dfrac{\Phi}{\sin\Phi}$ increases and its minimum is
achieved at the point $\Phi=0$, then it follows that on that
interval is greater than $1$. Therefore for each $\Phi \in [0,
\pi/2-\varphi_1]$, $\beta\mapsto \vn(\beta, \Phi)$ is an increasing
function. By elementary consideration for each $\beta \in (0,1/2)$
the function $\Phi\mapsto \vn(\beta, \Phi)$ is an increasing
function. Therefore
\[
\ve_2:=\left[\dfrac{\sqrt6}{2}\left(\dfrac{\pi}{2}-\arctan\dfrac{1}{\sqrt
2}\right)\right]^{2\beta}
\leq\dfrac{\sqrt6}{2}\left(\dfrac{\pi}{2}-\arctan\dfrac{1}{\sqrt
2}\right).
\]
Now, the function $(0,1/2)\ni\beta\mapsto \vj_1(\beta)$ is a
negative and increasing function. Therefore the supremum is at the
point $\beta=1/2$. Thus by elementary calculations, it follows that
\[
v(\pi/2)\leq \vv_1:={\bf -0.5727<0.\/}
\]
Therefore $v$ intersects $\varphi=\pi/2$ at a negative value.
%%%%%%%%%%%%%%%%%%%%%5%%%%%%%%%%%%%%%%%%%%%%%%%%%%%%%%%%%%%%%%%%%%%%%%%%%%%%%%%%%%%%%%%%%%%%%%%%%%%%%%%%%%%%%%%
\item[(1b)] Lower bound for $v(\pi/2)$ in the case $u>0$ and $\varphi \in (\varphi_1,
\pi/2)$.\\
Next we calculate a lower estimate for $v(\pi/2)$. For this we first
observe that
\[
\dfrac{dv}{d\varphi}=\dfrac12
\sqrt{\dfrac{2}{2^{2\beta}}\left[\dfrac{1}{\cos^{2\beta}\varphi}+
\dfrac{2^{1+\beta}}{(1+\sin^2\varphi)^\beta}\right]-v^2}\geq
\dfrac12\sqrt{\dfrac{2}{2^{2\beta}}\left[\dfrac{1}{\cos^{2\beta}\varphi}+\ve_3\right]-v^2},
\]
where $\ve_3:= \min_{\varphi \in [\varphi_1\,
\pi/2]}\vm(\beta,\varphi)=2$. Now we consider the interval
$J:=[7/16\pi, \pi/2]\subset[\varphi_1, \pi/2]$ and by a direct
calculation:
\begin{equation*}
\begin{aligned}
\dfrac{dv}{d\varphi}&\geq
\dfrac12\sqrt{\dfrac{2}{2^{2\beta}}\left[\dfrac{1}{\cos^{2\beta}(\varphi)}+2\right]-
v^2(\varphi_1)}\ \Rightarrow\\
&v(\pi/2) \geq v(7/16\pi) +
\dfrac12\int_{7/16\pi}^{\pi/2}\sqrt{\dfrac{2}{2^{2\beta}}\left[\dfrac{1}{\cos^{2\beta}\varphi}+2\right]-
6\,3^\beta \,8^{-\beta}}\ d\varphi\geq\\
& \geq v(\varphi_1) +
\int_{0}^{\pi/16}\dfrac{1}{2^{\beta+1/2}}\sqrt{\dfrac{1}{\Phi^{2\beta}}+\ve_4}\
d\Phi=  -\sqrt{6\, 3^\beta \, 8^{-\beta}} +
\int_{0}^{\pi/16}\dfrac{1}{2^{\beta+1/2}}\sqrt{\dfrac{1}{\Phi^{2\beta}}+\ve_4}\
d\Phi,
\end{aligned}
\end{equation*}
for $\ve_4:=2-v^2(\varphi_1)2^{2\beta}/2=2-3^{\beta+1}2^{-\beta}$;
moreover in the last inequality we used the fact that $\sin \Phi
\leq \Phi$. The function
\[
\vj_2(\beta):=
\int_{0}^{\pi/16}\dfrac{1}{2^{\beta+1/2}}\sqrt{\dfrac{1}{\Phi^{2\beta}}+\ve_4}
\]
is positive, increasing on the interval $[0, 1/2]$; moreover
$\vj_2(0)=0$. Thus
\[
v(\pi/2) \geq -\sqrt{6\, 3^\beta\, 2^{-3 \beta}}  \geq {\bf
-1.4993:=\vv_2\/}.
\]
 From (1a) \&(1b) it follows that $\vv_2 \leq
v(\varphi)\leq \vv_1$, for $\varphi \in [\varphi_1, \pi/2]$ and
$u>0$.
%%%%%%%%%%%%%%%%%%%%%%%%%%%%%%%%%%%%%%%%%%%%%%%%%%%%%%%%%%%%%%%%%%%%%%%%%%%%%%%%%%%
\item[(1c)] The continuation of the right branch through double
collision at $\varphi=\pi/2$ for $u<0$. We affirm that $v$ remains
negative in the interval $J$.

In fact, for $u<0$
\[
-\dfrac{dv}{d\varphi} \leq
\dfrac12\sqrt{\dfrac{2}{2^{2\beta}}\left[\dfrac{1}{\cos^{2\beta}\varphi}+\ve_5\right]}
\]
for $\ve_5:=\vm(\beta, 7\pi/16)$. Let $\vd_6:=\vn(\beta,\pi/16)$.
Integrating we get:
\begin{equation}
\begin{aligned}
v(7/16\pi)&\leq v(\pi/2)
+\dfrac{1}{2^{\beta+1/2}}\int_{7\pi/16}^{\pi/2}\sqrt{\dfrac{1}{\cos^{2\beta}\varphi}+\ve_5}\ d\varphi \leq\\
&\leq
v(\pi/2)+\dfrac{1}{2^{\beta+1/2}}\int_0^{\pi/16}\sqrt{\dfrac{\ve_6}{\Phi^{2\beta}}+\ve_5}\
d\Phi \leq \vv_1
+\dfrac{1}{2^{\beta1/2}}\int_0^{\pi/16}\sqrt{\dfrac{\ve_6}{\Phi}+\ve_5}\
d\Phi \\
&\leq -0.5727+
\dfrac{1}{2^{\beta1/2}}\int_0^{\pi/16}\sqrt{\dfrac{\ve_6}{\Phi}+\ve_5}\
d\Phi \leq {\bf-0.1004}<0
\end{aligned}
\end{equation}
Therefore
\[
\vv_2 \leq v(\pi/2)\leq v(\theta)\leq v(7\pi/16)<0, \qquad
\textrm{for}\ \ \theta \in [7\pi/16, \pi/2], \ \ \textrm{and}\ \
u<0.
\]
Therefore
\[
\vv_2 \leq v(\pi/2)\leq v(\theta)\leq v(7\pi/16)<0, \qquad
\textrm{for}\ \ \theta \in [7\pi/16, \pi/2], \ \ \textrm{and}\ \
u<0.
\]

%%%%%%%%%%%%%%%%%%%%%%%%%%%%%%%%%%%%%%%%%%%%%%%%%%%%%%%%%%%%%%%%%%%%%%%%%%%%%%%%%%%
\item[(1d)] From the above inequality we get the estimates $\vv_2^2\geq
v(\theta)^2$, for all $\theta \in [7\pi/16, \pi/2]$ and $u<0$. Using
this estimate  we can write:
\[
\begin{aligned}
-\dfrac{dv}{d\varphi}&=\dfrac12
\sqrt{2U(\varphi)-v^2}=\dfrac12\sqrt{\dfrac{2}{2^{2\beta}}\left[\dfrac{1}{\cos^{2\beta}\varphi}+
\dfrac{2^{1+\beta}}{(1+\sin^2\varphi)^\beta}\right] -v^2}\\
&\geq
\dfrac12\sqrt{\dfrac{2}{2^{2\beta}}\left[\dfrac{1}{\cos^{2\beta}\varphi}+\ve_7\right]
-v^2}\geq
\dfrac12\sqrt{\dfrac{2}{2^{2\beta}}\left[\dfrac{1}{\cos^{2\beta}\varphi}+\ve_7\right]
-\vv_2^2}
\end{aligned}
\]
where $\ve_7:=\min_{\varphi \in [7\pi/16, \pi/2]}\vm(\beta,
\varphi)=\vm(\beta, \pi/2)=2$. By a direct integration we have:
\begin{equation*}
\begin{aligned}
v(7\pi/16) &\geq v(\pi/2) +
\dfrac12\int_{7\pi/16}^{\pi/2}\sqrt{\dfrac{2}{2^{2\beta}}\left[\dfrac{1}{\cos^{2\beta}\varphi}+2\right]
-\vv_2^2}\ d\varphi \geq\\
&\geq v(\pi/2) + \dfrac12\int_{7\pi/16}^{\pi/2}
\sqrt{\dfrac{2}{2^{2\beta}}\left[\dfrac{1}{\cos^{2\beta}\varphi}+2\right]-\vv_2^2}\
d\varphi\geq\\
&\geq \vv_2+ \int_{7\pi/16}^{\pi/2}\dfrac{1}{2^{\beta+1/2}}
\sqrt{\dfrac{1}{\cos^{2\beta}\varphi}+\ve_8}\ d\varphi\\
&\geq \vv_2+ \int_{0}^{\pi/32}\dfrac{1}{2^{\beta+1/2}}
\sqrt{\dfrac{1}{\sin^{2\beta}\Phi}+\ve_8}\ d\Phi\\
& \geq \vv_2+ \int_{0}^{\pi/16}\dfrac{1}{2^{\beta+1/2}}
\sqrt{\dfrac{1}{\Phi^{2\beta}}+\ve_8}\ d\Phi \geq\vv_3 {\bf-1. 2078}
\end{aligned}
\end{equation*}
where $\ve_8:=2-\vv_2^2\,2^{2\beta-1}$.
\item[(1e)]To finish the estimate for the continuation of the right
branch through double collision, we now consider $\varphi \in
[0,7\pi/16]$ and $u<0$. We observe that in $\varphi=\varphi_1$,
$U(\theta)$ has a relative minimum and therefore for $\varphi \in
[0,7\pi/16]$ we have $U(\varphi)\geq U(\varphi_1)$. Now let $\bar
v:=\sqrt{2U(\varphi_1)}$; thus
\[
-\dfrac{dv}{d\varphi}\geq \dfrac12\sqrt{\bar v^2-v^2}
\]
and by a direct integration, we have
\[
-\arcsin\left(\dfrac{v}{\bar v}\right)\Big\vert_{v(0)}^{v(15\pi/32)}
\geq \dfrac12 \int_{0}^{7\pi/16} d \theta=7\pi/32.
\]
Therefore
\[
v(0)\geq \bar v\sin\left(7\pi/32+\arcsin(v(7\pi/16)/\bar
v)\right)\geq \bar v\sin\left(7\pi/32+\arcsin(\vv_2/\bar v)\right)
\geq {\bf0.4183}.
\]
\end{enumerate}
By symmetry, the same arguments holds also for the left branch of
$W^u(\vc_{-})$. Now the thesis easily follows by the fact that the
flow is of gradient-type with restpoints $\vc_{-}$ and $\vc_{+}$.
\begin{enumerate}
{\bf \item[(2)] The left branch of $W^u(\ve^1_{11})$ intersects
$\varphi=-\pi/2$ with at a negative value of $v$.\/}
\end{enumerate}
Now we consider the left branch of $W^u(\ve^-_{11})$. It is clear
that for the same value of $\varphi$, $v$ is greater along the left
branch of $W^u(\ve^-_{11})$ that  the corresponding value for the
left branch of $W^u(\ve_{12}^-)$, then due to the symmetries of the
problem it follows that the left branch of $W^u(\ve_{11}^-)$ meets
the section $v=0$ for some value of $\varphi$ which is negative. It
remains to prove that it actually meets $\varphi=-\pi/2$ at a
negative value of $v$. For $u<0$, we have:
\[
-\dfrac{dv}{d\varphi}=\sqrt{\va(\beta,\varphi)-\dfrac{v^2}{4}}
\]
where $\va(\beta,
\varphi):=\dfrac{U_\antiprisma(\varphi)}{2}=\dfrac{1}{2^{1+2\beta}}\left[\dfrac{1}{\cos^{2\beta}\varphi}+\dfrac{2^{1+\beta}}
{(1+\sin^2 \varphi)^\beta}\right]$. If $I_n:=[\tilde \varphi_{n-1},
\tilde \varphi_n]$ and $\va_n(\beta):=\max_{\varphi \in
I_n}\va(\beta,\varphi)$, then for $\varphi \in I_n$, we have:
\[
-\dfrac{dv}{d\varphi}\leq\sqrt{\va_n(\beta)-\dfrac{v^2}{4}}.
\]
By a direct integration, we can conclude that
\begin{equation}\label{eq:disricorsiva-2}
v(\tilde \varphi_n)\leq
2\sqrt{\va_n(\beta)}\sin\left(\dfrac12|\tilde \varphi_n-\tilde
\varphi_{n-1}|+\arcsin\left(\dfrac{v(\tilde
\varphi_{n-1})}{2\sqrt{\va_n(\beta)}}\right)\right).
\end{equation}
It is readily verified that $\va_n(\beta)=\dfrac{v(\tilde
\varphi_n)}{2}$ for any interval $I_n\subset [-\pi/2, -\varphi]\cup
[0, \varphi_1]:=\vA$, since $\varphi\mapsto \va(\beta,\varphi)$ is
positive and non-increasing on that subset; moreover
$\va_n(\beta)=\va(\beta, \tilde \varphi_{n-1})$ in $[\pi/2,
\pi/2]\setminus \vA:=\vB$ since on this set this function is
positive and non-decreasing. We define
\[
\varphi_n= \dfrac{(12-n)\varphi_1}{12}, \qquad \textrm{for each}\ \
n\, \in\, \n{36}.
\]
Solving recursively the inequality, we obtain
\[
v(\varphi_{32})=v(-5\varphi_1/3)\leq{\bf -1.6699.\/}
\]
On the other side, for $ \varphi\in [-\pi/2, 5\varphi_1/3]$,
\[
-\dfrac{dv}{d\varphi} \leq
\dfrac{1}{2^{\beta+1/2}}\sqrt{\dfrac{1}{\cos^{2\beta}\varphi}+\vd_9}
\]
where  $\vd_9:=\vm (\beta, -5\varphi_1/3)$. Again integrating (by
setting $\Phi:=\pi/2+\varphi$) and letting $\vd_{10}:= \vn(\beta,
\pi/2-5\varphi_1/3)$, we obtain:
\begin{equation}\label{eq:quasifine-2}
\begin{aligned}
-v(5\varphi_1/3)+ v(-\pi/2)&\leq
\dfrac{1}{2^{\beta+1/2}}\int_{-\pi/2}^{-5\varphi_1/3}
\sqrt{\dfrac{1}{\cos^{2\beta}\varphi}+\vd_9}\,d\varphi\\
&\leq \dfrac{1}{2^{\beta+1/2}}\int_0^{\pi/2-5\varphi_1/3}
\sqrt{\dfrac{1}{\sin^{2\beta}\Phi}+\vd_9}\,d\Phi\\
&\leq \dfrac{1}{2^{\beta+1/2}}\int_0^{\pi/2-5\varphi_1/3} \sqrt{\dfrac{\vd_{10}}{\Phi^{2\beta}}+\vd_9}\,d\Phi\\
&\leq {\bf 0.8721}.
\end{aligned}
\end{equation}
Now suppose that $v(-\pi/2)>0$; then from the inequality
\eqref{eq:quasifine-2} it  follows:
\[
v(-5\varphi_1/3)\geq {\bf -0.8721}
\]
which  is  a contradiction. We conclude that $v(-\pi/2)<0$ for the
left branch of $W^u(\ve_{11}^-)$.

\finedim

%%%%%======================================================================================
\subsection{Proof of Proposition \ref{prop:flussotetrab}}

The proof of this result will be divided into several steps.
\begin{enumerate}
\item {\bf There exists $\alpha_0 \in (1,2)$ such that for each $\alpha
\in (1, \alpha_0)$ the right branch of $W^u(\vc_{-})$ intersects
$\theta=\pi/2$ ($\theta=-\pi/2$, resp.) with $v<0$ and then
intersects the section $v=0$ with an angle $\theta \in (0,\pi/4)$.}
\item {\bf There exists $\alpha_1 \in (1,2)$ such that for each $\alpha
\in (\alpha_1, 2)$ the right branch of $W^u(\vc_{-})$ intersects
$\theta=\pi/2$ ($\theta=-\pi/2$, resp.) with $v<0$ and then
intersects the section $v=0$ with an angle $\theta \in (\pi/4,
\pi/2)$.}
\item {\bf There exists $\alpha_* \in (\alpha_0,\alpha_1)$ such that
the right branch of $W^u(\vc_{-})$ intersects the right branch of
$W^u(\vc_{+})$.}
\end{enumerate}
In order to prove item 1 we shall proceed as in the proof of the
previous proposition.
\begin{enumerate}
\item[(1a)] {\bf There exists $\alpha_0 \in (1,2)$ such that for
each $\alpha \in (1, \alpha_0)$, $v(\pi/2)<0$.\/}
\end{enumerate}
\begin{itemize}
\item  Upper bound for $v(\pi/2)$ in the case $u>0$ and $\theta \in (\pi/4,
\pi/2)$.
\begin{equation*}
\begin{aligned}
\dfrac{dv}{d\theta}&=\dfrac12
\sqrt{2^{1-\alpha}\left[1+\dfrac{1}{\cos^\alpha\theta}+\dfrac{1}{\sin^\alpha\theta}\right]-v^2}\leq
\dfrac{1}{2^{1/2+\beta}}\sqrt{\left[1+\dfrac{1}{\cos^\alpha\theta}+\dfrac{1}{\sin^\alpha\theta}\right]}\leq\\
&\leq \dfrac{1}{2^{1/2+\beta}}\sqrt{\dfrac{1}{\cos^\alpha
\theta}+\vb_2}
\end{aligned}
\end{equation*}
where
\[
\vb_2:=\max_{\theta \in [\pi/4, \pi/2]}h(\theta)=1+2^\beta, \qquad
h(\theta):= 1 + \dfrac{1}{\sin^\alpha \theta}.
\]
Thus
\begin{equation*}
\begin{aligned}
v(\pi/2)&\leq v(\pi/4)+ \int_{\pi/4}^{\pi/2}\dfrac{1}{2^{1/2+\beta}}
\sqrt{\dfrac{1}{\cos^\alpha\theta}+\vb_2}\ d\theta  =_{\Phi:=\pi/2-\theta}\\
&=v(\pi/4) + \dfrac{1}{2^{1/2+\beta}}\int_0^{\pi/4}
\sqrt{\dfrac{1}{\sin^\alpha\Phi}+\vb_2}\ d\Phi
\leq\\
&\leq v(\pi/4) +\dfrac{1}{2^{1/2+\beta}}\int_0^{\pi/4}
\sqrt{\dfrac{\vb_1}{\Phi^\alpha}+\vb_2} \ d\Phi=\\
&= -\dfrac{1}{2^{\beta-1/2}}\sqrt{1+2^{\beta+1}}
+\dfrac{1}{2^{1/2+\beta}}\int_0^{\pi/4}
\sqrt{\dfrac{\vb_1}{\Phi^{2\beta}}+\vb_2} \ d\Phi:=\vk_1(\beta),
\end{aligned}
\end{equation*}
where $\vb_1:=\max_{[0,\pi/4]}\left(\dfrac{\Phi}{\sin
\Phi}\right)^\alpha=\left(\dfrac{\pi}{4}\sqrt 2\right)^\alpha$ since
$g(\Phi):=\left(\dfrac{\Phi}{\sin\Phi}\right)^\alpha$ is an
increasing function. We observe that  $(1/2,1)\ni\beta\mapsto
\vk_1(\beta)$ is a continuous  increasing function on the interval
$(1/2, 1)$. However it changes sign in the interval $A_1:=(1.4,
1.6)$;
%%$A:=(0.73068, 0.73029)$;
in particular if $\alpha_0$ is any number in the open interval $(1,
1.46136]$, then we can conclude that
\[
v(\pi/2)\leq \vk_1(1/2)\approx {\bf-0.1659}:=\vw_1.
\]

\item Lower bound for $v(\pi/2)$ in the case $u>0$ and $\theta \in (\pi/4,
\pi/2)$. For this we first observe that
\[
\dfrac12
\sqrt{2^{1-\alpha}\left[1+\dfrac{1}{\cos^\alpha\theta}+\dfrac{1}{\sin^\alpha\theta}\right]-v^2}\geq
\dfrac{1}{2^{1/2+\beta}}\sqrt{\dfrac{1}{\cos^\alpha\theta}+\vb_3-2^{\alpha-1}v^2},
\]
where $\vb_3:= \min_{\theta \in [\pi/4, \pi/2]}h(\theta)=2$. Now we
consider the interval $J:=[\pi/4, 3\pi/8]$ and by a direct
calculation, we have:
\begin{equation*}
\begin{aligned}
\dfrac{dv}{d\theta}&\geq
\dfrac{1}{2^{1/2+\beta}}\sqrt{\dfrac{1}{\cos^\alpha\theta}+2-
v^2(\pi/4)}\ \Rightarrow\\ v(\pi/2) &\geq v(3\pi/8) +
\int_{3\pi/8}^{\pi/2}\dfrac{1}{2^{1/2+\beta}}\sqrt{\dfrac{1}{\cos^\alpha\theta}+2-
v^2(\pi/4)}\ d\theta\geq\\
&\geq v(3\pi/8) +
\int_{0}^{\pi/8}\dfrac{1}{2^{1/2+\beta}}\sqrt{\dfrac{1}{\Phi^{2\beta}}+1-2^{\beta+1}}\
d\Phi:=\vj_2(\beta).
\end{aligned}
\end{equation*}
Since:
\[
v(3\pi/8):=v_\beta(3\pi/8)=-\sqrt{2^{1-2\beta}\left[1+
\left(\dfrac{2\sqrt2}{1+\sqrt2}\right)^\beta+\left(\dfrac{2\sqrt2}{\sqrt2-1}\right)^\beta
\right]};
\]
and since both the functions
\[
(0,1/2)\ni\beta\mapsto v_\beta(3\pi/8)\ \textrm{and}\ (0, 1/2)\ni
\beta\mapsto
\int_{0}^{\pi/8}\dfrac{1}{2^{1/2+\beta}}\sqrt{\dfrac{1}{\Phi^{2\beta}}+1-2^{\beta+1}}
\]
are strictly increasing it immediately follows that also $\vj_2$ is
an increasing function and therefore
%\[
%v_{1/2}(3\pi/8)\approx \textcolor{red}{-2.449}\leq v(3\pi/8)\leq
%v_0(3\pi/8) \approx \textcolor{red}{-2.167}.
%\]
%To conclude the estimate, we observe that also the function:
%\[
%(0,1/2)\ni\beta\mapsto
%\vi(\beta):=\int_{0}^{\pi/8}\dfrac{1}{2^{1/2+\beta}}\sqrt{\dfrac{1}{\theta^{2\beta}}+1-
%2^{(\beta+1)}}\ d\theta\in \R^+
%\]
%is a positive and increasing function, with $\vi(0)=0$. Therefore,
we can  conclude that
\[
v(\pi/2)\geq \vj_1(1/2):=\vw_2\approx{\bf-1.6267}.
\]
We observe that for $\alpha \geq 1.1$ we have $ v(\pi/2)\geq
\vj_1(0.55):=\vw_3\approx{\bf-1.5285}$
Therefore $v$ intersects $\theta=\pi/2$ at a negative value.
%%%%%%%%%%%%%%%%%%%%%5%%%%%%%%%%%%%%%%%%%%%%%%%%%%%%%%%%%%%%%%%%%%%%%%%%%%%%%%%%%%%%%%%%%%%%%%%%%%%%%%%%%%%%%%%

\end{itemize}
%%%%%%%%%%%%%%%%%%%%%%%%%%%%%%%%%%%%%%%%%%%%%%%%%%%%%%%%%%%%%%%%%%%%%%%%%%%%%%%%%%%
\begin{enumerate}
\item[(1b)] {\bf There exists $\alpha_{00} \in (1, \alpha_0)$ such that the continuation of the right branch through double
collision at $\theta=\pi/2$ for $u<0$, remains negative in the
interval $[\pi/4, \pi/2]$.\/}
\end{enumerate}
In fact, for $u<0$
\[
-\dfrac{dv}{d\theta} \leq
\dfrac12\sqrt{2^{1-\alpha}\left[\dfrac{1}{\cos^\alpha\theta}+\vb_5\right]}
\]
for $\vb_5:=h(49\pi/100)$ and $h$ as above. Let $\vb_6:=g(\pi/100)$.
Integrating we get:
\begin{equation}
\begin{aligned}
v(49\pi/100)&\leq v(\pi/2)
+\dfrac12\int_{49\pi/100}^{\pi/2}\dfrac{1}{2^{1/2+\beta}}\sqrt{\left[\dfrac{1}{\cos^\alpha\theta}+\vb_5\right]}\ d\theta \leq\\
&\leq
v(\pi/2)+\dfrac12\int_0^{\pi/100}\dfrac{1}{2^{1/2+\beta}}\sqrt{\left[\dfrac{\vb_6}{\Phi^{2\beta}}+\vb_5\right]}\
d\Phi \leq\\
\leq& \vw_1
+\dfrac12\int_0^{\pi/100}\dfrac{1}{2^{1/2+\beta}}\sqrt{\left[\dfrac{\vb_6}{\Phi^{2\beta}}+\vb_5\right]}\
d\Phi\leq\\
&\leq -0.1659+
\dfrac12\int_0^{\pi/100}\dfrac{1}{2^{1/2+\beta}}\sqrt{\left[\dfrac{\vb_6}{\Phi^{2\beta}}+\vb_5\right]}\
d\Phi\\
%v_1
%+\dfrac{1}{\sqrt2}\int_0^{\pi/16}\sqrt{\dfrac{\va_6}{\theta}+\va_5}\
%d\theta \\
&\approx {\bf-0.0057},
\end{aligned}
\end{equation}
which holds for $\alpha \in (1,1.1]$. Thus the thesis follows by
choosing $\alpha_{00}=1.1$. Therefore
\[
\vw_2\leq v(\pi/2)\leq v(\theta)\leq v(49\pi/100)<0, \qquad
\textrm{for}\ \ \theta \in [49\pi/100, \pi/2], \ \ \textrm{and}\ \
u<0.
\]
From the above inequality we get the estimates $\vw_2^2\geq
v(\theta)^2$, for all $\theta \in [49\pi/100, \pi/2]$ and $u<0$.
Using this estimate  we can write:
\[
-\dfrac{dv}{d\theta}=\dfrac12
\sqrt{2^{1-\alpha}\left[1+\dfrac{1}{\sin^\alpha\theta}+\dfrac{1}{\cos^\alpha\theta}\right]-v^2}
\geq\dfrac12\sqrt{2^{1-\alpha}\left[\dfrac{1}{\cos\theta}+\va_7\right]
-\vw_2^2}
\]
where $\va_7:=\min_{\theta \in [49\pi/100,
\pi/2]}h(\theta)=h(\pi/2)=2$. By a direct integration we have:
\begin{equation*}
\begin{aligned}
v(49\pi/100) &\geq v(\pi/2) +
\dfrac12\int_{49\pi/100}^{\pi/2}\sqrt{2^{1-\alpha}\left[\dfrac{1}{\cos^\alpha\theta}+\va_7\right]
-\vw_2^2}\ d\theta \geq\\
&\geq v(\pi/2) + \dfrac12 \int_0^{\pi/100}
\sqrt{2^{1-\alpha}\left[\dfrac{1}{\Phi^\alpha}+2\right]-\vw_2^2}\
d\Phi\geq\\ &\geq \vw_2 + \dfrac12 \int_0^{\pi/100}
\sqrt{2^{1-\alpha}\left[\dfrac{1}{\Phi^\alpha}+2\right]-\vw_2^2}\
d\Phi\geq \\
&\geq%\vw_2 + \dfrac12 \int_0^{\pi/100}
%\sqrt{\left[\dfrac{1}{\Phi}+2\right]-\vw_2^2}\ d\Phi \approx
{\bf-1.1804}.
\end{aligned}
\end{equation*}
\begin{itemize}
\item Now we prove that the continuation of the right
branch through double collision is negative at the point
$\theta=\pi/4$.
\end{itemize}
In fact, for $u<0$
\[
-\dfrac{dv}{d\theta} \leq
\dfrac12\sqrt{2^{1-\alpha}\left[\dfrac{1}{\cos^\alpha\theta}+\vb_5\right]}
\]
for $\vb_5:=h(\pi/4)$ and $h$ as above. Let $\vb_6:=g(\pi/4)$.
Integrating we get:
\begin{equation}
\begin{aligned}
v(\pi/4)&\leq v(\pi/2)
+\dfrac12\int_{\pi/4}^{\pi/2}\dfrac{1}{2^{1/2+\beta}}\sqrt{\left[\dfrac{1}{\cos^\alpha\theta}+\vb_5\right]}\ d\theta \leq\\
&\leq
v(\pi/2)+\dfrac12\int_0^{\pi/4}\dfrac{1}{2^{1/2+\beta}}\sqrt{\left[\dfrac{\vb_6}{\Phi^{2\beta}}+\vb_5\right]}\
d\Phi \leq\\
\leq& \vw_1
+\dfrac12\int_0^{\pi/}\dfrac{1}{2^{1/2+\beta}}\sqrt{\left[\dfrac{\vb_6}{\Phi^{2\beta}}+\vb_5\right]}\
d\Phi\leq\\
&\leq \vk_1(\beta)+
\dfrac12\int_0^{\pi/4}\dfrac{1}{2^{1/2+\beta}}\sqrt{\left[\dfrac{\vb_6}{\Phi^{2\beta}}+\vb_5\right]}\
d\Phi\leq\\&\leq -0.6857 +
\dfrac12\int_0^{\pi/4}\dfrac{1}{2^{1/2+\beta}}\sqrt{\left[\dfrac{\vb_6}{\Phi^{2\beta}}+\vb_5\right]}\
d\Phi\leq\\
%v_1
%+\dfrac{1}{\sqrt2}\int_0^{\pi/16}\sqrt{\dfrac{\va_6}{\theta}+\va_5}\
%d\theta \\
&\leq {\bf-0.0713}.
\end{aligned}
\end{equation}

%
%
%For this we consider $\theta \in [\pi/4,49\pi/100]$ and $u<0$. We
%observe that in $\theta=49\pi/100$, $U(\theta)$ has a relative
%maximum and therefore for $\theta \in [\pi/4,49\pi/100 ]$ we have
%$U(\theta)\leq U(49\pi/100)$. Now let $\bar
%v:=-\sqrt{2U(49\pi/100)}$; thus
%\[
%-\dfrac{dv}{d\theta}\leq \dfrac12\sqrt{\bar v^2-v^2}
%\]
%and by a direct integration, we have
%\[
%-\arcsin\left(\dfrac{v}{\bar
%v}\right)\Big\vert^{v(49\pi/100)}_{v(\pi/4)} \leq \dfrac12
%\int_{\pi/4}^{49\pi/100} d \theta=3\pi/25.
%\]
%Therefore
%\[
%v(\pi/4)\leq \bar v\sin(3\pi/25+\arcsin(v(49\pi/100)/\bar v))\leq
%\textcolor{blue}{-0.4502}.
%\]
%Moreover  in $\theta=\pi/4$, $U(\theta)$ has a minimum and therefore
%for $\theta \in [\pi/8,3\pi/8 ]$ we have $U(\theta)\geq U(\pi/4)$.
%Now let $\bar v':=\sqrt{2U(\pi/4)}$; thus
%\[
%-\dfrac{dv}{d\theta}\geq \dfrac12\sqrt{\bar v'^2-v^2}
%\]
%and by a direct integration, we have
%\[
%-\arcsin\left(\dfrac{v}{\bar
%v}\right)\Big\vert^{v(3\pi/8)}_{v(\pi/8)} \geq \dfrac12
%\int_{\pi/8}^{3\pi/8} d \theta=\pi/8.
%\]
%Therefore
%\[
%v(\pi/8)\geq \bar v'\sin(\pi/8+\arcsin(v(3\pi/8)/\bar v'))\geq {\bf
%-0.3741}.
%\]
%These two estimates in particular imply that $v(\pi/4)<0$ by taking
%into account that for $u<0$ the solution $v$ is strictly decreasing
%on the interval $[\pi/8, 3\pi/8]$.

To finish the estimate for the continuation of the right branch
through double collision, we now consider $\theta \in [0, \pi/8]$
and $u<0$. We observe that in $\theta=\pi/8$, $U(\theta)$ has a
relative minimum and therefore for $\theta \in [0, \pi/8]$ we have
$U(\theta)\geq U(\pi/8)$. Now let $\tilde v:=\sqrt{2U(\pi/8)}$; thus
\[
-\dfrac{dv}{d\theta}\geq \dfrac12\sqrt{\tilde v^2-v^2}
\]
and by a direct integration, we have
\[
-\arcsin\left(\dfrac{v}{\tilde v}\right)\Big\vert^{v(\pi/8)}_{v(0)}
\geq \dfrac12 \int_{0}^{\pi/8} d \theta=\pi/16.
\]
Therefore
\[
v(0)\geq \tilde v\sin\left(\pi/16+
\arcsin\left(\dfrac{v(\pi/8)}{\tilde v}\right)\right)\geq \tilde
v\sin\left(\pi/16+ \arcsin\left(\dfrac{-0.3741}{\tilde
v}\right)\right) \approx {\bf0.4948}.
\]
By symmetry, the same arguments holds also for the left branch of
$W^u(\vc_{-})$. Now the thesis easily follows by the fact that the
flow is of gradient-type with restpoints $\vc_{-}$ and $\vc_{+}$.
\begin{enumerate}
\item[(1c)] {\bf There exists $\alpha_{01} \in (\alpha_{00}, \alpha_0)$
such that the continuation of the right branch through double
collision at $\theta=\pi/2$ for $u<0$, changes sign in the interval
$[\pi/4, \pi/2]$ for any $\alpha \in [\alpha_{01}, \alpha_0)$.\/}
\end{enumerate}
In fact, using this estimate  we can write:
\[
-\dfrac{dv}{d\theta}=\dfrac12
\sqrt{2^{1-\alpha}\left[1+\dfrac{1}{\sin^\alpha\theta}+\dfrac{1}{\cos^\alpha\theta}\right]-v^2}
\geq\dfrac12\sqrt{2^{1-\alpha}\left[\dfrac{1}{\cos\theta}+\va_7\right]
-\vw_2^2}
\]
where $\va_7:=\min_{\theta \in [49\pi/100,
\pi/2]}h(\theta)=h(\pi/2)=2$. By a direct integration we have:
\begin{equation*}
\begin{aligned}
v(49\pi/100) &\geq v(\pi/2) +
\dfrac12\int_{49\pi/100}^{\pi/2}\sqrt{2^{1-\alpha}\left[\dfrac{1}{\cos^\alpha\theta}+\va_7\right]
-\vw_2^2}\ d\theta \geq\\
&\geq v(\pi/2) + \dfrac12 \int_0^{\pi/100}
\sqrt{2^{1-\alpha}\left[\dfrac{1}{\Phi^\alpha}+2\right]-\vw_2^2}\
d\Phi\geq\\ &\geq \vw_2 + \dfrac12 \int_0^{\pi/100}
\sqrt{2^{1-\alpha}\left[\dfrac{1}{\Phi^\alpha}+2\right]-\vw_2^2}\
d\Phi\geq \\
&\geq%\vw_2 + \dfrac12 \int_0^{\pi/100}
%\sqrt{\left[\dfrac{1}{\Phi}+2\right]-\vw_2^2}\ d\Phi \approx
{\bf-1.1804}.
\end{aligned}
\end{equation*}
To finish the estimate for the continuation of the right branch
through double collision, we now consider $\theta \in [3\pi/8,
\pi/2]$ and $u<0$. We observe that in $\theta=\pi/4$, $U(\theta)$
has a relative minimum and therefore for $\theta \in [\pi/4,
49\pi/100]$ we have $U(\theta)\geq U(\pi/4)$. Now let $\tilde
v:=\sqrt{2U(\pi/4)}$; thus
\[
-\dfrac{dv}{d\theta}\geq \dfrac12\sqrt{\tilde v^2-v^2}
\]
and by a direct integration, we have
\[
-\arcsin\left(\dfrac{v}{\tilde
v}\right)\Big\vert^{v(49\pi/100)}_{v(\pi/4)} \geq \dfrac12
\int_{\pi/4}^{49\pi/100} d \theta=3\pi/25.
\]
Therefore
\[
v(\pi/4)\geq \tilde v\sin\left(3\pi/25+
\arcsin\left(\dfrac{v(3\pi/25)}{\tilde v}\right)\right)\geq \tilde
v\sin\left(3\pi/25+ \arcsin\left(\dfrac{-1.1804}{\tilde
v}\right)\right) \approx {\bf0.0511}.
\]

By symmetry, the same arguments holds also for the left branch of
$W^u(\vc_{-})$. Now the thesis easily follows by the fact that the
flow is of gradient-type with restpoints $\vc_{-}$ and $\vc_{+}$.
\begin{enumerate}
\item[(1d)] {\bf There exists $\alpha_{0*} \in (\alpha_{00}, \alpha_{01})$
such that for the continuation of the right branch through double
collision at $\theta=\pi/2$ for $u<0$, we have $v(\pi/4)=0$.\/}
\end{enumerate}
In order to prove this fact we consider the one parameter family of
initial value problems:
\begin{equation}
(P_\alpha) \qquad \left\{
\begin{array}{ll}
\dfrac{dv}{d\theta}= \dfrac12 \sqrt{2^{1-\alpha}\left[1+
\dfrac{1}{\cos^\alpha \theta}+ \dfrac{1}{\sin^\alpha
\theta}\right]-v^2}:=\vF_\alpha(\theta, v)\\
\\
v(\pi/4)=-\sqrt{2^{1-\alpha}(1+2^{\beta+1})},
\end{array}\right.
\end{equation}
which can be written in integral form as follows
\[
v_\alpha(\theta)=v(\pi/4)+\int_{\pi/4}^\theta \vF_\alpha(\bar
\theta, v(\bar\theta))\, d\bar\theta, \qquad \forall\, \theta \in
(\pi/4, \pi/2).
\]
If $v$ is a classical ($\mathscr C^1$)-solution of the ivp above,
then it is a $\mathscr C^0$-solution of the corresponding Volterra
integral equation; moreover the function
\[
\vI:[1,\alpha_0)\longrightarrow  \R: \alpha \longmapsto
\vI(\alpha):=\int_{\pi/4}^{\pi/2} \vF_\alpha(\theta, v(\theta))\,
d\theta,
\]
is continuous, by taking into account the theorem on integrals
depending on parameters. Thus the function
\[
\vV:(1,\alpha_0)\rightarrow \R:\alpha\mapsto
\vV(\alpha):=v_\alpha(\pi/2):= v(\pi/4)+ \vI(\alpha)
\]
is continuous. Moreover by the monotonicity of the integral of a
nonnegative function with respect to the domain of integration and
by the fact that the map $\alpha\mapsto \vF_\alpha(\theta, v)$
increases on $(1,2)$, it follows that the map $\vV$ also increases.
Furthermore $\vV(\alpha_{00})<0$ and $\vV(\alpha_{01})>0$; thus by
taking into account the theorem of zeros of a continuous and
increasing function, there exists a unique $\alpha_*$ in between
such that $V(\alpha_{0*})=0$, which means nothing but that
$v_\alpha(\pi/4)=0$.

\smallskip

In order to prove item 2, we proceed as follows.
\begin{enumerate}
\item[(2)] {\bf There exists $\alpha_1 \in (1,2)$ such that for
each $\alpha \in (\alpha_1, 2)$, $v(\pi/2)>0$.\/}
\end{enumerate}
\begin{itemize}
\item Lower bound for $v(\pi/2)$ in the case $u>0$ and $\theta \in (\pi/4,
\pi/2)$.
\end{itemize}

Next we calculate a lower estimate for $v(\pi/2)$. For this we first
observe that
\[
\dfrac12
\sqrt{2^{1-\alpha}\left[1+\dfrac{1}{\cos^\alpha\theta}+\dfrac{1}{\sin^\alpha\theta}\right]-v^2}\geq
\dfrac{1}{2}\sqrt{2^{1-\alpha}\left[\dfrac{1}{\cos^\alpha\theta}+\vb_3\right]-v^2},
\]
where $\vb_3:= \min_{\theta \in [\pi/4, \pi/2]}h(\theta)=2$. Now we
consider the interval $J:=[3\pi/8, \pi/2]$ and by a direct
calculation, we have:
\begin{equation*}
\begin{aligned}
\dfrac{dv}{d\theta}&\geq
\dfrac{1}{2}\sqrt{2^{1-\alpha}\left[\dfrac{1}{\cos^\alpha\theta}+2\right]-v^2(\pi/4)}\
\Rightarrow\\ v(\pi/2) &\geq v(3\pi/8) +
\int_{3\pi/8}^{\pi/2}\dfrac{1}{2^{1/2+\beta}}\sqrt{\dfrac{1}{\cos^\alpha\theta}+2-
2^{\alpha-1}v^2(\pi/4)}\ d\theta\geq\\
&\geq v(3\pi/8) +
\int_{0}^{\pi/8}\dfrac{1}{2^{1/2+\beta}}\sqrt{\dfrac{1}{\Phi^{2\beta}}+1-2^{\beta+1}}\
d\Phi:=\vk_2(\beta).
\end{aligned}
\end{equation*}
We observe that:
\[
v(3\pi/8):=v_\beta(3\pi/8)=-\sqrt{2^{1-2\beta}\left[1+
\left(\dfrac{2\sqrt2}{1+\sqrt2}\right)^\beta+\left(\dfrac{2\sqrt2}{\sqrt2-1}\right)^\beta
\right]}.
\]
Moreover, the function below
\[
(1/2,1)\ni\beta\mapsto v_\beta(3\pi/8)\ \textrm{and}\ (0, 1/2)\ni
\beta\mapsto
\int_{0}^{\pi/8}\dfrac{1}{2^{1/2+\beta}}\sqrt{\dfrac{1}{\Phi^{2\beta}}+1-2^{\beta+1}}
\]
as well as $\vk_2$ are strictly increasing functions. %it immediately
%follows that also $\vk_2$ is an increasing function and therefore
%\[
%v_{1/2}(3\pi/8)\approx \textcolor{red}{-2.449}\leq v(3\pi/8)\leq
%v_0(3\pi/8) \approx \textcolor{red}{-2.167}.
%\]
%To conclude the estimate, we observe that also the function:
%\[
%(0,1/2)\ni\beta\mapsto
%\vi(\beta):=\int_{0}^{\pi/8}\dfrac{1}{2^{1/2+\beta}}\sqrt{\dfrac{1}{\theta^{2\beta}}+1-
%2^{(\beta+1)}}\ d\theta\in \R^+
%\]
%is a positive and increasing function, with $\vi(0)=0$. Therefore,
However it changes sign and in particular it is positive in the
interval  $A_2:=[1.7,2)$. Therefore by choosing any number $\alpha_1
\in A_2$ we can  conclude that
\[
v(\pi/2)\geq \vk_2(1.7):=\vw_2\approx {\bf 0.1055}>0.
\]

By symmetry, the same arguments holds also for the left branch of
$W^u(\vc_{-})$. Now the thesis easily follows by the fact that the
flow is of gradient-type with restpoints $\vc_{-}$ and $\vc_{+}$.

\smallskip
\begin{enumerate}
\item[(3)] {\bf There exists a unique $\alpha_* \in (\alpha_0,\alpha_1)$ such that
$v_\alpha(\pi/2)=0$.\/}
\end{enumerate}
In order to prove this fact we consider the one parameter family of
initial value problems:
\begin{equation}
(P_\alpha) \qquad \left\{
\begin{array}{ll}
\dfrac{dv}{d\theta}= \dfrac12 \sqrt{2^{1-\alpha}\left[1+
\dfrac{1}{\cos^\alpha \theta}+ \dfrac{1}{\sin^\alpha
\theta}\right]-v^2}:=\vF_\alpha(\theta, v)\\
\\
v(\pi/4)=-\sqrt{2^{1-\alpha}(1+2^{\beta+1})},
\end{array}\right.
\end{equation}
which can be written in integral form as follows
\[
v_\alpha(\theta)=v(\pi/4)+\int_{\pi/4}^\theta \vF_\alpha(\bar
\theta, v(\bar\theta))\, d\bar\theta, \qquad \forall\, \theta \in
(\pi/4, \pi/2).
\]
If $v$ is a classical ($\mathscr C^1$)-solution of the ivp above,
then it is a $\mathscr C^0$-solution of the corresponding Volterra
integral equation; moreover the function
\[
\vI:[1,2)\longrightarrow  \R: \alpha \longmapsto
\vI(\alpha):=\int_{\pi/4}^{\pi/2} \vF_\alpha(\theta, v(\theta))\,
d\theta,
\]
is continuous, by taking into account the theorem on integrals
depending on parameters. Thus the function
\[
\vV:(1,2)\rightarrow \R:\alpha\mapsto \vV(\alpha):=v_\alpha(\pi/2):=
v(\pi/4)+ \vI(\alpha)
\]
is continuous. Moreover by the monotonicity of the integral of a
nonnegative function with respect to the domain of integration and
by the fact that the map $\alpha\mapsto \vF_\alpha(\theta, v)$
increases on $(1,2)$, it follows that the map $\vV$ also increases.
Furthermore $\vV(\alpha_0)<0$ and $\vV(\alpha_1)>0$; thus by taking
into account the theorem of zeros of a continuous and increasing
function, there exists a unique $\alpha_*$ in between such that
$V(\alpha_*)=0$, which means nothing but that $v_\alpha(\pi/2)=0$.

\end{document}